\documentclass[12pt]{amsart}
\usepackage{amsfonts,amssymb,amscd,amsmath,epsfig}
\usepackage{latexsym}
\usepackage{mathrsfs}
\usepackage{tocvsec2}
\usepackage[%dvipdfm,
hyperindex=true,pagebackref=true,bookmarks=true,
colorlinks=true,linkcolor=blue,citecolor=red]{hyperref}
\def\~{{\rm --}}
\begin{document}

%MACOS FOR LECTURES ON DAHA
\renewcommand{\tilde}{\widetilde}
\renewcommand{\hat}{\widehat}

\newcommand{\BR}{{\mathbb R}}
\newcommand{\BQ}{{\mathbb Q}}
\newcommand{\BC}{{\mathbb C}}
\newcommand{\BP}{{\mathbb P}}
\newcommand{\BZ}{{\mathbb Z}}
\newcommand{\BN}{{\mathbb N}}
\newcommand{\BS}{{\mathbb S}}

\newcommand{\cH}{{\mathcal H}}
\newcommand{\cA}{{\mathcal A}}
\newcommand{\cB}{{\mathcal B}}
\newcommand{\ccF}{{\mathfrak F}}
\newcommand{\cD}{{\mathcal D}}
\newcommand{\cL}{{\mathcal L}}
\newcommand{\cF}{{\mathcal F}}
\newcommand{\cP}{{\mathcal P}}
\newcommand{\cX}{{\mathcal X}}
\newcommand{\cY}{{\mathcal Y}}
\newcommand{\cS}{{\mathcal S}}
\newcommand{\cSol}{\hbox{$\mathcal Sol$}}
\newcommand{\cT}{\hbox{$\mathcal T$}}

\newcommand{\Z}{{\mathbb Z}}
\newcommand{\Q}{{\mathbb Q}}
\newcommand{\N}{{\mathbb N}}
\newcommand{\C}{{\mathbb C}}
\newcommand{\R}{{\mathbb R}}
\newcommand{\X}{{\mathbb X}}
\newcommand{\Y}{{\mathbb Y}}

\newcommand{\CH}{{\mathcal H}}
\newcommand{\CA}{{\mathcal A}}

\def\HH{\mbox{${\mathcal H}$\kern-5.2pt${\mathcal H}$}}

\newcommand{\binomial}[2]{\genfrac{(}{)}{0pt}{}{ #1 }{ #2 }}
\newcommand{\qbinomial}[2]{\genfrac{[}{]}{0pt}{}{ #1 }{ #2 }_q }
\newcommand{\qbinom}[3]{\genfrac{[}{]}{0pt}{}{ #1 }{ #2 }_{ #3 } }

%%SPECIAL SEC 1.0

\def\der{\partial}
\def\tensor{\otimes}
\def\gam{\gamma} \def\Gam{\Gamma}
\def\del{\delta} \def\Del{\Delta}
\def\kap{\kappa}
\def\lam{\lambda} \def\Lam{\Lambda}
\def\Comp{{\mathbb C}}
\def\sM{{\mathcal M}}

\newtheorem{theorem}{Theorem}[section]
\newtheorem{maintheorem}[theorem]{Main Theorem}
\newtheorem{proposition}[theorem]{Proposition}
\newtheorem{definition}[theorem]{Definition}
\newtheorem{lemma}[theorem]{Lemma}
\newtheorem{corollary}[theorem]{Corollary}
\newtheorem{notation}[theorem]{Notation}
\newtheorem{remark}[theorem]{Remark}
\newtheorem{example}[theorem]{Example}

\newtheorem{theorem }{Theorem}[section]
\newtheorem{maintheorem }[theorem]{Main Theorem}
\newtheorem{proposition }[theorem]{Proposition}
\newtheorem{definition }[theorem]{Definition}
\newtheorem{lemma }[theorem]{Lemma}
\newtheorem{corollary }[theorem]{Corollary}
\newtheorem{notation }[theorem]{Notation}
\newtheorem{remark }[theorem]{Remark}
\newtheorem{example }[theorem]{Example}

\newtheorem{ maintheorem }[theorem]{Main Theorem}
\newtheorem{ theorem}{Theorem}[section]
\newtheorem{ proposition}[theorem]{Proposition}
\newtheorem{ definition}[theorem]{Definition}
\newtheorem{ lemma}[theorem]{Lemma}
\newtheorem{ corollary}[theorem]{Corollary}
\newtheorem{ notation}[theorem]{Notation}
\newtheorem{ remark}[theorem]{Remark}
\newtheorem{ example}[theorem]{Example}

\newtheorem{thm}{Theorem}[section]
\newtheorem{prop}[thm]{Proposition}
\newtheorem{lem}[thm]{Lemma}
\newtheorem{cor}[thm]{Corollary}
\newtheorem{conj}[thm]{Conjecture}
\newtheorem{con}[thm]{Conjecture}
\newtheorem{dfn}[thm]{Definition}
\newtheorem{df}[thm]{Definition}
 \newcommand{\rem}{{\bf Comment.\ }}
 \newcommand{\rmk}{{\bf Comment.\ }}
 \newcommand{\exmp}{{\bf Example.\ }}
 \newcommand{\ex}{{\bf Example.\ }}
 \newcommand{\prob}{{\bf Problem.\ }}

\newtheorem{note}{Note} 
\renewcommand{\thenote}{}
\newtheorem*{acka}{Acknowledgments}
\newtheorem{ack}{Acknowledgments}
\renewcommand{\theack}{}
\renewcommand{\appendixname}{\bf Appendix}
\renewcommand{\proof}{{\em Proof.\ }}

\hyphenation{
ap-pen-dix as-ymp-tot-ic at-trib-uted at-trib-ut-able
Bry-li-n-sky com-mu-ta-tion de-ge-ne-rate
de-riv-a-tive dis-trib-ute equi-vari-ant ex-tra-or-di-nary  
geo-met-ric griev-ance griev-ous grad-ed ho-lo-no-my ho-mo-thetic
in-fin-ite-ly in-fin-i-tes-i-mal Ha-rish Cha-n-dra mul-ti-plic-able 
non-euclid-ean non-iso-mor-phic non-smooth par-a-digm 
par-a-bol-ic pa-rab-o-loid pa-ram-e-trize phe-nom-e-non 
post-script pseu-do-dif-fer-en-tial pseu-do-fi-nite 
qua-drat-ics quad-ra-ture Han-kel rec-tan-gle semi-def-i-nite 
set-up wide-spread Euler-ian Feb-ru-ary Gauss-ian Grothen-dieck 
Hamil-ton-ian Her-mi-t-ian her-mi-t-ian Jan-u-ary 
Japan-ese Ka-shi-wa-ra Kor-te-weg Le-gendre No-vem-ber Rie-mann-ian 
Sep-tem-ber Za-mo-lo-d-chi-kov Kni-zh-nik quan-tum Op-dam
Mac-do-nald Ca-lo-ge-ro Su-ther-land Mo-ser 
Ol-sha-net-sky  Pe-re-lo-mov in-de-pen-dent ope-ra-tors 
cy-clo-to-mic ra-tio-nal de-gen-er-a-tion 
in-ter-est-ing de-for-ma-tions de-for-ma-tion pro-ce-dure 
fol-lows ope-ra-tors  pre-serve suf-fices ap-proach 
for-mu-las con-sider its com-ple-tion cor-re-spond-ing 
au-to-mor-phism be-cause pro-por-tional fi-nal-ly let-ting 
equi-v-a-lence ge-n-er-al-ized Mac-do-nald iden-ti-ties 
cor-re-s-pond sub-dia-grams par-ti-tion na-t-u-ral-ly 
or-dered stan-dard de-for-ma-tion ar-gu-ment com-bined 
sphe-r-i-cal rep-re-sen-ta-tions tri-go-no-me-t-ric
ge-n-er-al-ly speak-ing pri-m-it-ive ir-re-du-cible 
sum-ma-tion  rep-re-sen-ta-tives pro-por-ti-o-na-li-ty
ultra-sphe-ri-cal Ro-gers}

\def\ffor{\quad\hbox{ for }\quad}
\def\wwhen{\quad\hbox{ when }\quad}
\def\wwhere{\quad\hbox{ where }\quad}
\def\aand{\quad\hbox{ and }\quad}
\def\for{\  \hbox{ for } \ }
\def\iif{ \ \hbox{ if } \ }
\def\when{ \ \hbox{ when } \ }
\def\where{\  \hbox{ where } \ }
\def\and{\  \hbox{ and } \ }
\def\and{\  \hbox{ and } \ }
\def\oor{\  \hbox{ or } \ }
\def\proof{{\em Proof. \  }}

\def\equal{\stackrel{\,\mathbf{def}}{= \kern-3pt =}}

\def\la{\lambda}
\def\La{\Lambda}
\def\om{\omega}
\def\Om{\Omega}
\def\Th{\Theta}
\def\th{\theta}
\def\al{\alpha}
\def\be{\beta}
\def\ga{\gamma}
\def\ep{\epsilon}
\def\up{\upsilon}
\def\Up{\Upsilon}
\def\de{\delta}
\def\De{\Delta}
\def\ka{\kappa}
\def\kapp{\hbox{\bf \ae}}
\def\si{\sigma}
\def\Si{\Sigma}
\def\Ga{\Gamma}
\def\ze{\zeta}
\def\io{\iota}
\def\bio{b^\iota}
\def\aio{a^\iota}
\def\twio{\tilde{w}^\iota}
\def\hwio{\hat{w}^\iota}
\def\gio{\g^\iota}
\def\Bio{B^\iota}

\def\del{\delta}
\def\pa{\partial}
\def\vp{\varphi}
\def\ve{\varepsilon}
\def\inf{\infty}

\def\vph{\varphi}
\def\vps{\varpsi}
\def\vPh{\varPhi}
\def\vep{\varepsilon}
\def\vpi{{\varpi}}
\def\vth{{\vartheta}}
\def\vsi{{\varsigma}}
\def\vrh{{\varrho}}

\def\bph{\bar{\phi}}
\def\bsi{\bar{\si}}
\def\bvp{\bar{\varphi}}

\newcommand{\bS}{{\mathbf S}}
\newcommand{\bH}{{\mathbf H}}
\newcommand{\bF}{{\mathbf F}}
\newcommand{\bE}{{\mathbf E}}

\def\tal{\tilde{\alpha}}
\def\tbe{\tilde{\beta}}
\def\tde{\tilde{\delta}}
\def\tpi{\tilde{\pi}}
\def\txi{\tilde{\xi}}
\def\tPi{\tilde{\Pi}}
\def\tPhi{\tilde{\Phi}}
\def\tV{\tilde{V}}
\def\tJ{\tilde{J}}
\def\tla{\tilde{\lambda}}
\def\tga{\tilde{\gamma}}
\def\tGa{\tilde{\Gamma}}
\def\tvs{\tilde{{\varsigma}}}
\def\tu{\tilde{u}}
\def\tU{\tilde{U}}
\def\tw{\widetilde w}
\def\tW{\widetilde W}
\def\tB{\tilde B}
\def\tv{\tilde v}
\def\tV{\tilde V}
\def\tz{\tilde z}
\def\tb{\tilde b}
\def\ta{\tilde a}
\def\tih{\tilde h}
\def\trh{\tilde {\rho}}
\def\tx{\tilde x}
\def\tf{\tilde f}
\def\tg{\tilde g}
\def\tG{\tilde G}
\def\tk{\tilde k}
\def\tl{\tilde l}
\def\tL{\tilde L}
\def\tD{\tilde D}
\def\tR{\tilde R}
\def\tP{\tilde P}
\def\tH{\tilde H}
\def\tp{\tilde p}

\def\hH{\hat{H}}
\def\hh{\hat{h}}
\def\hR{\hat{R}}
\def\hY{\hat{Y}}
\def\hX{\hat{X}}
\def\hP{\hat{P}}
\def\hT{\hat{T}}
\def\hV{\hat{V}}
\def\hG{\hat{G}}
\def\hF{\hat{F}}
\def\hw{\widehat{w}}
\def\hW{\widehat{W}}
\def\hu{\hat{u}}
\def\hs{\hat{s}}
\def\hv{\hat{v}}
\def\hb{\hat{b}}
\def\hB{\widehat{B}}
\def\hze{\hat{\zeta}}
\def\hsi{\hat{\sigma}}
\def\hrh{\hat{\rho}}
\def\hth{\hat{\theta}}
\def\hy{\hat{y}}
\def\hx{\hat{x}}
\def\hz{\hat{z}}
\def\hg{\hat{g}}
\def\he{\hat{e}}
\def\hE{\widehat{E}}

\def\B{\mathbf{B}}
\def\I{\mathbf{I}}
\def\P{\mathbf{P}}
\def\G{\mathbf{G}}
\def\S{\mathbf{S}}
\def\F{\mathbf{F}}
\def\one{\mathbf{1}}
\def\Sn{\mathbf{S}_n}
\def\0{\mathbf{0}}
\def\H{\mathbf{H}}
\def\V{\mathbf{V}}

\def\f{\mathcal{F}}
\def\çF{\mathcal{F}}
\def\o{\mathcal{O}}
\def\t{\mathcal{T}}
\def\r{\mathcal{R}}
\def\l{\mathcal{L}}
\def\m{\mathcal{M}}
\def\k{\mathcal{K}}
\def\n{\mathcal{N}}
\def\d{\mathcal{D}}
\def\p{\mathcal{P}}
\def\cP{\mathcal{P}}
\def\a{\mathcal{A}}
\def\h{\mathcal{H}}
\def\c{\mathcal{C}}
\def\y{\mathcal{Y}}
\def\e{\mathcal{E}}
\def\v{\mathcal{V}}
\def\z{\mathcal{Z}}
\def\x{\mathcal{X}}
\def\s{\mathcal{S}}
\def\g{\mathcal{G}}
\def\u{\mathcal{U}}
\def\w{\mathcal{W}}
\def\i{\mathcal{I}}
\def\j{\mathcal{J}}
\def\b{\mathcal{B}}

\def\lan{\langle}
\def\llb{(\!(}
\def\ran{\rangle}
\def\rrb{)\!)}
 \def\dim{{\hbox{\rm dim}}_{\mathbb C}\,}
\def\lng{\hbox{\rm{\tiny lng}}}
\def\sht{\hbox{\rm{\tiny sht}}}
\def\sph{\hbox{\rm{\tiny sph}}}
\def\inv{\hbox{\rm{\tiny inv}}}

\def\br#1{\langle #1 \rangle}

\def\rank{\hbox{rank}}
\def\gl{\mathfrak{gl}_N}
%\def\sgn{\hbox{sgn}}
%\font\germ=eufb10 %at 12pt 
%\def\mathfrak#1{\hbox{\germ #1}}

\newcommand{\Aut}{\operatorname{Aut}}
\newcommand{\Hom}{\operatorname{Hom}}
\newcommand{\End}{\operatorname{End}}
\newcommand{\Ind}{\operatorname{Ind}}
\newcommand{\ad}{\operatorname{ad}}
\newcommand{\pr}{\operatorname{pr}}
\newcommand{\aweyl}{\tilde{\mathbb S}_n}
\newcommand{\hec}{{\mathcal H}^t_n}
\newcommand{\Func}{{\mathcal F}({\mathbb C}^n,{\mathcal H}^t_n)}
\newcommand{\tr}{\operatorname{tr}}
\newcommand{\Out}{\operatorname{Out}}
\newcommand{\Rad}{\operatorname{Rad}}
\newcommand{\Spec}{\operatorname{Spec}}
\newcommand{\id}{\operatorname{id}}
\newcommand{\Int}{\operatorname{Int}}
\newcommand{\ct} {\operatorname{ct}}

\newcommand{\rat}{{\mathbb Q}}
\newcommand{\real}{{\mathbb R}}
\newcommand{\cplx}{{\mathbb C}}
\newcommand{\zint}{{\mathbb Z}}

\newcommand{\sq}{\phantom{1}\hfill$\qed$}
\newcommand{\Rea}{\Re}
\newcommand{\Ima}{\Im}

\newcommand{\st}{\bowtie}
\newcommand{\modd}{\mbox{\,mod\,}}
\newcommand{\lr}{\langle}
\newcommand{\rr}{\rangle}
\newcommand{\eps}{\varepsilon}
\newcommand{\phk}{\phi^{(k)}}
\newcommand{\psk}{\psi^{(k)}}
\newcommand{\Res}{\mbox{Res}\;}
\newcommand{\sgn}{\mbox{sgn}}
\newcommand{\mn} {\left\{ \begin{array}{c}m\\
n\end{array}\right\}}

\def\sX{\mathscr{X}}
\def\sH{\mathscr{H}}
\def\sY{\mathscr{Y}}
\def\TT{\mathfrak{T}}
\def\JJ{\mathfrak{J}}
\def\HH{\mathfrak{H}}
\def\FF{\mathfrak{F}}
\def\GG{\mathfrak{G}}
\def\CC{\mathfrak{C}}
\def\LL{\mathfrak{L}}

\def\BB{\mathfrak{B}}
\def\AA{\mathfrak{A}}
\def\ZZ{\mathfrak{Z}}
\def\HH{\hbox{${\mathcal H}$\kern-5.2pt${\mathcal H}$}}
\def\HHH{\hbox{${\mathbb H}$\kern-4.2pt${\mathbb H}$}}
\def\tHH{\widetilde{\HH\ }}

\font\smm=msbm10 at 12pt 
\def\symbol#1{\hbox{\smm #1}}
\def\lsmash{{\symbol n}}
\def\rsmash{{\symbol o}}
\def\#{\sharp}

\font\tenbf=cmbx10
\font\tenrm=cmr10
\font\tenit=cmti10
\font\ninebf=cmbx9
\font\ninerm=cmr9
\font\nineit=cmti9
\font\eightbf=cmbx8
\font\eightrm=cmr8
\font\eightit=cmti8
\font\sevenrm=cmr7
\font\sevenbf=cmbx7

%END MACROS

\title [DAHA and iterated torus knots]
{DAHA and iterated torus knots}
\author[Ivan Cherednik]{Ivan Cherednik $^\dag$}
\author[Ivan Danilenko] {Ivan Danilenko}
%\date{February 2, 2014}

\begin{abstract}
The theory of DAHA-Jones polynomials is extended
from torus knots to their arbitrary iterations
(for any reduced root systems and weights), which 
incudes the polynomiality, duality and other properties
of the DAHA superpolynomials. Presumably they 
coincide with the reduced stable Khovanov-Rozansky 
polynomials in the case of non-negative coefficients.
The new theory matches well the classical theory of 
algebraic knots and (unibranch) plane curve singularities; 
the Puiseux expansion naturally emerges. We expect the
DAHA superpolynomials for such knots to coincide with 
the reduced ones in the Oblomkov-Shende-Rasmussen 
Conjecture upon its further generalization  
to arbitrary dominant weights. For instance, the DAHA 
uncolored superpolynomials at $a=0, q=1\,$  are conjectured 
to provide the Betti numbers of the Jacobian factors of 
the corresponding singularities.
\end{abstract}

\thanks{$^\dag$ \today.
\ \ \ Partially supported by NSF grant
DMS--1363138 and the Simons Foundation}

\address[I. Cherednik]{Department of Mathematics, UNC
Chapel Hill, North Carolina 27599, USA\\
chered@email.unc.edu}

\address[I. Danilenko]{MIPT, 9 Institutskiy per., 
Dolgoprudny, Moscow Region 141700, Russia \\
ITEP, 25 Bolshaya Cheremushkinskaya, Moscow 117218, Russia \\ 
danilenko\_i@mail.ru}

 \def\sht{\raisebox{0.4ex}{\hbox{\rm{\tiny sht}}}}
 \def\bysame{{\bf --- }}
 \def\~{{\bf --}}
 \def\rr{{\mathsf r}}
 \def\ss{{\mathsf s}}
 \def\mm{{\mathsf m}}
 \def\pp{{\mathsf p}}
 \def\ll{{\mathsf l}}
 \def\aa{{\mathsf a}}
 \def\bb{{\mathsf b}}
 \def\NS{\hbox{\tiny\sf ns}}
 \def\ssum{\hbox{\small$\sum$}}
\newcommand{\comment}[1]{}
\renewcommand{\tilde}{\widetilde}
\renewcommand{\hat}{\widehat}
\renewcommand{\V}{\mathbb{V}}
\renewcommand{\F}{\mathbb{F}}
\newcommand{\dagx}{\hbox{\tiny\mathversion{bold}$\dag$}}
\newcommand{\ddagx}{\hbox{\tiny\mathversion{bold}$\ddag$}}
\newtheorem{conjecture}[theorem]{Conjecture}
\newcommand*\toeq{
\raisebox{-0.15 em}{\,\ensuremath{
\xrightarrow{\raisebox{-0.3 em}{\ensuremath{\sim}}}}\,}
}
\newcommand{\unknot}{\hbox{\tiny\!\raisebox{0.2 em}{$\bigcirc$}}}

\vskip -0.0cm
%\par
%{\centering
%\medskip
%\par}
%\vskip -0.0cm
\maketitle
\smallskip

{\tiny
\centerline{{\bf MSC} (2010): 14H50, 17B22, 17B45, 20C08,
20F36, 33D52, 30F10, 55N10, 57M25, 57R58}
}
\smallskip

\vskip -0.0cm
\renewcommand{\baselinestretch}{1.2}
{\textmd
\tableofcontents
}
\renewcommand{\baselinestretch}{1.0}
%}
%$\mathfrak{a,b,c,d,e,f,g,h,i,j,k,l,m,n,o,p,q,r,s,t,e^u,e^v,e^w}$
\vfill\eject

\renewcommand{\natural}{\wr}

\setcounter{section}{-1}
\setcounter{equation}{0}
\section{\sc Introduction}
The theory of {\em DAHA-Jones polynomials\,} of torus
knots \cite{CJ,GoN,CJJ}
is fully extended in this paper to arbitrary iterated knots 
for any reduced root systems and dominant weights, 
which incudes the polynomiality, duality and other 
properties of the {\em DAHA superpolynomials}. We conjecture
the coincidence of the latter with the stable reduced
{\em Khovanov-Rozansky 
polynomials\,} \cite{KhR1,KhR2,Ras} for 
{\em pseudo-algebraic\,} knots, for instance for all 
algebraic knots.

A similar connection is also expected with  
the physics superpolynomials associated with the 
{\em BPS states\,} (see e.g. \cite{DGR,AS,DMMSS}) and 
those related to the {\em rational DAHA\,} 
(see \cite{GORS,GoN} and references therein). 
We note that using rational DAHA for this is restricted so
far only to the torus knots and the fundamental 
representation (a generalization to its symmetric/wedge
powers is in progress). 
\smallskip
 
The new theory matches well the classical theory of algebraic 
knots and unibranch {\em plane curve singalarities\,},
though we see no {\em a priori} reasons for this.  
The {\em Newton pairs\,} in
the theory of {\em Puiseux expansion\,} naturally emerge in
our approach; see e.g. \cite{EN}.  
For instance, the iterations that are trivial topologically 
result in interesting algebraic symmetries of the DAHA
superpolynomials, for instance in uniform formulas 
for the torus knots in the form $T(\rr,m\rr+\ss)$
for {\em any\,} $m>0$. 
%The number of operations they require is roughly comparable 
%with that needed for $T(\rr,\ss)$ for {\em any\,} $m>0$. 

We expect that our superpolynomials add $\,t\,$ to
the {\em Oblomkov-Shende Conjecture\,} \cite{ObS,Ma}
and add weights to Conjecture 2 from \cite{ORS} for 
algebraic knots; see \cite{Ma}, Section 1.5.
The {\em reduced ORS Conjecture\,}   
connects the Poincar\'e polynomial 
$\mathscr{P}_{\hbox{\tiny KhR}}$
of the triply-graded {\em reduced\,} Khovanov-Rozansky 
homology of an algebraic link 
with $\mathscr{P}_{\hbox{\tiny alg}}$ describing
the cohomology of the nested Hilbert scheme of 
the corresponding plane curve singularity under the 
{\em weight filtration\,}. 
We conjecture them to coincide with the uncolored DAHA 
superpolynomial $\mathscr{P}_{\hbox{\tiny \!DAHA}}$. 

The coincidence with $\mathscr{P}_{\hbox{\tiny KhR}}$ is expected
for any {\em pseudo-algebraic\,} knots (those with 
positive coefficients of DAHA superpolynomials). The reduction
of $\mathscr{P}_{\hbox{\tiny \!DAHA}}$ to the {\em HOMFLYPT}
polynomials is conjectured for any iterated knots and weights
(we prove it for $sl_2$ and checked in all the examples we
calculated). The link $\mathscr{P}_{\hbox{\tiny \!DAHA}}\sim 
\mathscr{P}_{\hbox{\tiny KhR}}$ for arbitrary  
iterated knots is a challenge (we provide examples).

%See \cite{GORS} and Conjecture \ref{CONCONJ} below
%concerning the connection of the DAHA $q,t,a$\~parameters
%and the ``standard" ones.

Finding $\mathscr{P}_{\hbox{\tiny KhR}}$
is an involved task; one of the reasons is that
the {\em skein relations\,} generally cannot be extended 
from the HOMFLYPT polynomials to the HOMFLYPT homology. 
The polynomial $\mathscr{P}_{\hbox{\tiny alg}}$ 
is very sophisticated too. For instance,
its portion of minimal $a$\~degree requires knowing the 
so-called {\em perverse filtration\,} in the cohomology of 
the {\em Jacobian factors\,}. Even the Betti numbers of
the latter are known only in few cases.  See  
Theorem 22 and Conjecture 23 in \cite{Pi} concerning the torus 
knots and those for {\em Puiseux exponents\,} $(4,2u,v)$.  

By contrast,  $\mathscr{P}_{\hbox{\tiny \!DAHA}}$ can be
calculated in an entirely formal way without 
any (theoretical) limitations; we expect that this construction
will eventually include arbitrary iterated {\em links\,}. 
Thus the main obstacle with stating Parts $(ii,iii)$ of Conjecture 
\ref{CONCONJ} below for arbitrary weights is the absence of 
$\mathscr{P}_{\hbox{\tiny KhR}}$ and
$\mathscr{P}_{\hbox{\tiny alg}}$ in such a generality. 
Concerning $\mathscr{P}_{\hbox{\tiny \!DAHA}}\sim 
\mathscr{P}_{\hbox{\tiny alg}}$, 
the connection of DAHA with the $K$\~theory of the
{\em Hilbert schemes\,} of $\C^2$ \cite{SV} and the 
{\em affine Springer fibers\,} \cite{VV,Yun} must be mentioned;
the Jacobian factors of isolated plane curve 
singularities are certain affine Springer fibers \cite{La}.
\smallskip

In the examples we calculated, the $t$\~coefficients of the
DAHA superpolynomials evaluated at $a\!=\!0, q\!=\!1$ coincide 
with the Betti numbers of the corresponding Jacobian factors 
when these numbers are known/conjectured. Their sum, the 
{\em Euler number\,}, can be calculated via
the HOMFLYPT polynomial; see Section 7 of \cite{ObS} and \cite{Ma}.
 
We provide conjectural Betti numbers in the 
case of Puiseux exponents $(6,8,9)$ and $(6,9,10)$ and quite a few 
examples beyond the technique used in \cite{Pi}
In the case of  $(6,8,9)$, we calculate {\em formal\,} 
DAHA counterparts of such Betti numbers for $2\omega_1,\omega_2$. 
Also, we do this for several uncolored {\em non-algebraic\,} 
iterated knots; for instance for
``non-algebraic" $v=3,5$, preceding 
the formulas from \cite{Pi} 
for $(4,6,v\ge 7)$. 
\smallskip

In all considered examples (not only for the pseudo-algebraic
iterated knots), our superpolynomials 
reduce to the corresponding HOMFLYPT polynomials under $t\mapsto q,
a\mapsto -a$ for arbitrary dominant weights and any iterated knots. 
They match the (uncolored) Khovanov 
polynomials  for $sl_2$ in the sense of 
Conjecture 3.7 from \cite{CJ} for quite a few considered
{\em pseudo-algebraic\,} knots, when all coefficients of 
$\mathscr{P}_{\hbox{\tiny \!DAHA}}$ are positive. Generally,
this connection can be very involved.
\smallskip

We note that our paper was triggered by paper \cite{Sam},
though it is restricted to $A_1$ and we do not quite understand 
its approach.
%It was suggested in \cite{Sam} to combine applying the
%automorphisms from $GL_2(\Z)$ to Macdonald polynomials 
%as in \cite{CJ} with certain intermediate projections, 
%but the implementation of this idea there and its connection
%with our construction remain unclear to us. 
The Newton pairs and other features of our construction do not 
appear in \cite{Sam}. Also, the polynomials $J_n$ there have
significant $q,t$\~denominators (even in 
the uncolored case); the polynomiality of DAHA-Jones polynomials 
is the key in our theory. The $J_2$\~polynomials for
$C\!ab(\pm 5,2)T(3,2)$ from Section 5.2 are very different from
our ones in these cases (and we do not understand how they
were obtained); see (\ref{jon5-2}) below. The first part of 
\cite{Sam} is devoted  to the {\em skein modules\,}, which 
continues \cite{BeS}; this direction is of obvious importance.   
\medskip

\setcounter{section}{0}
\setcounter{equation}{0}
\section{\sc DAHA and Macdonald polynomials}
\subsection{\bf Affine root systems}
Let $R=\{\al\}   \subset \R^n$ be a root system of type
$A_n,\ldots,\!G_2$
with respect to a euclidean form $(z,z')$ on $\R^n
\ni z,z'$,
$W$ the Weyl group 
generated by the reflections $s_\al$,
$R_{+}$ the set of positive  roots
corresponding to fixed simple 
roots $\al_1,...,\al_n;$ $R_-=-R_+$. 
The form is normalized
by the condition  $(\al,\al)=2$ for 
{\em short\,} roots. 
The weight lattice is
$P=\oplus^n_{i=1}\Z \om_i$, 
where $\{\om_i\}$ are fundamental weights:
$ (\om_i,\al_j^\vee)=\de_{ij}$ for the
coroots $\al^\vee=2\al/(\al,\al).$
Replacing $\Z$ by $\Z_{\pm}=\{m\in\Z, \pm m\ge 0\}$, we obtain
$P_\pm.$ See  e.g., \cite{Bo} or \cite{C101}. 

Setting 
$\nu_\al\equal (\al,\al)/2$,
the vectors $\ \tal=[\al,\nu_\al j] \in
\R^n\times \R \subset \R^{n+1}$
for $\al \in R, j \in \Z $ form the
{\em twisted affine root system\,}
$\tR \supset R$ ($z\in \R^n$ are identified with $ [z,0]$).
We add $\al_0 \equal [-\vth,1]$ to the simple
roots for the {\em maximal short root\,} $\vth\in R_+$.
The corresponding set
$\tR_+$ of positive roots is %%%% tR_+BOOK!!!
$R_+\cup \{[\al,\nu_\al j],\ \al\in R, \ j > 0\}$.

The set of the indices of the images of $\al_0$ by all
automorphisms of the affine Dynkin diagram will be denoted by 
$O$ ($O=\{0\} \for E_8,F_4,G_2$). 
Let $O'\equal\{r\in O, r\neq 0\}$.
The elements $\om_r$ for $r\in O'$ are  
{\em minuscule weights\,}, defined by the inequalities 
$(\om_r,\al^\vee)\le 1$ for all $\al \in R_+$. We set $\om_0=0$
for the sake of uniformity.
\smallskip

{\sf Affine Weyl groups.}
Given $\tal=[\al,\nu_\al j]\in \tR,  \ b \in P$, let
\begin{align}
&s_{\tal}(\tz)\ =\  \tz-(z,\al^\vee)\tal,\
\ b'(\tz)\ =\ [z,\ze-(z,b)]
\label{ondon}
\end{align}
for $\tz=[z,\ze] \in \R^{n+1}$.
The
{\em affine Weyl group\,} $\tW=\lan s_{\tal}, \tal\in \tR_+\ran)$ 
is the semidirect product $W\lsmash Q$ of
its subgroups $W=$ $\lan s_\al,
\al \in R_+\ran$ and $Q$, where $\al$ is identified with
\begin{align*}
& s_{\al}s_{[\al,\,\nu_{\al}]}=\
s_{[-\al,\,\nu_\al]}s_{\al}\for
\al\in R.
\end{align*}

The {\em extended Weyl group\,} $ \hW$ is $W\lsmash P$, where
the corresponding action is 
\begin{align}
&(wb)([z,\ze])\ =\ [w(z),\ze-(z,b)] \for w\in W, b\in P.
\label{ondthr}
\end{align}
It is isomorphic to $\tW\lsmash \Pi$ for $\Pi\equal P/Q$. 
The latter group consists of $\pi_0=$id\, and the images $\pi_r$
of minuscule $\om_r$ in $P/Q$. 
%Note that 
%$\pi_r^{-1}$ is $\pi_{r^*}$, where $^*$ is the standard involution 
%(sometimes trivial) of the {\em nonaffine\,} Dynkin diagram.

The group $\Pi$
is naturally identified with the subgroup of $\hW$ of the
elements of the length zero; the {\em length\, } is defined as 
follows:
\begin{align*}
&l(\hw)=|\la(\hw)| \for \la(\hw)\equal\tR_+\cap \hw^{-1}(-\tR_+).
\end{align*}
One has $\om_r=\pi_r u_r$ for $r\in O'$, where $u_r$ is the 
element $u\in W$ of minimal length such that $u(\om_r)\in P_-$. 
\smallskip

Setting $\hw = \pi_r\tw \in \hW$ for $\pi_r\in \Pi,\, \tw\in \tW,$
\,$l(\hw)$ coincides with the length of any reduced decomposition
of $\tw$ in terms of the simple reflections
$s_i, 0\le i\le n.$ 
\medskip

\subsection{\bf Definition of DAHA}
We follow \cite{CJJ,CJ,C101}.
Let $\mm,$ be the least natural number
such that  $(P,P)=(1/\mm)\Z.$  Thus
$\mm=|\Pi|$ unless 
$\mm=2 \for D_{2k}$ and $\ \mm=1 \for B_{2k},C_{k}.$ 

The double affine Hecke algebra, {\em DAHA\,}, depends
on the parameters
$q, t_\nu\, (\nu\in \{\nu_\al\})\,$ and will be defined
over the ring
$Z_{q,t}\equal\Z[q^{\pm 1/\mm},t_\nu^{\pm 1/2}]$
formed by
polynomials in terms of $q^{\pm 1/\mm}$ and
$\{t_\nu^{1/2}\}.$ Note that the coefficients of the 
Macdonald polynomials will belong to 
$
\Q(q,t_\nu).
$

For $\tal=[\al,\nu_\al j] \in \tR,\ 0\le i\le n$, we set
\begin{align*}
&   t_{\tal} =t_{\al}=t_{\nu_\al}=q_\al^{k_\nu} ,\ \, 
q_{\tal}=q^{\nu_\al}, \ \, t_i = t_{\al_i},  q_i=q_{\al_i},
\end{align*}

Also, using here (and below) {\em\small\, sht,\ lng\,} instead 
of $\nu$, we set
\begin{align*}
\rho_k\equal \frac{1}{2}\!\sum_{\al>0} k_\al \al=
k_{\sht}\rho_{\sht}\!+\!k_{\lng}\rho_{\lng},\ \,
\rho_\nu=\frac{1}{2}\!\sum_{\nu_\al=\nu} \al=
\!\!\sum_{\nu_i=\nu,\, i>0}  \om_i.
\end{align*}

For pairwise commutative $X_1,\ldots,X_n,$
\begin{align}
& X_{\tb}\ \equal\ \prod_{i=1}^nX_i^{l_i} q^{ j}
\iif \tb=[b,j],\ \hw(X_{\tb})\ =\ X_{\hw(\tb)},
\label{Xdex}\\
&\hbox{where\ } b=\sum_{i=1}^n l_i \om_i\in P,\ j \in
\frac{1}{ m}\Z,\ \hw\in \hW.
\notag \end{align}
For instance, $X_0\equal X_{\al_0}=qX_\vth^{-1}$.
\medskip

Recall that 
$\om_r=\pi_r u_r$ for $r\in O'$ (see above).
We will use that $\pi_r^{-1}$ is $\pi_{\iota(i)}$, where
$\iota$ is the standard involution  of the nonaffine 
Dynkin diagram,
induced by $\al_i\mapsto -w_0(\al_i)$. Generally,
$\iota(b)=-w_0(b)=b^\iota$, where $w_0$ is the longest element 
in $W$.  Finally, we set $m_{ij}=2,3,4,6$
when the number of links between $\al_i$ and $\al_j$ in the affine 
Dynkin diagram is $0,1,2,3$.

\begin{definition}
The double affine Hecke algebra $\HH\ $
is generated over $\Z_{q,t}$ by
the elements $\{ T_i,\ 0\le i\le n\}$,
pairwise commutative $\{X_b, \ b\in P\}$ satisfying
(\ref{Xdex})
and the group $\Pi,$ where the following relations are imposed:

(o)\ \  $ (T_i-t_i^{1/2})(T_i+t_i^{-1/2})\ =\
0,\ 0\ \le\ i\ \le\ n$;

(i)\ \ \ $ T_iT_jT_i...\ =\ T_jT_iT_j...,\ m_{ij}$
factors on each side;

(ii)\ \   $ \pi_rT_i\pi_r^{-1}\ =\ T_j \iif
\pi_r(\al_i)=\al_j$;

(iii)\  $T_iX_b \ =\ X_b X_{\al_i}^{-1} T_i^{-1} \iif
(b,\al^\vee_i)=1,\
0 \le i\le  n$;

(iv)\ $T_iX_b\ =\ X_b T_i\ $ if $\ (b,\al^\vee_i)=0
\for 0 \le i\le  n$;

(v)\ \ $\pi_rX_b \pi_r^{-1}\ =\ X_{\pi_r(b)}\ =\
X_{ u^{-1}_r(b)}
 q^{(\om_{\iota(r)},b)},\  r\in O'$.
\label{double}
\end{definition}

Given $\tw \in \tW, r\in O,\ $ the product
\begin{align}
&T_{\pi_r\tw}\equal \pi_r T_{i_l}\cdots T_{i_1},\where
\tw=s_{i_l}\cdots s_{i_1} \for l=l(\tw),
\label{Twx}
\end{align}
does not depend on the choice of the reduced decomposition
Moreover,
\begin{align}
&T_{\hv}T_{\hw}\ =\ T_{\hv\hw}\  \hbox{ whenever\,}\
 l(\hv\hw)=l(\hv)+l(\hw) \for
\hv,\hw \in \hW. \label{TTx}
\end{align}
In particular, we arrive at the pairwise
commutative elements 
\begin{align}
& Y_{b}\equal
\prod_{i=1}^nY_i^{l_i} \iif
b=\sum_{i=1}^n l_i\om_i\in P,\ 
Y_i\equal T_{\om_i},b\in P.
\label{Ybx}
\end{align}
When acting in $\v$, they are called {\em difference
Dunkl operators.}.
\smallskip

\subsection{\bf The automorphisms}\label{sect:Aut}
The following maps can be (uniquely) extended to
an automorphism of $\HH\,$, where 
$q^{1/(2\mm)}$ must be added to $\Z_{q,t}$
(see \cite{C101}, (3.2.10)-(3.2.15)):
\begin{align}\label{tauplus}
& \tau_+:\  X_b \mapsto X_b, \ T_i\mapsto T_i\, (i>0),\
\ Y_r \mapsto X_rY_r q^{-\frac{(\om_r,\om_r)}{2}}\,,
\\
& \tau_+:\ T_0\mapsto  q^{-1}\,X_\vth T_0^{-1},\
\pi_r \mapsto q^{-\frac{(\om_r,\om_r)}{2}}X_r\pi_r\
(r\in O'),\notag\\
& \label{taumin}
\tau_-:\ Y_b \mapsto \,Y_b, \ T_i\mapsto T_i\, (i\ge 0),\
\ X_r \mapsto Y_r X_r q^\frac{(\om_r,\om_r)}{ 2},\\
&\tau_-(X_{\vth})= 
q T_0 X_\vth^{-1} T_{s_{\vth}}^{-1};\ \
\si\equal \tau_+\tau_-^{-1}\tau_+\, =\,
\tau_-^{-1}\tau_+\tau_-^{-1},\notag\\
&\si(X_b)=Y_b^{-1},\   \si(Y_b)=
T_{w_0}^{-1}X_{b^\iota}^{-1}T_{w_0},\ \si(T_i)=T_i (i>0).
\label{taux}
\end{align}
These automorphisms fix $\ t_\nu,\ q$
and their fractional powers, as well as the
following {\em anti-involution\,}:
\begin{align}
&\vph:\ 
X_b\mapsto Y_b^{-1},\, Y_b\mapsto X_b^{-1},\
T_i\mapsto T_i\ (1\le i\le n).\label{starphi}
%\\
%&\phi(\tau_+)\equal \phi\circ \tau_+\circ \phi\ =\ \tau_-\,,\ 
%\ \, \phi(\tau_-)\ =\ \tau_+.\notag 
\end{align} 

We will also need the involution
\begin{align}
\eta:\ &T_i\mapsto T_i^{-1},\ X_b\mapsto X_b^{-1},\
\pi_r\mapsto \pi_r \,(0\le i\le n),
\label{etatxpi}
\end{align}
which ``conjugates" $t,q$; namely, 
$\,t^{\frac{1}{2\mm}}_\nu \mapsto t_\nu^{-\frac{1}{2\mm}},\, 
q^{\frac{1}{2\mm}}\mapsto  q^{-\frac{1}{2\mm}}$.
As above, $b\in P,\ r\in O'$.
The involution $\eta\,$ extends the
{\em Kazhdan\~Lusztig involution\,} in the affine Hecke theory;
see \cite{C101}, (3.2.19-22). Note that 
$$
\vph\tau_\pm\vph=\tau_\mp=\si\tau_\pm^{-1}\si^{-1},\ 
\eta\tau_{\pm}\eta=\tau_{\pm}^{-1}.
$$

Let us list the matrices corresponding to the automorphisms and
anti-automorphisms above upon the natural projection 
onto $GL_2(\Z)$, corresponding to  
$t^{\frac{1}{2\mm}}_\nu=1=q^{\frac{1}{2\mm}}$. 
The matrix {\tiny 
$\begin{pmatrix} \al & \be \\ \ga & \de\\ \end{pmatrix}$}
will represent the map $X_b\mapsto X_b^\al Y_b^\ga,
Y_b\mapsto X_b^\be Y_b^\de$ for $b\in P$. One has:
\smallskip

$\!\!\!\tau_+\mapsto$ 
{\tiny 
$\begin{pmatrix}1 & 1 \\0 & 1 \\ \end{pmatrix}$},\ 
$\tau_-\mapsto$ 
{\tiny 
$\begin{pmatrix}1 & 0 \\1 & 1 \\ \end{pmatrix}$},\
$\si\mapsto$ 
{\tiny 
$\begin{pmatrix}0 & 1 \\-1 & 0 \\ \end{pmatrix}$},\
$\vph\mapsto$ 
{\tiny 
$\begin{pmatrix}0 & -1 \\-1 & 0 \\ \end{pmatrix}$},\
$\eta\mapsto$ 
{\tiny 
$\begin{pmatrix}-1 & 0 \\0 & 1 \\ \end{pmatrix}$}.\
\smallskip

{\sf Projective $GL_2$.}
We define the {\em projective\,} $GL_{\,2}^\wedge(\Z)\,$ as
the group generated by $\tau_{\pm}, \eta$ subject to the
relations $\tau_+\tau_-^{-1}\tau_+\, =\,
\tau_-^{-1}\tau_+\tau_-^{-1},\ \eta^2=1$ and 
$\eta\tau_{\pm}\eta=\tau_{\pm}^{-1}.$ The span 
of $\tau_{\pm}$ is the projective $PSL_{\,2}^\wedge(\Z)\,$
(due to Steinberg), which is isomorphic to the braid 
group $B_3$.

\medskip

\subsection{\bf Macdonald polynomials}
Following \cite{C101}, 
we use the PBW Theorem to express any $H\in \HH$ in the form 
\,$\sum_{a,w,b} c_{a,w,b}\, X_a T_{w} Y_b$\, for $w\in W$,
$a,b\in P$ (this presentation is unique). Then we substitute:
\begin{align}\label{evfunct}
\{\,\}_{ev}:\ X_a \ \mapsto\  q^{-(\rho_k,a)},\ 
Y_b \ \mapsto\  q^{(\rho_k,b)},\ 
T_i \ \mapsto\  t_i^{1/2}. 
\end{align}

{\sf Polynomial representation.}
The functional $\,\HH\ni H\mapsto \{H\}_{ev}$, 
called {\em coinvariant\,}, acts via the projection 
$H\mapsto H\!\Downarrow \equal H(1)$ of $\HH\,$
onto the {\em polynomial representation \,}$\v$, which is
the $\HH$\~module induced from the one-dimensional
character $T_i(1)=t_i^{-1/2}=Y_i(1)$ for $1\le i\le n$ and
$T_0(1)=t_0^{-1/2}$. Recall  that $t_0=t_{\sht}$; 
see \cite{C101,CJ}.
\smallskip

In detail, the polynomial representation $\v$
is isomorphic to $\Z_{q,t}[X_b]$
as a vector space and the action of $T_i(0\le i\le n)$ there 
is given by 
the {\em Demazure-Lusztig operators\,}:
\begin{align}
&T_i\  = \  t_i^{1/2} s_i\ +\
(t_i^{1/2}-t_i^{-1/2})(X_{\al_i}-1)^{-1}(s_i-1),
\ 0\le i\le n.
\label{Demazx}
\end{align}
The elements $X_b$ become the multiplication operators 
and  $\pi_r (r\in O')$ act via the general formula
$\hw(X_b)=X_{\hw(b)}$ for $\hw\in \hW$. Note that $\tau_-$ 
and $\eta$ naturally act in the polynomial representation. 
For the latter,
\begin{align}\label{eta-conj}
\eta(f)\!=\!f^\star, \hbox{\, where\, } X_b^\star\!=\!X_{-b}, 
(q^\upsilon)^\star\!=\!q^{-\upsilon}, (t^v)^\star\!=\!t^{-v} 
\hbox{\, for\, } \upsilon\in \Q.
\end{align}

One has the following relations:
\begin{align}\label{evsym}
&\{\,\vph(H)\,\}_{ev}=\{\,H\,\}_{ev},\ 
\{\,\eta(H)\,\}_{ev}=\{\,H\,\}_{ev}^\star,\ 
\end{align}

{\sf Macdonald polynomials.}
The Macdonald polynomials $P_b(X)$ (due to Kadell
for the classical root systems) are uniquely defined
as follows. Let $c_+$ be the unique element 
such that $c_+\in W(c)\cap P_+$. For $b\in P_+$,
\begin{align*}
&P_b\! -\!\!\!\!\sum_{b'\in W(b)}\!\!\! X_{b'} 
\in\, \oplus_{b_+\neq c_+\in b+Q_+}\Q(q,t_\nu) X_c 
\hbox{\, and\, }
CT\bigl(P_b X_{c^\iota}\,\mu(X;q,t)\bigr)\!=\!0
\\
&\hbox{for such $c$,\, where\, } 
\mu(X;q,t)\equal\!\prod_{\al \in R_+}
\prod_{j=0}^\infty \frac{(1\!-\!X_\al q_\al^{j})
(1\!-\!X_\al^{-1}q_\al^{j+1})
}{
(1\!-\!X_\al t_\al q_\al^{j})
(1\!-\!X_\al^{-1}t_\al^{}q_\al^{j+1})}\,.
\end{align*}
Here $CT$ is the constant term; 
$\mu$ is considered
a Laurent series of $X_b$ with 
the coefficients expanded in terms of
positive powers of $q$. The coefficients of
$P_b$ belong to the $\Q(q,t_\nu)$.
One has:
\begin{align}\label{macdsym}
&P_b(X^{-1})\,=\,P_{b^\iota}(X)\,=\,
P_b^\star(X),\ \, P_{b}(q^{-\rho_k})=
P_{b}(q^{\rho_k})\\
\label{macdeval}
=&(P_{b}(q^{-\rho_k}))^\star=
q^{-(\rho_k,b)}
\prod_{\al>0}\prod_{j=0}^{(\al^{\!\vee},b)-1}
\Bigl(
\frac{
1- q_\al^{j}t_\al X_\al(q^{\rho_k})
 }{
1- q_\al^{j}X_\al(q^{\rho_k})
}
\Bigr).
\end{align}
Recall that $\iota(b)=b^\iota=-w_0(b)$ for $b\in P$.

DAHA provides an important alternative (operator) 
approach to $P$\~polynomials; namely, they satisfy
the (defining) relations
\begin{align}\label{macdopers}
L_f(P_b)=f(q^{\rho_k+b})P_b,\  L_f\equal f(X_a\mapsto Y_a)
\end{align}
for any symmetric ($W$\~invariant) polynomial 
$f\in \C[X_a,a\in P]^W$. The following {\em
evaluation formula\,} will be important below:

\begin{align*}
P_{b}(q^{-\rho_k})=&
q^{-(\rho_k,b)}\,\Pi_R^b\,
\prod_{\al>0}\prod_{j=1}^{(\al^{\!\vee},b)-1}
\Bigl(
\frac{
1- q_\al^{j}t_\al X_\al(q^{\rho_k})
 }{
1- q_\al^{j}X_\al(q^{\rho_k})
}\Bigr),\\
&\Pi_R^b=\!\!\!\prod_{\al>0,(\al,b)>0}
\frac{
1- t_\al X_\al(q^{\rho_k})
 }{
1- X_\al(q^{\rho_k})
}.
\end{align*}
The product $\Pi_R^b$ becomes the {\em Poincar\'e polynomial\,} 
$\Pi_R$ for generic $b\in P_+$. See formula (3.3.23) from 
\cite{C101}.
We set  $P_b^\circ\equal P_b/P_b(q^{-\rho_k})$ for $b\in P_+$,
this is the so-called {\em spherical normalization\,}.
\medskip

\setcounter{equation}{0}
\section{\sc DAHA-Jones theory}
\subsection{\bf Iterated torus knots}\label{sec:ITER-KNOTS}
The torus knots $T(\rr,\ss)$ are defined for any
integers assuming that \,gcd$(\rr,\ss)=1$.
One has the symmetries $\,T(\rr,\ss)=T(\ss,\rr)=T(-\rr,-\ss)$,
where we use ``$=$" for the ambient isotopy equivalence.
Also $\,T(\rr,\ss)=\unknot\,$ if $\,|\rr|\le 1$ or $|\ss|\le 1$
for the {\em unknot\,}, denoted by $\,\unknot$ here and below.
See e.g. \cite{RJ,EN} or Knot Atlas for the
details and the corresponding theory of the
invariants. 
 
The {\em iterated torus knots\,}, denoted in this paper by 
$\t(\vec{\rr},\vec{\ss})$,
will be associated with two sequences of
integers {\em of any signs\,}: 
\begin{align}\label{iterrss} 
\vec{\rr}=\{\rr_1,\ldots \rr_\ell\}, \ 
\vec{\ss}=\{\ss_1,\ldots \ss_\ell\} \hbox{\, such that\, 
gcd}(\rr_i,\ss_i)=1;
\end{align}
$\ell$ will be called the {\em length\,} of $\,\vec\rr,\vec\ss$.
The pairs $\{\rr_i,\ss_i\}$ will become 
{\em characteristic\,} or {\em Newton pairs\,}
for algebraic knots. We will use this name, but
generally consider arbitrary iterated knots in this paper. 

These are combinatorial data. The (topological) definition
of iterated torus knots requires one more sequence: 
\begin{align}\label{Newtonpair}
\aa_1=\ss_1,\,\aa_{i}=\aa_{i-1}\rr_{i-1}\rr_{i}+\ss_{i}\,\ 
(i=2,\ldots,m). 
\end{align}
See  e.g. \cite{EN} and \cite{Pi}. Then, in terms of
the cabling defined below:
\begin{align}\label{Knotsiter}
\t(\vec{\rr},\!\vec{\ss})\equal
C\!ab(\vec{\aa},\!\vec{\rr})(\unknot)=
\bigl(C\!ab(\aa_\ell,\rr_\ell)\cdots C\!ab(\aa_2,\rr_2)\bigr)
\bigl(T(\aa_1,\!\rr_1)\bigr).
\end{align}
Note that the first iteration (application of 
$C\!ab$) is for $\{\aa_1,\rr_1\}$ (not for the last pair!).
Recall that $\,\unknot$ is the
{\em unknot\,} and knots are considered up to 
{\em ambient isotopy\,} (we use \,``$=$" for it).
\smallskip

{\sf Cabling}.
The {\em cabling\,} $C\!ab(\aa,\bb)(K)$ of any oriented
knot $K$ in (oriented) $S^3$ is defined as follows; 
see e.g. \cite{Mo,EN} and references therein. 
We consider a small $2$\~dimensional torus
around $K$ and put there the torus knot $T(\aa,\bb)$
in the direction of $K$, 
which is $C\!ab(\aa,\bb)(K)$ (up to ambient isotopy).
We will sometimes set $C\!ab(\vec\aa,\vec\rr)= 
C\!ab(\vec\aa,\vec\rr)(\unknot)$.  

This procedure depends on 
the order of $\aa,\bb$ and orientation of
$K$. We choose the latter in the standard way (compatible 
with almost all sources, including  the Mathematica package
``KnotTheory"); the parameter $\aa$ gives the number of
turns around $K$. This construction also depends
on the framing of the cable knots; we take the natural
one, associated with the parallel copy of the torus where 
a given cable knot sits (its parallel copy 
has zero linking number with this knot). 

\comment{
We will sometimes use the following straightforward 
extension of the definition
from (\ref{Knotsiter}) to arbitrary knots $K$\,:
\begin{align}\label{TCKnots}
\t(\vec{\rr},\!\vec{\ss})(K)\!=\!
C\!ab(\vec{\aa},\!\vec{\rr})(K)\equal
\bigl(C\!ab(\aa_\ell,\rr_\ell)\cdots C\!ab(\aa_1,\rr_1)\bigr)
\bigl(K\bigr).
\end{align}
}
\smallskip

By construction,
$C\!ab(\aa,0)(K)=\unknot$ and  $C\!ab(\aa,1)(K)=K$ for any
knot $K$ and $\aa\neq 0$.
Accordingly, we have the following {\em reduction cases\,}:
\begin{align}\label{iter-triv}
&\hbox{When\, }
\rr_i=0\,\ (i<\ell),\ \, \t(\vec{\ss},\vec{\rr})=\\
\t(\{\rr_{i+1},\cdots,\rr_{\ell}&\},
\{\ss_{i+1},\cdots,\ss_{\ell}\}),\  
\t(\vec{\rr},\vec{\ss})=\unknot
\for i=\ell.\notag\\
&\hbox{When\, } \rr_i=1,\ \ss_i=0,\ 
\t(\vec\rr,\vec\ss)=\label{iter-s0}\\
\t(\{\rr_1,\cdots,\rr_{i-1}&,
\rr_{i+1},\cdots,\rr_\ell\},\{\ss_1,\cdots,\ss_{i-1},
\ss_{i+1},\cdots,\ss_\ell\}).\notag\\
&\hbox{When\, } \rr_i=1,\, \ss_i\neq 0,\ \,
\t(\vec\rr,\vec\ss)=\label{iter-s1}\\
\t(\{\rr_1,\cdots,\rr_{i-1}&,
\rr_{i+1},\cdots,\rr_\ell\},\{\cdots,\ss_{i-1},
\,\ss'_{i+1}\,,\ss_{i+2},\cdots\}),\notag\\
\hbox{where\,\,\, } \ss_{i+1}'\!=\!\ss&_{i+1}\!+\!\ss_{i}\rr_{i+1} 
\hbox{\,\, if\,\, } i <\ell\hbox{\,\, (no\, $s'_{\ell+1}$\, 
for\, } i=\ell). 
\end{align}
\smallskip

Also $T(\rr,\ss)=T(\ss,\rr)$, which results in
the {\em transposition property\,}:
\begin{align}\label{iter-sym}
&\t(\vec{\rr},\vec{\ss})=
\t(\{\ss_1,\rr_{2},\cdots,\rr_{\ell}\},
\{\rr_{1},\ss_2\cdots,\ss_{\ell}\}).
\end{align}
The counterparts of these properties hold for the 
DAHA-Jones polynomials, though due to entirely algebraic
reasons.
\smallskip

\comment{
The change of the orientation of $K$, denoted by $-K$,
is as follows:
\begin{align}\label{iterorienz}
C\!ab(\aa,\bb)(-K)=-&C\!ab(\aa,-\bb)(K),\ 
C\!ab(\vec\aa,\{-\rr_1,\ldots,-\rr_{i},
\rr_{i+1},\ldots \rr_{n}\})\notag\\
=\,C\!ab(\{\aa_{i+1},\ldots,\aa_n&\},\{\rr_{i+1},\ldots,\aa_n\})
\bigl(
-C\!ab(\{\aa_{1},\ldots,\aa_i\},\{\rr_1,\ldots,\rr_i\})\bigr),
\notag\\
\hbox{which gives}&\ 
\t(\{-\rr_1,\ldots,-\rr_{i},
\rr_{i+1}\ldots \rr_{n}\},\,\vec\ss)\\
=\,\t(\{\rr_{i+1}\ldots&,\rr_n\},\{\ss_{i+1}\ldots,\ss_{n}\})
\bigl(
-\t(\{\rr_1,\ldots,\rr_i\},\{\ss_{1},\ldots,\ss_i\})\bigr).\notag
\end{align}
Replacing $K$ by its {\em mirror image\,}, 
denoted by $K^\star$, one has:
\begin{align}\label{iterorien}
C\!ab(\aa,\bb)(K^\star)=
&\bigl(C\!ab(-\aa,\bb)(K)\bigr)^\star \hbox{\,\, for any }
\aa,\bb \hbox{ with\,\,   gcd}(\aa,\bb)\!=\!1,\notag\\
C\!ab(\{-\aa_1,\ldots,&-\aa_{i},\aa_{i+1},\ldots,\aa_n\},
\{\rr_{1},\ldots \rr_{n}\})\notag\\
=C\!ab(\{\aa_{i+1},\ldots,&\aa_n\},\!\{\rr_{i+1},\ldots,\rr_n\})
\bigl(
C\!ab(\{\aa_{1},\ldots,\aa_i\},
\!\{\rr_1,\ldots,\rr_i\})^\star\bigr),
\notag\\
\t(\vec\rr,\, \{-\ss_1,&\ldots,-\ss_{i},
\ss_{i+1}\ldots \ss_{n}\})\\
=\,\t(\{\rr_{i+1},\ldots&,\rr_n\},\{\aa_{i+1},\ss_{i+2},
\ldots,\ss_{n}\})\,
\bigl(\,
\t(\{\rr_1,\ldots,\rr_i\},\{\ss_{1},\ldots,\ss_i\})^\star\,\bigr).
\notag
\end{align}
}
\smallskip

Switching from $K$ to its {\em mirror image\,}, 
denoted by $K^\star$, one has:
\begin{align}\label{iterorien}
C\!ab(\aa,\bb)(K^\star)&=
\bigl(C\!ab(-\aa,\bb)(K)\bigr)^\star \hbox{\,\, for any }
\aa,\bb \hbox{ with\,\,   gcd}(\aa,\bb)\!=\!1,\notag\\
\hbox{and\, }&C\!ab(-\vec\aa,\vec\rr)\!=\!
\bigl(C\!ab(\vec\aa,\vec\rr)\bigr)^\star,\ 
\t(\vec\rr, -\vec\ss)\!=\!
\bigl(\t(\vec\rr,\vec\ss)\bigr)^\star.
\end{align}

Changing the {\em orientation\,}, denoted by ``$-$", at the 
$i$\hbox{\tiny th}\, step, we obtain 
that for any $1\le i\le \ell$,
\begin{align}\label{iter2min}
-C\!ab(\vec\aa,\vec\rr)\!=\!
C\!ab&(\{\ldots\!,\aa_{i\!-\!1},\!-\aa_i,\aa_{i\!+\!1},\!\ldots\},
\{\ldots\!,\rr_{i\!-\!1},\!-\rr_i,\rr_{i\!+\!1},\!\ldots\}),
\notag\\
-\t(\vec\rr,\vec\ss)\!=\!\t&(\{\ldots,\rr_{i\!-\!1},
\!-\rr_i,\rr_{i\!+\!1},\ldots\},
\{\ldots\!,\ss_{i\!-\!1},\!-\ss_i,\ss_{i\!+\!1},\!\ldots\} ).
\end{align} 

We note that the Jones and HOMFLYPT polynomials for 
$K^\star$ are obtained from those for $K$ by the formal
conjugation of the parameters, which is 
$q\mapsto q^{-1}, a\mapsto a^{-1}$. This will hold
for the DAHA-Jones polynomials and
DAHA superpolynomials (the conjugation of $t$ will
be added to that of $q,a$). The orientation does not
influence our construction.
\medskip

\subsection{\bf Algebraic knots}\label{sec:ALG-KNOTS}
There exists a deep connection of the iterated torus knots
for strictly positive $\rr_i,\ss_i$
with the germs of isolated {\em plane curve singularities\,}. 
Its origin is the Newton's successive approximations 
for $y$ in terms
of $x$ satisfying a polynomial equation $f(x,y)=0$
in a neighborhood of $0=(x=0,y=0)$. The main claim is that
the knot $\t(\vec\rr,\vec\ss)$ for the 
{\em characteristic pairs\,} $\{\rr_i,\ss_i\}$
with $\rr_i,\ss_i>0\,$ is the link 
of the {\em germ of the singularity\,} 
\begin{align}\label{yxcurve}
y = x^{\frac{\ss_1}{\rr_1}}
(c_1 + x^{\frac{\ss_2}{\rr_1\rr_2}}
\bigl(c_2 + \ldots + 
x^{\frac{\ss_\ell}{\rr_1\rr_2\cdots\rr_\ell}}\bigr))
\hbox{\, at\, } 0,
\end{align} 
which is the intersection of the corresponding plane curve 
$f(x,y)=0$ with a small $3$-dimensional sphere in 
$\C^2$ around $0$. We will always
assume that this germ is {\em unibranch}.

The inequality  $\ss_1<\rr_1$ is commonly
imposed here (otherwise $x$ and $y$ can be switched).
Formula (\ref{yxcurve}) is the celebrated 
{\em Newton-Puiseux expansion}, though Puiseux 
performed the multiplication here:
\begin{align*}
y\!=\!b_1\, x^{\frac{\mm_1}{\rr_1}}\!+\!
b_2\, x^{\,\frac{\mm_2}{\rr_1\rr_2}}\!+\!
b_3\, x^{\,\frac{\mm_3}{\rr_1\rr_2\rr_3}}\!+\ldots,\ \,
\mm_1\!=\!\ss_1,\, \mm_i\!=\!\ss_i\!+\!\rr_i\mm_{i-1},\ldots\,,
\end{align*} 
where $\{\rr_i,\mm_i\}$ are called the  
{\em Puiseux pairs\,} (which we will not use).

We will denote the germ of
the curve in (\ref{yxcurve}) by $\,\c_{\vec\rr,\vec\ss}$. 
If it has more than one branch (near $0$), then
its complete topological description requires 
knowing the sequences $\{\vec{\rr},\vec{\ss}\}$ for all components 
together with the pairwise linking numbers. See e.g. \cite{Mi}
and pg. 49 of \cite{EN}. All algebraic knots can be obtained
in such a way; the positivity of $\,\rr_i,\ss_i\,$ is 
necessary for $\,\t(\vec\rr,\vec\ss)\,$ to be an algebraic link.
\smallskip

{\sf Jacobian factors.}
One can associate with a unibranch $\,\c_{\vec\rr,\vec\ss}\,$ the
{\em Jacobian factor\,} $J(\c_{\vec\rr,\vec\ss})$.
Up to a homeomorphism, it can be introduced as the canonical 
compactification of the {\em generalized Jacobian\,} of an
integral {\em rational\,} planar curve with 
$\,\c_{\vec\rr,\vec\ss}\,$
as its {\em only\,} singularity \cite{La}; there is a
purely local definition as well. Its dimension is 
the {\em $\de$\~invariant\,} of the singularity
$\,\c_{\vec\rr,\vec\ss},$ also called its 
{\em arithmetic genus}.

The Jacobian factors are in the focus 
of many studies; see e.g. \cite{Bea}, \cite{FGV} (Theorem 1), 
\cite{Wa}, \cite{Pi} and \cite{La}.  The latter paper
were the connection with the 
{\em affine Springer fibers\,} was established,
provides a link to the theory of DAHA modules \cite{VV,Yun},
though we cannot relate this so far to our using DAHA. 
An important connection  
of $\,J(\c_{\vec\rr,\vec\ss})\,$ and the {\em rational DAHA\,} 
was established in \cite{ORS,GORS} (see references there),
which must be mentioned, but this is not connected with
our approach and that from \cite{CJ,CJJ} so far.

Calculating the {\em Euler number\,}  
$e(J(\c_{\vec\rr,\vec\ss}))$, the topological Euler
characteristic of $J(\c_{\vec\rr,\vec\ss}),$
and the corresponding {\em Betti numbers\,} 
in terms of $\,\vec\rr,\vec\ss\,$ is a challenging problem. 
For torus knots $T(\rr,\ss)$, one has 
$e(J(\c_{\rr,\ss}))=$ $\frac{1}{\rr+\ss}
\binom{\rr+\ss}{\rr}$ due to \cite{Bea}. The 
Catalan numbers here are related to the perfect modules 
of rational DAHA; see e.g. \cite{GM}. 

\smallskip
The Euler numbers of $J(\c_{\vec\rr,\vec\ss})$ were calculated 
in \cite{Pi} (the Main Theorem) for the following triples of 
{\em Puiseux characteristic exponents\,}: 
\begin{align}\label{PiEuler}
(4,2u,v),\, (6,8,v)\, \hbox{\, and\, } 
(6,10,v)  \hbox{\, for odd\, } u,v>0,
\end{align}
where $4\!<\!2u\!<v$, and $8\!<\!v, 10\!<\!v$ respectively.
Here $\de\!=$dim\,$J(\c_{\vec\rr,\vec\ss})$
is $\frac{(\rr-1)(\ss-1)}{2}$ for $T(\rr,\ss)$ and 
$2u+(v-1)/2-1$ for the series $(4,2u,v)$. Generally,
$\de$ equals the cardinality
$|\N\setminus \Gamma|$, where $\Gamma$ is the {\em valuation
semigroup\,} associated with $\c_{\vec\rr,\vec\ss}$;
see \cite{Pi} and \cite{Za}. The Euler numbers of
the Jacobian factors (any isolated plane curve singularities)
can be also calculated
now via the HOMFLYPT polynomials of 
the corresponding links (see below) due to \cite{ObS,Ma}.
\smallskip

Concerning the Betti numbers for the torus knots and
the series $(4,2u,v)$, the odd (co)homology 
of $\,J(\c_{\vec\rr,\vec\ss})\,$ vanishes and
the formulas for Betti numbers 
$\,h_{2k}=\hbox{rank\,}(H_{2k}(J(\c_{\vec\rr,\vec\ss})))$
are known for some (sufficiently close to $0$ or $\de$) 
values of $k$, where $0\le k\le \de$. Not much is 
known/expected beyond these two series. 
Paper \cite{Pi} contains several results and 
conjectures for the Betti numbers of Jacobian factors; 
for instance there are quite a few formulas there for 
$T(\rr,\ss)$ with min$(\rr,\ss)\le 4$. 
\smallskip
 
{\sf What DAHA provides.}
Conjecturally, the coefficient of the DAHA superpolynomial
$\h_{\{\vec\rr,\vec\ss\,\}}$ (to be defined below)
at $a=0$ and for $ q=1$ in the uncolored case 
is $\sum_{k=0}^{\de}h_{2k}t^k$ for {\em any\,} unibranch 
$\c_{\vec\rr,\vec\ss}$,
which matches the tables/conjecture from \cite{Pi}. 
This includes the {\em van Straten-Warmt conjecture\,}, 
which claims that all odd Betti numbers 
of $J(\c_{\vec\rr,\vec\ss})$ vanish, since
the DAHA superpolynomials contain only integral powers of $q,t,a$.    
In the opposite direction, this would
imply the positivity of $\h_{\{\vec\rr,\vec\ss\,\}}$ at
$a=0,q=1$. The whole (uncolored) $\h_{\{\vec\rr,\vec\ss\,\}}$
describes the geometry of Hilbert schemes of 
$\,\c_{\vec\rr,\vec\ss}\,$ due to 
Conjecture 2 from \cite{ORS}.

The general positivity conjecture claims that the
coefficient of the whole $\h_{\{\vec\rr,\vec\ss\,\}}$ 
are positive for any rectangle Young diagrams and algebraic 
knots. Assuming such a positivity
(the corresponding knots will be called
{\em pseudo-algebraic\,}), $\h_{\{\vec\rr,\vec\ss\,\}}$  
is expected to coincide with the {\em stable reduced 
Khovanov-Rozansky polynomial\,} of
$\t(\vec\rr,\vec\ss)$, producing (by definition) the 
Khovanov-Rozansky polynomials for $sl_n$ upon $a=-t^{n+1}$
for sufficiently large $n$. The substitution $a=-1$ here
is expected to result in the {\em Heegard-Floer homology\,}.
This is stated in the DAHA parameters $q,t,a$;  
see Parts $(ii,iii)$  of Conjecture \ref{CONCONJ}. 
We restrict ourselves there to the uncolored case, 
because this is needed in \cite{ORS} and in the definition  
of the Khovanov-Rozansky triply-graded homology. 
Though see \cite{Web} and other
papers on the {\em categorification\,} (via Quantum Group).
%The theory of DAHA superpolynomials is for any weights.  

%Also, let us mention that the so-called 
%{\em Puiseux pairs\,} are somewhat different from the 
%Newton pairs, though connected with the latter by a simple
%algebraic transformations. 
\smallskip

\subsection{\bf DAHA-Jones polynomials}
The following theorems are mainly an extension of 
Theorem 1.2 in \cite{CJJ} on DAHA-Jones polynomials 
from torus knots  to arbitrary iterated torus knots.

%We enlarge $PSL_{\,2}^{\wedge}(\Z)$ used
%there to $GL_{\,2}^{\wedge}(\Z)$, considering in detail 
%the symmetry $(\rr,\ss)\mapsto (\ss,\rr)$  and
%mirroring (\ref{jones-sym}), and add important color-exchange 
%relations (\ref{jones-bc}).
%Last but not least, the justifications are provided.

Torus knots $T(\rr,\ss)$ are naturally represented   
by $\ga_{\rr,\ss}\in GL_{\,2}(\Z)$ with the 
first column $(\rr,\ss)^{tr}$ ($tr$ is the transposition)
for $\,\rr,\ss\in \Z$ 
assuming that \,gcd$(\rr,\ss)=1$. 
Then let $\hat{\ga}_{\rr,\ss}\in GL_{\,2}^{\wedge}(\Z)$ be
any pullback of $\ga_{\rr,\ss}$.

Note that $\,(\rr,\ss)\,$ can be obviously lifted
to $\,\ga\,$ of determinant $1$ and, accordingly, 
to $\hat{\ga}\,$ from the subgroup $PSL^\wedge_{\,2}(\Z)$
generated by $\{\tau_{\pm}\}$, i.e. without using
$\eta$. This is actually sufficient for the construction 
of the DAHA-Jones polynomials below.
% using $GL_2(\Z)$ 
%instead of $PSL_2(\Z)\subset GL_2(\Z)$ 
%is actually equivalent to the conjugation relation 
%from (\ref{jones-sym-iter}).
\smallskip

For a polynomial $F$ in terms of 
fractional powers of $q$ and $t_\nu$, 
the {\em tilde-normalization}
$\tilde{F}$ will be the result of the division of $F$
by the lowest $q,t_\nu$\~monomial, assuming that it 
is well defined. We put $q^\bullet t^\bullet$ for
a monomial factor (possibly fractional)
in terms of $q,t_\nu$.
%thus $\tilde{F}$ contains only non-negative powers (possibly
%rational) of $q,t_\nu$. 
%We will constantly
%use the tilde-normalization below, though not all properties of 
%the DAHA-Jones polynomials really require such a 
%normalization (see \cite{CJJ}). 

\begin{theorem}\label{JONITER}
Recall that $H\mapsto H\!\!\Downarrow\, \equal H(1)$, where the
action of $H\in \HH$ in $\v$ is used.
Given two sequences $\vec\rr,\vec\ss$ of length $\ell$ as 
in (\ref{iterrss}), a reduced root system $R$ and a weight 
$b\in P_+\,$, we lift $(\rr_i,\ss_i)^{tr}$ to
$\ga_i$ and then to $\hat{\ga_i}\in GL_{\,2}^{\wedge}(\Z)$ as above.
The {\sf DAHA-Jones polynomial} is
\begin{align}\label{jones-dit}
&J\!D_{\,\vec\rr,\,\vec\ss\,}^{R}(b\,;\,q,t)\ =\ 
J\!D_{\,\vec\rr,\,\vec\ss\,}(b\,;\,q,t)\, \equal\\ 
\Bigl\{\hat{\ga_{1}}&\Bigr(
\cdots\Bigl(\hat{\ga}_{\ell-1}
\Bigl(\bigl(\hat{\ga}_\ell(P_b)/P_b(q^{-\rho_k})\bigr)\!\Downarrow
\Bigr)\!\Downarrow\Bigr) \cdots\Bigr)\Bigr\}_{ev}.\notag
\end{align}
It does not depend
on the particular choice of the lifts \,$\ga_i\in GL_2(\Z)$
and $\hat{\ga}_i\in GL_{\,2}^{\wedge}(\Z)$\, for $1\le  i
\le\ell.$ Assuming that $\rr_i,\ss_i>0$ for all $i$,
i.e. $\t(\vec\rr,\vec\ss)$ is algebraic,
the tilde-normalization 
$\tilde{J\!D}_{\vec\rr,\vec\ss}\,(b\,;\,q,t)$ is well defined
and is a polynomial in terms of $\,q,t_\nu$ with
constant term $1$. 
\sq
\end{theorem}

The proofs of this and the next theorem almost exactly
follow those in \cite{CJ, CJJ} (for torus knots).

%\smallskip
\begin{theorem} \label{MAINTHM} 
(i) The polynomials
$J\!D_{\vec\rr,\vec\ss}(b\,;\,q,t)$  
coincide up to $q^\bullet t^\bullet$
for the pairs from (\ref{iter-triv}-\ref{iter-s1}), i.e.
for $\{\vec{\rr},\vec{\ss}\}$ 
with isotopic $\t(\vec\rr,\vec\ss)$. Also,
(\ref{iter-sym}) corresponds to
\begin{align}\label{iter-sym-j}
&J\!D_{\,\vec{\rr},\,\vec{\ss}}\,(b\,;\,q,t)=
J\!D_{\{\ss_1,\rr_{2},\cdots,\rr_{\ell}\},
\{\rr_{1},\ss_2\cdots,\ss_{\ell}\}}\,(b\,;\,q,t),
\end{align}
and the following DAHA counterparts of
(\ref{iterorien}) and (\ref{iter2min}) hold:
\begin{align}
\label{jones-sym-iter} 
&J\!D_{\,\vec\rr,\,-\vec\ss}\,(b\,;\,q,t) = 
\bigr(J\!D_{\,\vec\rr,\,\vec\ss}\,
(b\,;\,q,t)\bigl)^\star \hbox{\,for\, } \star \hbox{\, 
from (\ref{eta-conj})},
\\
&J\!D_{\,\vec\rr,\,\vec\ss}=
J\!D_{\{\ldots,\rr_{i\!-\!1},
-\rr_i,\rr_{i\!+\!1},\ldots\},
\{\ldots,\ss_{i\!-\!1},-\ss_i,\ss_{i\!+\!1},\ldots\}}
\hbox{\,\ for\,\, } 1\!\le\! i\!\le\! \ell.\notag
\end{align}

(ii) Let us assume that for $b,c\in P_+$ and certain 
\,id\,$\neq w\in W$,  
\begin{align}\label{bcw-rel}
q^{\,(b+\rho_k-w(\rho_k)-w(c)\,,\,\al)\,}=1=
q_{\al}^{\,(b-w(c)\,,\,\al^{\!\vee})}\,
t_{\sht}^{\,(\rho^w_{\sht}\,,\,\al^{\!\vee})}\,
t_{\lng}^{\,(\rho^w_{\lng}\,,\,\al^{\!\vee})}
\end{align}
for any $\al\in R_+$, where we set\,
$\rho_{\nu}^w\equal w(\rho_{\nu})-\rho_{\nu}$. Then 
\begin{align}\label{jones-bc}
J\!D_{\,\vec\rr,\,\vec\ss}\,(b\,;\,q,t)=
J\!D_{\,\vec\rr,\,\vec\ss}\,(c\,;\,q,t)
\hbox{\, for such \, } q,\{t_\nu\} 
\end{align}
and any $\vec\rr,\vec\ss$. Also, for 
$b=\hbox{\small$\sum$}_{i=1}^n b_i \om_i$,
\begin{align}\label{jones-eval}
J\!D_{\,\vec\rr,\,\vec\ss}
\,\bigl(b\,;\,q\!=\!1,t
\bigr)\!=\!
\hbox{\small$\prod$}_{i=1}^n J\!D_{\,\vec\rr,\,\vec\ss}\,
(\om_i\,;\,q\!=\!1,t)^{b_i} \hbox{\, for any \,} \vec\rr,\vec\ss.
\end{align}
\vskip -1.25cm \sq
 \end{theorem}

Extending the connection conjectures from \cite{CJ}, we
expect that 
$J\!D_{\,\vec\rr,\,\vec\ss}\,(b\,;\,q, t_\nu
\!\mapsto\! q_\nu)$ 
coincide up to $q^\bullet$ 
with the reduced Quantum Group (WRT) invariants for the
corresponding  $\t(\vec\rr,\vec\ss)$ and any $b\in P_+$.
The Quantum group is associated with the root system 
$\tR$; see \cite{CJ}. The {\em shift operator\,} was used there
to deduce this coincidence from \cite{Ste} in the case
of $A_n$ and torus knots; quite a few confirmations were provided 
for other root systems (including special ones). 

Using the shift operator for the iterated knots is almost 
exactly the same as for the torus knots. However, the 
Quantum Group invariants for the
iterated knots are generally much less studied than those
for the torus knots. In the case of $A_1$, we will check this
using the formulas from \cite{RJ,Mo}. We hope 
that $A_n$ will be doable too. 
%See 
%Conjecture \ref{CONCONJ} below and the discussion there. 
\medskip

\subsection{\bf DAHA superpolynomials}
Following \cite{CJ,GoN,CJJ},
Theorems \ref{JONITER}, \ref{MAINTHM} can
be extended to the {\em DAHA- superpolynomials\,}
the result of the stabilization  of 
$J\!D^{A_n}(b;q,t)$ with respect to $n$.
This stabilization was announced 
in \cite{CJ}; its proof was published in \cite{GoN}. 
Both approaches use \cite{SV} and can be 
extended to arbitrary iterated knots, the Duality Conjecture 
proposed in \cite{CJ})
and proven in \cite{GoN}.

We mainly omit the justifications in this paper; with few
reservations, they are similar to those in \cite{GoN,CJJ}
(for torus knots). The focus of this paper is on the
connections with the Khovanov-Rozansky theory and the
geometry of (germs of) plane curve singularities.
\smallskip

The sequences $\vec\rr,\vec\ss$ will be from the
previous sections ($\ell$ is the length), as well 
%However, {\em from now on we assume that $\,\rr_i \aa_i\ge 0\,$
%for all $1\le i\le \ell$}. This is (much) weaker than
%the positivity of $\rr_i,\ss_i$ for all $\,i\,$
%necessary and sufficient for $\t(\vec\rr,\vec\ss)$ to be algebraic.
as the iterated knots $\t(\vec\rr,\,\vec\ss)$,
DAHA-Jones polynomials $\tilde{J\!D}_{\,\vec\rr,\,\vec\ss\,}$
and (later) the Jacobian factors $J(\c_{\vec\rr,\vec\ss})$
of germs $\,\c_{\,\vec\rr,\,\vec\ss}\,$ 
of unibranch isolated plane curve singularities 
(Section \ref{sec:ALG-KNOTS}).

\begin{theorem}\label{STABILIZ}
We switch to $A_n$, setting $t=t_{\sht}=q^k$. 
Let us consider $P_+\ni b=$ $\sum_{i=1}^n b_i \om_i$ as 
a (dominant) weight for any $A_N$ with $N\ge n-1$,
where we set $\om_{n}=0$ upon the restriction to $A_{n-1}$.

(i) {\sf Stabilization.} Given $\t(\vec\rr,\vec\ss)$, there exists 
$\h_{\,\vec\rr,\,\vec\ss}\,(b\,;\,q,t,a)$,
a polynomial from $\Z[q,t^{\pm 1},a]$, satisfying the relations
%such that
%its coefficient of $a^0$ is tilde-normalized 
%(i.e. in the form $\sum_{u,v\ge 0}C_{u,v}q^u t^v$ 
%with $C_{0,0}=1$) and the following holds:
\begin{align}\label{jones-sup}
\h_{\,\vec\rr,\,\vec\ss}\,(b;q,t,a\!=\!-t^{N\!+\!1})\!=\!
\pm q^\bullet t^\bullet\, 
J\!D_{\,\vec\rr,\,\vec\ss\,}^{A_{N}}(b;q,t)
\hbox{\, for any } N\!\ge\! n\!-\!1, 
\end{align} 
and normalized as follows. The polynomial
$\,\h(a\!=\!0)\,$ is assumed minimal, i.e. not divisible by 
$\,q,\,t\,$ and 
any prime number; the coefficient of the minimal power of $t$ in
$\,\h(q\!=\!0,a\!=\!0)\,$ is selected positive.

(ii) {\sf Symmetries.}  The corresponding normalization of
$J\!D$\~polynomials, denoted by $ \tilde{J\!D}$, 
extends the tilde-normalization for algebraic knots
from Theorem \ref{JONITER}. For sufficiently large $N$
(for all $N\ge n-1$ if $\rr_i,\ss_i>0$),  
\begin{align}\label{tilde-sup}
\h_{\,\vec\rr,\,\vec\ss}\,(b;q,t,a\!=\!-t^{N\!+\!1})\!=\! 
\tilde{J\!D}_{\,\vec\rr,\,\vec\ss\,}^{A_{N}}(b;q,t).
\end{align} 
The symmetries from Part (i) of Theorem \ref{MAINTHM}
hold for $\h_{\,\vec\rr,\,\vec\ss}$ under such a normalization,
but (\ref{jones-sym-iter}); the latter holds  
up to $\,a^\bullet q^\bullet t^\bullet\,$
for $a^\star=a^{-1}$.
%For the latter we set 
%$\,a^\star= a^{-1}\,$ and then multiply 
%$\h_{\,\vec\rr,\,\vec\ss}^\star$ by $\,a^{deg}$ and proper
%$q^\bullet t^\bullet$,
%where\, deg=deg$_a$\, is the $a$\~degree of  
%$\h_{\,\vec\rr,\,\vec\ss}$.

(iii) {\sf Color Exchange.} Imposing (\ref{bcw-rel}), 
we consider $\,w\,$ 
%from Part (ii) of Theorem  \ref{MAINTHM} 
as an element of $\S_{N+1}$ for every
$N\ge n$ (the Weyl group for $A_N$) naturally acting in 
the corresponding $P$. Then for any \, $\vec\rr,\vec\ss$,
\begin{align}\label{jones-bca}
\h_{\,\vec\rr,\,\vec\ss}\,(b\,;\,q,t,a) = 
\h_{\,\vec\rr,\,\vec\ss}\,(c\,;\,q,t,a)
\hbox{\, for such \, } q,t.
\end{align}
In particular, let $w=s_i=(i,i+1)$ with $i<n$
and $\,b=c-\bigl(k+(c,\al_i)\bigr)\al_i$ for
a dominant $c$. 
Then the inequalities
$\,c_i/2\!\le\! -k\!\le\! c_i+\min\{c_{i\pm 1>0}\}\,$
for $k\!\in\! -\Z_+$ are sufficient.

(iv) {\sf Specialization and deg${}_a$}.
Making $q=1$, one has:
\begin{align}\label{jones-seval}
\h_{\,\vec\rr,\,\vec\ss}\,(b\,;\,q\!=\!1,t,a)\!=\!
\prod_{i=1}^n\h_{\,\vec\rr,\,\vec\ss}\,
(\om_i\,;\,q\!=\!1,t,a)^{b_i}
\hbox{\, for\, } b\!=\!\sum_{i=1}^n b_i \om_i.
\end{align}
Assuming  that $|\rr_1|>|\ss_1|$, the $a$-degree
deg${}_a \h_{\,\vec\rr,\,\vec\ss}\,(b\,;\,q,t,a)\,$
is no greater than $(|\ss_1\rr_2\cdots\rr_\ell|-1)$ times 
\, ord $(\la_b)$, the 
number of boxes in the Young diagram $\la_{b}$
associated with $b\in P_+$. 

(v) {\sf Super-duality \cite{GS,CJ,GoN}.} 
Up to a power of $q$ and $t$,
\begin{align}\label{iter-duality}
\h_{\,\vec\rr,\,\vec\ss}\,(\la\,;q,t,a)=q^{\bullet}t^{\bullet}
\h_{\,\vec\rr,\,\vec\ss}\,(\la^{tr}\,;t^{-1},q^{-1},a),
\end{align}
where we switch from (arbitrary) dominant
weights $\,b\,$ to the corresponding Young diagrams $\la=\la_b$,
and $\la^{tr}$ is the transposition of $\la$. 
\end{theorem} 
\vskip -1.25cm \sq
\vskip 0.2cm

{\sf Color exchange.}
See \cite{CJJ} for a systematic consideration of the
color exchange condition (\ref{jones-bca}). In terms
of the Young diagrams it is as follows.
We associate with
$c=\sum_{i=1}^n c_i\om_i$ in Part $(iii)$ the Young diagram
$$
\la_c=\{m_1\!=\!c_1+\ldots+c_n,m_2\!=\!c_1+\ldots+c_{n-1},
\ldots, m_n\!=\!c_n,0,0,\ldots\}.
$$
Then we switch to $\la_c'=\{m'_i=m_i-k(i-1)\}$, apply $w\in W$ to
$\la_c'$  and finally obtain 
\begin{align}\label{lamw}
\la_b=\{m'_{w(i)}+k(i-1),\}=\{m_{w(i)}+k(i-w(i))\}.
\end{align} 
Here 
$w$ transforms the rows of $\la$ and we set 
$w\{m_1,m_2,\ldots,m_n\}=\{m_{w(1)},m_{w(2)},\ldots,m_{w(n)}\}$.
Given $k<0$ (it can be fractional), $\la_b$ 
{\em must be a Young diagram\,}; this condition 
determines which $\,w\,$ can be used.
\smallskip

{\sf Toward $C^\vee C$.}
The second part of \cite{CJ} is devoted to the DAHA-Jones
polynomials for the root system $C^\vee C_1$, which
prepared the extension of the DAHA-Jones theory
to all $C^\vee C_n$.
We expect here a relatively straightforward switch to 
arbitrary iterated knots.  There is
an interesting application of DAHA-Jones polynomials
for $C^\vee C_1$ to the superpolynomials of type $A$
for knots $T(2\mm+1,2)$
(i.e. any $A_n$). We see relations of this kind
for certain iterated knots, say 
$C\!ab(2\mm+1,2)T(3,2)$,
but do not know how far this can go. Cf. \cite{DMMSS}.
\smallskip

{\sf The $a$-degree.}
We expect that \,deg${}_a \h_{\,\vec\rr,\,\vec\ss}\,
(b\,;\,q,t,a)\,$
from Part $(iv)$ equals 
$\,(\ss_1\rr_2\cdots\rr_\ell-1)\,\hbox{ord}(\la_b)$
for algebraic knots (there is no such a
coincidence in general). 
Due to (\ref{jones-seval}), it suffices to check the
equality for the fundamental weights only.
This equality is likely to hold
for the HOMFLYPT polynomials of algebraic 
knots (the inequality from Part $(iv)$ holds for all
iterated knots), which would formally result in such
an equality for $\h$ assuming Part $(i)$ of 
Conjecture \ref{CONCONJ} below.
\medskip

\subsection{\bf Connection Conjecture}
Let us briefly discuss the HOMFLYPT and Khovanov-Rozansky
polynomials. See \cite{CJ} for some details. 

After papers \cite{GoN,CJJ},
the last unproven ``DAHA-intrinsic" 
conjecture from \cite{CJ} is the positivity of 
all $q,t,a$\~coefficients of  (reduced) 
$\,\h_{\rr,\ss}(b\,;\,q,t,a)\,$ for any rectangle $\la_b$. 
We expect this 
to hold for $\h_{\,\vec\rr,\,\vec\ss}\,(b\,;\,q,t,a)$ for 
any {\em algebraic\,} knots.

However we need to restrict Part $(ii)$ of the conjecture below
to the uncolored case, since the
Khovanov-Rozansky triply-graded homology \cite{KhR1,KhR2,Kh,Ras}
requires this. The corresponding
{\em reduced\,} Poincar\`e polynomial, the so-called {\em stable 
Khovanov-Rozansky polynomial\,}, 
will be denoted by  $K\!h\!R_{\hbox{\tiny stab}}$ in what will
follow.
%\smallskip

The polynomials $K\!h\!R_{\hbox{\tiny stab}}$
are expected to coincide with the
(reduced) physics {\em superpolynomials}  
based on the {\em BPS states\,} \cite{DGR,AS,FGS,GGS}
(their theory is not rigorous)   
and those obtained in terms of {\em rational DAHA\,}
\cite{GORS,GoN} for torus knots. The latter are
known only in the uncolored case so far; the case of symmetric 
powers of the fundamental representation is in progress, 
see \cite{GGS}. %We will not touch these connections.

\smallskip

We mention that {\em categorification theory\,} 
can generally deal with arbitrary 
weights, but the categorification of the HOMFLYPT 
polynomials is not done in such a generality. 
See e.g. \cite{Kh,Rou,Web} and references there.
Also our construction is {\em reduced\,} 
(the normalization is by $1$ at $\unknot$, necessary for the  
polynomiality of DAHA-Jones polynomials),
which creates additional difficulties for the categorification.
\smallskip

\begin{conjecture}{\sf Connection Conjecture.}\label{CONCONJ}
Let 
$\,\h_{\,\vec\rr,\,\vec\ss}\,(b;q,t,a)_{st}\,$
denote  $\,\h_{\,\vec\rr,\,\vec\ss}\,(b;q,t,a)\,$ expressed
in terms of the standard parameters (in the Khovanov-Rozansky 
theory; see \cite{CJ} and Section 1 in \cite{ORS}):
\begin{align}\label{qtareli}
&t=q^2_{st},\  q=(q_{st}t_{st})^2,\  a=a_{st}^2 t_{st},\notag\\
&q^2_{st}=t,\  t_{st}=\sqrt{q/t},\  a_{st}^2=a\sqrt{t/q}.
\end{align}

(i) We conjecture that for any $\,\vec{\rr},\vec{\ss}\,$
as above (possibly negative), 
\begin{align}\label{conjhomfly}
&\h_{\,\vec\rr,\,\vec\ss}\,(b\,;
\,q,t\!\mapsto\! q,a\!\mapsto\! -a)=
\tilde{
\hbox{\small H\!O\!M}}_{\,\vec\rr,\,\vec\ss}\,(\la_b\,;\,q,a),
\end{align}
where the latter is the tilde-normalization of the (reduced)
HOMFLYPT polynomial; $b$ is an arbitrary dominant weight 
represented by the Young diagram $\la_b$. 
We conjecture that the $a$-degrees of these polynomials 
equal to $(\ss_1 \rr_2\cdots \rr_\ell-1)\hbox{\,ord\,}(\la_b)$
provided the positivity of $\ss_i,\rr_i$
and $\ss_1<\rr_1$.

(ii) For\, uncolored {\sf pseudo-algebraic knots\,} (those
with positive coefficients of uncolored DAHA superpolynomials), 
one has:
\begin{align}\label{khrconj}
\h_{\,\vec\rr,\,\vec\ss}\,(\square\,;\,q,t,a)_{st}=
\tilde{K\!h\!R}_{\hbox{\tiny stab}}(q_{st},t_{st},a_{st}),\
\hbox{\, where\ \,} \square=\hbox{unknot},
\end{align}
$\tilde{K\!h\!R}_{\hbox{\tiny stab}}$ is (reduced)
$K\!h\!R_{\hbox{\tiny stab}}$ divided  
by the smallest power of $a_{st}$ and by
$q_{st}^\bullet t_{st}^\bullet\,$ such that
$\,\tilde{K\!h\!R}_{\hbox{\tiny stab}}(a_{st}\!=\!0)
\in \Z_+[q_{st},t_{st}]$ with the constant term $1$.

(iii) Following Section \ref{sec:ALG-KNOTS}, let
$\,\c_{\vec\rr,\vec\ss}$ be the (unibranch) germ
at $0$ of the plain curve singularity corresponding
to $\{\rr_i,\ss_i>0\}$ and $\,\{h_j\}\,$ the Betti numbers of 
the Jacobian factor $\,J(\c_{\vec\rr,\vec\ss})$. Then
Conjecture 2 from \cite{ORS} and claim $(ii)$ above result in
the coincidence of the tilde-normalization of reduced
$\mathscr{P}_{\! \hbox{\tiny alg}}$ there
%,that is the tilde-normalization of
%
%$$\frac{q^2_{st}-1}{1+a_{st}^2 t_{st}}
%\,\sum_{l,m\ge 0}q_{st}^{2l}a_{st}^{2m}
%t_{st}^{m^2}\mathfrak{w}(\c_{\vec\rr,\vec\ss}^{[l\le l+m]})$$
with our 
$\h_{\,\vec\rr,\,\vec\ss}\,(\square\,;\,q,t,a)_{st}$;
see also formula (\ref{overlinep}) below.
Moreover, $h_{\hbox{\tiny odd}}=0$ (the van Straten-Warmt 
conjecture)  and
\begin{align}\label{bettit}
\h_{\,\vec\rr,\,\vec\ss}\,(\square\,;\,q=1,t,a=0)=
\hbox{\small$\sum$}_{j=0\,}^{\de} h_{2j}t^j, \hbox{\, for\, }
\de\!=\!\hbox{dim}\,J(\c_{\vec\rr,\vec\ss}).
\end{align}
\vskip -1.cm
\sq
\end{conjecture}
\smallskip

{\sf Comments.}
We hope that we will be able to justify Part $(i)$
in the Connection Conjecture following Proposition 2.3
from \cite{CJ}, where we used \cite{Ste} (for torus knots).
We note that the skein relations 
for the HOMFLYPT polynomials for iterated knots can be found
in \cite{AM,AMM} (see also \cite{Ma}). 
The case of $A_1$ was already checked:
\begin{align}\label{homjon2}
\h(m\om_1;\,q,q,-q^2)=
\tilde{\hbox{\small\em H\!O\!M}}(m\om_1;\,q,q^2)
=\tilde{J\!D}^{A_1}(m\om_1;\,q,q),
\end{align}
where the latter is the reduced tilde-normalized
Jones polynomial for any iterated knot 
$\t(\,\vec\rr,\,\vec\ss)$. We use
the explicit formulas from \cite{RJ,Mo}.
\smallskip

{\sf Khovanov-Rozansky polynomials.}
Due to $(ii)$, the conjectural relation of the
(uncolored) DAHA- superpolynomials 
to {\em Khovanov-Rozansky polynomials\,} for $sl_{n\!+1}$
is $a\mapsto\! t^{n\!+1}\sqrt{t/q}$ 
in DAHA parameters in the stable case (for sufficiently large
$n$); then one uses the substitution 
from (\ref{qtareli}). 
Also, the relation to the 
{\em Heegard- Floer homology\,} is expected
for $n\!+\!1\!=\!0$; we will not discuss here the
details. 

It is difficult to
calculate almost any
Khovanov-Rozansky polynomials for $sl_{n+1}$,
denoted below by $K\!h\!R_{n+1}$ (the reduced ones), 
unless $n=1$, i.e. in the case of the celebrated
Khovanov polynomials.
The polynomials $K\!h\!R_{\hbox{\tiny stab}}$\, are more
algebraic but not too much simpler.
%\smallskip

Even if the stable $K\!hR$\~polynomial is known,
the problem of recovering individual $K\!h\!R_{n+1}$ is quite
a challenge. Its theoretical solution is provided by the theory of 
{\em differentials\,} $\partial_{n+1}$; the corresponding homology,
$Ker(\partial_{n+1})/Im(\partial_{n+1})$ results in $K\!h\!R_{n+1}$
for any $n\ge 1$. These differentials are generally 
``transcendental". Their certain algebraic simplification,
suggested in \cite{DGR} and developed further in \cite{CJ},
works surprisingly well for small knots.
Among the torus knots, the ``smallest" counterexample for $n=1$ 
found in \cite{CJ} was $T(12,7)$.

We will denote the output of this procedure by
$D\!A\!H\!A' K\!h\!R_{n+1}$. The assumption in \cite{CJ} is 
that  $\partial_{n+1}$ are ``as surjective as possible" starting 
with small $n$; see Conjecture 2.7 and Section 3.6 there.
The reduction to $D\!A\!H\!A' K\!h\!R_{n+1}$ is
basically as follows. If $\h$ (in the DAHA-parameters) contains,
say, 
$q^{i}(t^j a^{m}+ t^{j-n-1}a^{m+1}+t^{j-2n-2}a^{m+2})$, where
$q^{i}t^i a^m$ was not involved in the previous
reductions, then $q^i t^{j-2n-2}a^{m+2}$ will be taken
from this triple to $D\!A\!H\!A' K\!h\!R_{n+1}$ (subject to
possible further reductions). However if the actual (topological)  
$\partial_{n+1}$ is not onto when acting from the space 
associated with  $q^i t^{j-n-1}a^{m+1}$ to that 
for $q^i t^{j}a^m$, then $q^{i}t^j a^{m}$ must be taken here. 

Therefore
$\tilde{K\!h\!R}_2-D\!A\!H\!A' K\!h\!R_{2}$ in the standard 
parameters can be expected a linear combination of 
the corrections $q_{st}^i t_{st}^j(1+t_{st})$ and 
$q_{st}^i t_{st}^j(1-t_{st}^2)$ for relatively simple knots. 
Such a difference, especially if it is with positive
coefficients, can be considered a confirmation of 
Part $(ii)$, which can be extended to   
{\em non-pseudo-algebraic\,}  iterated knots.

Changing $-1$ in $D\!A\!H\!A' K\!h\!R_2$ by 
$t_{st}^{-1}$ can be sufficient to 
obtain $\tilde{K\!h\!R}_2$ for small knots;
$C\!ab(7,2)T(4,3)$ is the simplest clear counterexample
we found.
We employ \,http://katlas.org/wiki/KnotTheory\, to calculate
Khovanov polynomials and use 
the presentations as in  Figure \ref{bngt}.

%\smallskip
\comment{
Here one can speculate that if the ``deepest" \, $K\!h\!R_2$\,
(counted from stable $n$) coincides with 
$D\!A\!H\!A' K\!h\!R_2$, then the procedure from \cite{CJ},
we are going to use, is likely to give correct answers for all 
``intermediate" $K\!h\!R_n$, but this we do not know (and do
not conjecture). 
}

%\medskip

\setcounter{equation}{0}
\section{\sc Numerical confirmations}

We will consider examples (mainly numerical) confirming
the Connection Conjecture. We selected only relatively
simple and instructional examples; however, some formulas are long.
Since they are expected to contain a lot
of geometric information (including the Betti numbers
of Jacobian factors), we believe that they are necessary in 
this paper. We do not consider torus knots here; see \cite{CJ,CJJ}.
Figure \ref{bngt} contains the braids for the main iterated
knots we will consider below.

%\vskip 1.5in
\begin{figure}[htbp]
\begin{center}
\vskip -0.0in
\hskip -0.5in
\includegraphics[scale=0.8]{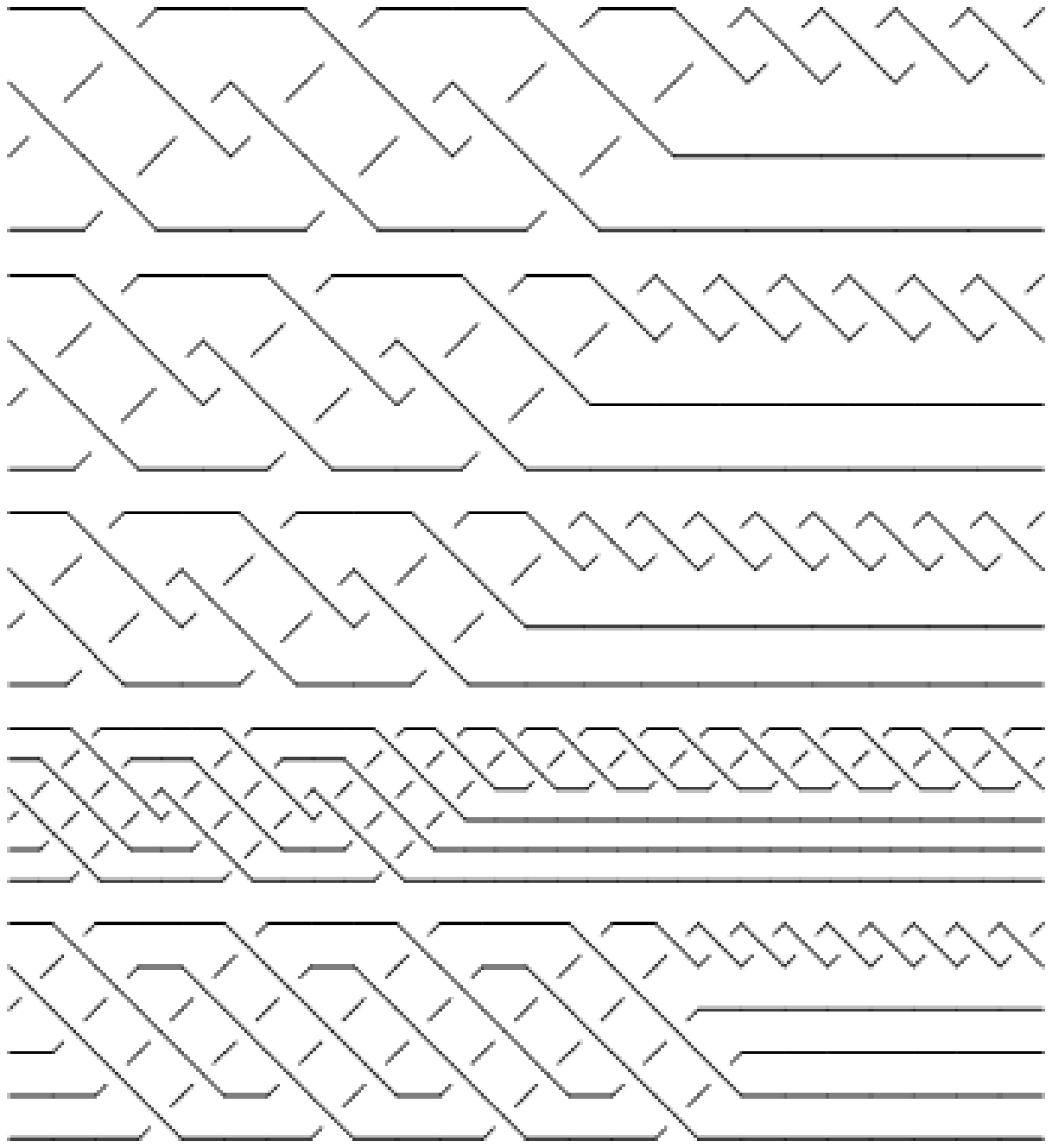}
\vskip +0.in
\caption{
{Cables (11,2),(13,2),(15,2) of T(3,2),}
\newline 
{Cab(19,3)T(3,2) and Cab(25,2)T(4,3)}}
\label{bngt}
\end{center}
\end{figure}
\medskip

\subsection{\bf Simplest algebraic iterations}
$$
(\mathbf 1):\ \vec\rr=\{3,2\},\,
\vec\ss=\{2,1\},\,
\t= C\!ab(13,2)T(3,2);\ 
\h_{\,\vec\rr,\,\vec\ss}\,(\square\,;\,\,q,t,a)=
$$
\renewcommand{\baselinestretch}{0.5} 
{\small
\(
1+q t+q^2 t+q^3 t+q^2 t^2+q^3 t^2+2 q^4 t^2+q^3 t^3+q^4 t^3+
2 q^5 t^3+q^4 t^4+q^5 t^4+2 q^6 t^4+q^5 t^5+q^6 t^5+q^7 t^5+
q^6 t^6+q^7 t^6+q^7 t^7+q^8 t^8+a^3 \bigl(q^6+q^7 t+q^8 t^2\bigr)+
a^2 \bigl(q^3+q^4+q^5+q^4 t+2 q^5 t+2 q^6 t+q^5 t^2+2 q^6 t^2
+2 q^7 t^2+q^6 t^3+2 q^7 t^3+q^8 t^3+
q^7 t^4+q^8 t^4+q^8 t^5\bigr)+
a \bigl(q+q^2+q^3+q^2 t+2 q^3 t+3 q^4 t+q^5 t+q^3 t^2+2 q^4 t^2+
4 q^5 t^2+q^6 t^2+q^4 t^3+2 q^5 t^3+4 q^6 t^3+q^7 t^3+q^5 t^4+
2 q^6 t^4+3 q^7 t^4+q^6 t^5+2 q^7 t^5+q^8 t^5+q^7 t^6+q^8 t^6
+q^8 t^7\bigr).
\)
}
\renewcommand{\baselinestretch}{1.2} 
\smallskip

{\sf Betti numbers.}
We will systematically provide the values of 
the DAHA superpolynomials at $a=0,q=1$; the corresponding
Betti numbers (conjecturally coinciding with their 
$t$\~coefficients) are known only in $2$ cases 
from the eight uncolored algebraic 
knots we consider. We will also give the formula for such a 
{\em DAHA-Bette polynomial\,} for $2\om_1$ and 
$\t_{\{3,2\},\{2,1\}}= 
C\!ab(13,2)T(3,2)$, which is  
the square of $\h(q\!=\!1,t,a\!=\!0)$ in (\ref{bett1}), that
for $\om_2$, and  calculate the {\em DAHA-Bette polynomials\,}
for non-algebraic knots. As for the latter, the coefficients
in (\ref{bett6},\ref{bett7}) 
are clear predecessors of the Betti numbers 
for $\,C\!ab(13,2)T(3,2)$, 
the first algebraic knot in the family $\,C\!ab(2\mm+1,2)T(3,2)\,$
for $\,\mm\in \Z$.

For $\vec\rr=\{3,2\},\, \vec\ss=\{2,1\}\,$, 
the DAHA-Bette polynomial equals
\begin{align}\label{bett1}
\h(q\!=\!1,t,a\!=\!0)=
1\!+\!3 t\!+\!4 t^2\!+\!4 t^3\!+\!4 t^4\!+\!3 t^5\!+\!2 t^6
\!+\!t^7\!+\!t^8,
\end{align}
which gives exactly the first line of
Betti numbers from the table before Conjecture 23 of
\cite{Pi}; the corresponding Puiseux exponents are
$(4,6,7)$. The next example gives the second
line in this table, for the exponents $(4,6,9)$; see
formula (\ref{bett2}) below. The corresponding two
$\c_{\vec\rr,\vec\ss}$ are represented by the local
rings $\r=\C[[t^4,t^6+t^{7,9}]]$;
$\de=8,9$. Note that {\em the arithmetic genus $\de$ 
always coincides with the $t$\~degree of the
Bette polynomial\,} (in the examples below).
Any other $\h(q\!=\!1,t,a\!=\!0)$ provided below are beyond 
known (or conjectured) formulas for the Betti numbers
of $\,J(\c_{\vec\rr,\vec\ss})$
(as far as we know). 
 
\medskip
{\sf HOMFLYPT Polynomial.}
It coincides with the  corresponding DAHA-generated 
tilde-normalized HOMFLYPT polynomial:
\begin{align*}
&\h(q,t\!\mapsto q,a\!\mapsto\!-a)=
1+q^2+q^3+2 q^4+q^5+3 q^6+q^7+3 q^8+q^9\\
&+3 q^{10}+q^{11}+2 q^{12}+q^{13}+q^{14}+q^{16}\,\, -a^3 
\bigl(q^6+q^8+q^{10}\bigr)\\
&+a^2 \bigl(q^3+q^4+2 q^5+2 q^6+3 q^7+2 q^8
+3 q^9+2 q^{10}+2 q^{11}+q^{12}+q^{13}\bigr)\\
&-a \bigl(q+q^2+2 q^3+2 q^4+4 q^5+3 q^6+5 q^7+3 q^8+5 q^9\\
&\ \ \ \ \ \ \ \ +3 q^{10}+4 q^{11}+2 q^{12}+2 q^{13}
+q^{14}+q^{15}\bigr).
\end{align*}
Up to $q^\bullet$, this polynomial is equal to the one considered 
in Section 7 of \cite{ObS} upon its division by $a^{16}$ 
and the substitution $z=q-q^{-1}$. We need to 
change $a^2,q^2$ there by our $a,q$, i.e. switch
from the {\em standard parameters\,} to our ones.

All superpolynomials we obtained
(we provide here only some), including non-algebraic
iterated knots, 
satisfy the HOMFLYPT part of the Connection
Conjecture \ref{CONCONJ}. 
The HOMFLYPT formulas will be omitted below. We used the 
software by S.~Artamonov for calculating colored 
HOMFLYPT polynomials, connected with papers \cite{AM,AMM}.
\smallskip

{\sf Khovanov-Rozansky polynomials.}
We will provide the corresponding DAHA-Khovanov 
polynomials for all examples below in {\em the
standard parameters\,}; for instance, $q,t$ in the next
two formulas are $q_{st},t_{st}$ from the conjecture.
If the algorithm from \cite{CJ} results in 
$D\!A\!H\!A'K\!h\!R_{2}$ matching the Khovanov polynomial
(for $sl_2$), then this is an indication that
$D\!A\!H\!A'K\!h\!R_{n+1}$ can match $\tilde{K\!h\!R}_{n+1}$ 
for all $n$.
 
For this and the next knot, we will provide 
$D\!A\!H\!A'K\!h\!R_3$, but will omit   
them in all other examples (they were calculated
for any $n$); $\tilde{K\!h\!R}_{\ge 3}$ are unknown for the 
knots we consider here. 
\smallskip

$$(1):\ D\!A\!H\!A'K\!h(C\!ab(13,2)T(3,2))=$$
\renewcommand{\baselinestretch}{0.5} 
{\small
\(
1+q^4 t^2+q^6 t^3+q^6 t^4+q^{10} t^5+q^8 t^6+q^{12} t^7
+2 q^{12} t^8+2 q^{14} t^9+q^{16} t^{10}
+2 q^{18} t^{11}+2 q^{20} t^{12}+q^{22} t^{13},
\)
}
\renewcommand{\baselinestretch}{1.2} 

$$\ \ \ \ D\!A\!H\!A'K\!h\!R_3(C\!ab(13,2)T(3,2))=$$
\renewcommand{\baselinestretch}{0.5} 
{\small
\(
1+q^4 t^2+q^8 t^3+q^6 t^4+q^8 t^4+q^{10} t^5
+q^{12} t^5+q^8 t^6+q^{10} t^6+2 q^{14} t^7
+q^{16} t^7+2 q^{12} t^8+q^{14} t^8+q^{18} t^8
+2 q^{16} t^9+2 q^{18} t^9+2 q^{16} t^{10}
+q^{22} t^{10}+3 q^{20} t^{11}+q^{22} t^{11}
+2 q^{20} t^{12}+q^{24} t^{12}+q^{26} t^{12}
+3 q^{24} t^{13}+q^{24} t^{14}
+q^{28} t^{14}+2 q^{28} t^{15}+q^{32} t^{16}.
\)
}
\renewcommand{\baselinestretch}{1.2} 

\smallskip

This knot, the next one, and the last
two knots (6,7 below) have their $DAHA$-Khovanov
polynomials coinciding with the actual
Khovanov polynomials
(under the tilde-normalization). 
\medskip

{\mathversion{bold}$\h_{\{3,2\},\{2,3\}}$.}
Let us switch from $C\!ab(13,2)$ to
the ``next" $C\!ab(15,2)$.
$$
(\mathbf 2):\ \vec\rr=\{3,2\},\,
\vec\ss=\{2,3\},\,
\t= C\!ab(15,2)T(3,2);\ 
\h_{\,\vec\rr,\,\vec\ss}\,(\square\,;\,\,q,t,a)=
$$
\renewcommand{\baselinestretch}{0.5} 
{\small
\(
1+q t+q^2 t+q^3 t+q^2 t^2+q^3 t^2+2 q^4 t^2+q^3 t^3+q^4 t^3
+2 q^5 t^3
+q^4 t^4+q^5 t^4+2 q^6 t^4+q^5 t^5+q^6 t^5+2 q^7 t^5+q^6 t^6+q^7 t^6
+q^8 t^6+q^7 t^7+q^8 t^7+q^8 t^8+q^9 t^9+a^3 \bigl(q^6+q^7 t+q^8 t^2
+q^9 t^3\bigr)+a^2 \bigl(q^3+q^4+q^5+q^4 t+2 q^5 t+2 q^6 t+q^5 t^2
+2 q^6 t^2+2 q^7 t^2+q^6 t^3+2 q^7 t^3+2 q^8 t^3+q^7 t^4+2 q^8 t^4
+q^9 t^4+q^8 t^5+q^9 t^5+q^9 t^6\bigr)+a \bigl(q+q^2+q^3+q^2 t
+2 q^3 t+3 q^4 t+q^5 t+q^3 t^2+2 q^4 t^2+4 q^5 t^2+q^6 t^2
+q^4 t^3+2 q^5 t^3+4 q^6 t^3+q^7 t^3+q^5 t^4+2 q^6 t^4
+4 q^7 t^4+q^8 t^4+q^6 t^5+2 q^7 t^5+3 q^8 t^5+q^7 t^6
+2 q^8 t^6+q^9 t^6+q^8 t^7+q^9 t^7+q^9 t^8\bigr).
\)
}
\renewcommand{\baselinestretch}{1.2} 
\smallskip

{\sf Betti Numbers.}
\begin{align} \label{bett2}
\h(q\!=\!1,t,a\!=\!0)\!=\!
1\!+\!3 t\!+\!4 t^2\!+\!4 t^3\!+\!4 t^4\!+\!4 t^5\!+\!3 t^6
\!+\!2 t^7\!+\!t^8\!+\!t^9.
\end{align}
\smallskip

{\sf Khovanov-Rozansky polynomials.}
$$(2):\ D\!A\!H\!A'K\!h(C\!ab(15,2)T(3,2))=$$
\renewcommand{\baselinestretch}{0.5} 
{\small
\(
1+q^4 t^2+q^6 t^3+q^6 t^4+q^{10} t^5
+q^8 t^6+q^{12} t^7+2 q^{12} t^8
+2 q^{14} t^9+q^{18} t^{11}+q^{20} t^{12},
\)
}
\renewcommand{\baselinestretch}{1.2} 
$$\ \ \ \ D\!A\!H\!A'K\!h\!R_3(C\!ab(15,2)T(3,2))=$$
\renewcommand{\baselinestretch}{0.5} 
{\small
\(
1+q^4 t^2+q^8 t^3+q^6 t^4+q^8 t^4+q^{10} t^5+q^{12} t^5
+q^8 t^6+q^{10} t^6+2 q^{14} t^7
+q^{16} t^7+2 q^{12} t^8+q^{14} t^8
+q^{18} t^8+2 q^{16} t^9+2 q^{18} t^9
+q^{16} t^{10}+q^{22} t^{10}+2 q^{20} t^{11}
+q^{22} t^{11}+q^{24} t^{12}+q^{26} t^{12}.
\)
}
\renewcommand{\baselinestretch}{1.2} 
\smallskip

\subsection{\bf Degree 5 and beyond}
The $a$\~degree was $3$ in the previous examples;
let us consider two examples of uncolored DAHA superpolynomials
of degree $5$ (which is $\ss_1\rr_2-1$ in these cases) and
then some further examples of DAHA-Betti polynomials, which
are beyond known or conjectured formulas. We note 
that $D\!A\!H\!A'K\!h\!R$ do not coincide with 
$\tilde{K\!h}$ in these two examples, but the difference
matches our conjecture.

$$
(\mathbf 3):\ \vec\rr=\{3,3\},\,
\vec\ss=\{2,1\},\,
\t= C\!ab(19,3)T(3,2);\ 
\h_{\,\vec\rr,\,\vec\ss}\,(\square\,;\,\,q,t,a)=
$$
\renewcommand{\baselinestretch}{0.5} 
{\small
\(
1+q t+q^2 t+q^3 t+q^4 t+q^5 t+q^2 t^2+q^3 t^2+2 q^4 t^2
+2 q^5 t^2+3 q^6 t^2+2 q^7 t^2+q^8 t^2+q^3 t^3+q^4 t^3
+2 q^5 t^3+3 q^6 t^3+4 q^7 t^3+4 q^8 t^3+4 q^9 t^3
+q^{10} t^3+q^4 t^4+q^5 t^4+2 q^6 t^4+3 q^7 t^4
+5 q^8 t^4+5 q^9 t^4+6 q^{10} t^4+4 q^{11} t^4
+q^{12} t^4+q^5 t^5+q^6 t^5+2 q^7 t^5+3 q^8 t^5
+5 q^9 t^5+6 q^{10} t^5+7 q^{11} t^5+6 q^{12} t^5
+3 q^{13} t^5+q^6 t^6+q^7 t^6+2 q^8 t^6+3 q^9 t^6
+5 q^{10} t^6+6 q^{11} t^6+8 q^{12} t^6+7 q^{13} t^6
+3 q^{14} t^6+q^{15} t^6+q^7 t^7+q^8 t^7+2 q^9 t^7
+3 q^{10} t^7+5 q^{11} t^7+6 q^{12} t^7+8 q^{13} t^7
+8 q^{14} t^7+3 q^{15} t^7+q^8 t^8+q^9 t^8+2 q^{10} t^8
+3 q^{11} t^8+5 q^{12} t^8+6 q^{13} t^8+8 q^{14} t^8
+7 q^{15} t^8+3 q^{16} t^8+q^9 t^9+q^{10} t^9
+2 q^{11} t^9+3 q^{12} t^9+5 q^{13} t^9+6 q^{14} t^9
+8 q^{15} t^9+6 q^{16} t^9+q^{17} t^9+q^{10} t^{10}
+q^{11} t^{10}+2 q^{12} t^{10}+3 q^{13} t^{10}
+5 q^{14} t^{10}+6 q^{15} t^{10}+7 q^{16} t^{10}
+4 q^{17} t^{10}+q^{11} t^{11}+q^{12} t^{11}
+2 q^{13} t^{11}+3 q^{14} t^{11}+5 q^{15} t^{11}
+6 q^{16} t^{11}+6 q^{17} t^{11}+q^{18} t^{11}
+q^{12} t^{12}+q^{13} t^{12}+2 q^{14} t^{12}
+3 q^{15} t^{12}+5 q^{16} t^{12}+5 q^{17} t^{12}
+4 q^{18} t^{12}+q^{13} t^{13}+q^{14} t^{13}
+2 q^{15} t^{13}+3 q^{16} t^{13}+5 q^{17} t^{13}
+4 q^{18} t^{13}+q^{19} t^{13}+q^{14} t^{14}
+q^{15} t^{14}+2 q^{16} t^{14}+3 q^{17} t^{14}
+4 q^{18} t^{14}+2 q^{19} t^{14}+q^{15} t^{15}
+q^{16} t^{15}+2 q^{17} t^{15}+3 q^{18} t^{15}
+3 q^{19} t^{15}+q^{16} t^{16}+q^{17} t^{16}
+2 q^{18} t^{16}+2 q^{19} t^{16}+q^{20} t^{16}
+q^{17} t^{17}+q^{18} t^{17}+2 q^{19} t^{17}
+q^{20} t^{17}+q^{18} t^{18}+q^{19} t^{18}
+q^{20} t^{18}+q^{19} t^{19}+q^{20} t^{19}
+q^{20} t^{20}+q^{21} t^{21}
\)

\vfil
\noindent
\(
+a^5 \bigl(q^{15}
+q^{16} t+q^{17} t+q^{17} t^2+q^{18} t^2
+q^{19} t^2+q^{18} t^3+q^{19} t^3+q^{19} t^4
+q^{20} t^4+q^{20} t^5+q^{21} t^6\bigr)
\)

\vfil
\noindent
\(
+a^4 \bigl(q^{10}+q^{11}+q^{12}+q^{13}+q^{14}
+q^{11} t+2 q^{12} t+3 q^{13} t+3 q^{14} t
+3 q^{15} t+q^{16} t+q^{12} t^2+2 q^{13} t^2
+4 q^{14} t^2+5 q^{15} t^2+5 q^{16} t^2+3 q^{17} t^2
+q^{18} t^2+q^{13} t^3+2 q^{14} t^3+4 q^{15} t^3
+6 q^{16} t^3+6 q^{17} t^3+3 q^{18} t^3+q^{19} t^3
+q^{14} t^4+2 q^{15} t^4+4 q^{16} t^4+6 q^{17} t^4
+6 q^{18} t^4+3 q^{19} t^4+q^{15} t^5+2 q^{16} t^5
+4 q^{17} t^5+6 q^{18} t^5+5 q^{19} t^5+q^{20} t^5
+q^{16} t^6+2 q^{17} t^6+4 q^{18} t^6+5 q^{19} t^6
+3 q^{20} t^6+q^{17} t^7+2 q^{18} t^7+4 q^{19} t^7
+3 q^{20} t^7+q^{21} t^7+q^{18} t^8+2 q^{19} t^8
+3 q^{20} t^8+q^{21} t^8+q^{19} t^9+2 q^{20} t^9
+q^{21} t^9+q^{20} t^{10}+q^{21} t^{10}+q^{21} t^{11}\bigr)
\)

\vfil
%\smallskip
\noindent
\(
+a^3 \bigl(q^6+q^7+2 q^8+2 q^9+2 q^{10}+q^{11}+q^{12}
+q^7 t+2 q^8 t+4 q^9 t+6 q^{10} t+7 q^{11} t+6 q^{12} t
+4 q^{13} t+2 q^{14} t+q^8 t^2+2 q^9 t^2+5 q^{10} t^2
+8 q^{11} t^2+12 q^{12} t^2+12 q^{13} t^2+10 q^{14} t^2
+5 q^{15} t^2+2 q^{16} t^2+q^9 t^3+2 q^{10} t^3
+5 q^{11} t^3+9 q^{12} t^3+14 q^{13} t^3+17 q^{14} t^3
+15 q^{15} t^3+9 q^{16} t^3+3 q^{17} t^3+q^{18} t^3
+q^{10} t^4+2 q^{11} t^4+5 q^{12} t^4+9 q^{13} t^4
+15 q^{14} t^4+19 q^{15} t^4+18 q^{16} t^4+10 q^{17} t^4
+3 q^{18} t^4+q^{11} t^5+2 q^{12} t^5+5 q^{13} t^5
+9 q^{14} t^5+15 q^{15} t^5+20 q^{16} t^5+18 q^{17} t^5
+9 q^{18} t^5+2 q^{19} t^5+q^{12} t^6+2 q^{13} t^6
+5 q^{14} t^6+9 q^{15} t^6+15 q^{16} t^6+19 q^{17} t^6
+15 q^{18} t^6+5 q^{19} t^6+q^{13} t^7+2 q^{14} t^7
+5 q^{15} t^7+9 q^{16} t^7+15 q^{17} t^7+17 q^{18} t^7
+10 q^{19} t^7+2 q^{20} t^7+q^{14} t^8+2 q^{15} t^8
+5 q^{16} t^8+9 q^{17} t^8+14 q^{18} t^8+12 q^{19} t^8
+4 q^{20} t^8+q^{15} t^9+2 q^{16} t^9+5 q^{17} t^9
+9 q^{18} t^9+12 q^{19} t^9+6 q^{20} t^9+q^{21} t^9
+q^{16} t^{10}+2 q^{17} t^{10}+5 q^{18} t^{10}
+8 q^{19} t^{10}+7 q^{20} t^{10}+q^{21} t^{10}
+q^{17} t^{11}+2 q^{18} t^{11}+5 q^{19} t^{11}
+6 q^{20} t^{11}+2 q^{21} t^{11}+q^{18} t^{12}
+2 q^{19} t^{12}+4 q^{20} t^{12}+2 q^{21} t^{12}
+q^{19} t^{13}+2 q^{20} t^{13}+2 q^{21} t^{13}
+q^{20} t^{14}+q^{21} t^{14}+q^{21} t^{15}\bigr)
\)

\vfil
%\smallskip
\noindent
\(
+a^2 \bigl(q^3+q^4+2 q^5+2 q^6+2 q^7+q^8+q^9+q^4 t
+2 q^5 t+4 q^6 t+6 q^7 t+8 q^8 t+7 q^9 t+6 q^{10} t
+3 q^{11} t+q^{12} t+q^5 t^2+2 q^6 t^2+5 q^7 t^2
+8 q^8 t^2+13 q^9 t^2+15 q^{10} t^2+15 q^{11} t^2
+10 q^{12} t^2+5 q^{13} t^2+q^{14} t^2+q^6 t^3
+2 q^7 t^3+5 q^8 t^3+9 q^9 t^3+15 q^{10} t^3
+20 q^{11} t^3+24 q^{12} t^3+19 q^{13} t^3
+11 q^{14} t^3+4 q^{15} t^3+q^{16} t^3+q^7 t^4+2 q^8 t^4
+5 q^9 t^4+9 q^{10} t^4+16 q^{11} t^4+22 q^{12} t^4
+29 q^{13} t^4+27 q^{14} t^4+16 q^{15} t^4+6 q^{16} t^4
+q^{17} t^4+q^8 t^5+2 q^9 t^5+5 q^{10} t^5+9 q^{11} t^5
+16 q^{12} t^5+23 q^{13} t^5+31 q^{14} t^5+30 q^{15} t^5
+19 q^{16} t^5+6 q^{17} t^5+q^{18} t^5+q^9 t^6
+2 q^{10} t^6+5 q^{11} t^6+9 q^{12} t^6+16 q^{13} t^6
+23 q^{14} t^6+32 q^{15} t^6+30 q^{16} t^6+16 q^{17} t^6
+4 q^{18} t^6+q^{10} t^7+2 q^{11} t^7+5 q^{12} t^7
+9 q^{13} t^7+16 q^{14} t^7+23 q^{15} t^7+31 q^{16} t^7
+27 q^{17} t^7+11 q^{18} t^7+q^{19} t^7+q^{11} t^8
+2 q^{12} t^8+5 q^{13} t^8+9 q^{14} t^8+16 q^{15} t^8
+23 q^{16} t^8+29 q^{17} t^8+19 q^{18} t^8+5 q^{19} t^8
+q^{12} t^9+2 q^{13} t^9+5 q^{14} t^9+9 q^{15} t^9
+16 q^{16} t^9+22 q^{17} t^9+24 q^{18} t^9+10 q^{19} t^9
+q^{20} t^9+q^{13} t^{10}+2 q^{14} t^{10}+5 q^{15} t^{10}
+9 q^{16} t^{10}+16 q^{17} t^{10}+20 q^{18} t^{10}
+15 q^{19} t^{10}+3 q^{20} t^{10}+q^{14} t^{11}
+2 q^{15} t^{11}+5 q^{16} t^{11}+9 q^{17} t^{11}
+15 q^{18} t^{11}+15 q^{19} t^{11}+6 q^{20} t^{11}
+q^{15} t^{12}+2 q^{16} t^{12}+5 q^{17} t^{12}
+9 q^{18} t^{12}+13 q^{19} t^{12}+7 q^{20} t^{12}
+q^{21} t^{12}+q^{16} t^{13}+2 q^{17} t^{13}+5 q^{18} t^{13}
+8 q^{19} t^{13}+8 q^{20} t^{13}+q^{21} t^{13}+q^{17} t^{14}
+2 q^{18} t^{14}+5 q^{19} t^{14}+6 q^{20} t^{14}
+2 q^{21} t^{14}+q^{18} t^{15}+2 q^{19} t^{15}
+4 q^{20} t^{15}+2 q^{21} t^{15}+q^{19} t^{16}
+2 q^{20} t^{16}+2 q^{21} t^{16}+q^{20} t^{17}
+q^{21} t^{17}+q^{21} t^{18}\bigr)
\)

\vfil
%\smallskip
\noindent
\(
+a \bigl(q+q^2+q^3+q^4+q^5+q^2 t+2 q^3 t+3 q^4 t
+4 q^5 t+5 q^6 t+4 q^7 t+2 q^8 t+q^9 t+q^3 t^2
+2 q^4 t^2+4 q^5 t^2+6 q^6 t^2+9 q^7 t^2+10 q^8 t^2
+9 q^9 t^2+5 q^{10} t^2+2 q^{11} t^2+q^4 t^3
+2 q^5 t^3+4 q^6 t^3+7 q^7 t^3+11 q^8 t^3+14 q^9 t^3
+16 q^{10} t^3+13 q^{11} t^3+6 q^{12} t^3+2 q^{13} t^3
+q^5 t^4+2 q^6 t^4+4 q^7 t^4+7 q^8 t^4+12 q^9 t^4
+16 q^{10} t^4+20 q^{11} t^4+20 q^{12} t^4+13 q^{13} t^4
+4 q^{14} t^4+q^{15} t^4+q^6 t^5+2 q^7 t^5+4 q^8 t^5
+7 q^9 t^5+12 q^{10} t^5+17 q^{11} t^5+22 q^{12} t^5
+24 q^{13} t^5+17 q^{14} t^5+7 q^{15} t^5+q^{16} t^5
+q^7 t^6+2 q^8 t^6+4 q^9 t^6+7 q^{10} t^6+12 q^{11} t^6
+17 q^{12} t^6+23 q^{13} t^6+26 q^{14} t^6+18 q^{15} t^6
+7 q^{16} t^6+q^{17} t^6+q^8 t^7+2 q^9 t^7+4 q^{10} t^7
+7 q^{11} t^7+12 q^{12} t^7+17 q^{13} t^7+23 q^{14} t^7
+26 q^{15} t^7+17 q^{16} t^7+4 q^{17} t^7+q^9 t^8
+2 q^{10} t^8+4 q^{11} t^8+7 q^{12} t^8+12 q^{13} t^8
+17 q^{14} t^8+23 q^{15} t^8+24 q^{16} t^8+13 q^{17} t^8
+2 q^{18} t^8+q^{10} t^9+2 q^{11} t^9+4 q^{12} t^9
+7 q^{13} t^9+12 q^{14} t^9+17 q^{15} t^9+22 q^{16} t^9
+20 q^{17} t^9+6 q^{18} t^9+q^{11} t^{10}+2 q^{12} t^{10}
+4 q^{13} t^{10}+7 q^{14} t^{10}+12 q^{15} t^{10}
+17 q^{16} t^{10}+20 q^{17} t^{10}+13 q^{18} t^{10}
+2 q^{19} t^{10}+q^{12} t^{11}+2 q^{13} t^{11}
+4 q^{14} t^{11}+7 q^{15} t^{11}+12 q^{16} t^{11}
+16 q^{17} t^{11}+16 q^{18} t^{11}+5 q^{19} t^{11}
+q^{13} t^{12}+2 q^{14} t^{12}+4 q^{15} t^{12}
+7 q^{16} t^{12}+12 q^{17} t^{12}+14 q^{18} t^{12}
+9 q^{19} t^{12}+q^{20} t^{12}+q^{14} t^{13}
+2 q^{15} t^{13}+4 q^{16} t^{13}+7 q^{17} t^{13}
+11 q^{18} t^{13}+10 q^{19} t^{13}+2 q^{20} t^{13}
+q^{15} t^{14}+2 q^{16} t^{14}+4 q^{17} t^{14}
+7 q^{18} t^{14}+9 q^{19} t^{14}+4 q^{20} t^{14}
+q^{16} t^{15}+2 q^{17} t^{15}+4 q^{18} t^{15}
+6 q^{19} t^{15}+5 q^{20} t^{15}+q^{17} t^{16}
+2 q^{18} t^{16}+4 q^{19} t^{16}+4 q^{20} t^{16}
+q^{21} t^{16}+q^{18} t^{17}+2 q^{19} t^{17}
+3 q^{20} t^{17}+q^{21} t^{17}+q^{19} t^{18}
+2 q^{20} t^{18}+q^{21} t^{18}+q^{20} t^{19}
+q^{21} t^{19}+q^{21} t^{20}\bigr).
\)
}
\renewcommand{\baselinestretch}{1.2} 
\vfil
\smallskip

{\sf  Betti numbers} (conjecturally):
\begin{align} \label{bett3}
&\h(q\!=\!1,t,a\!=\!0)=
1 + 5 t + 12 t^2 + 20 t^3 + 28 t^4 + 34 t^5 
+ 37 t^6 \\
&+ 37 t^7
+ 36 t^8 + 33 t^9 + 29 t^{10}
+ 25 t^{11} + 21 t^{12} 
+ 17 t^{13} + 13 t^{14}\notag\\
&+ 10 t^{15} + 7 t^{16} + 5 t^{17} + 3 t^{18} + 2 t^{19} 
+ t^{20} + t^{21}.\notag
\end{align}

The germ $\c_{\vec\rr,\vec\ss}$ is 
given here by the local ring $\r=\C[[t^6,t^{9}+t^{10}]]$ with 
the {\em Puiseux exponents\,} $(6,9,10)$, the valuation
semigroup $\Gamma$ generated by $6,9,19$
and  $\de=|\N\setminus \Ga|=21$,
which matches the top $t$\~degree in the
DAHA-Betti polynomial from (\ref{bett3}). Thus
the Euler number of $J(\c_{\vec\rr,\vec\ss})$ of
$\c_{\vec\rr,\vec\ss}$  is $365$ (the value in 
(\ref{bett3}) at $t=1$), which follows from
%from the reduction of the superpolynomial
%above to the HOMFLYPT one (under $a\mapsto -a,t=q$) and 
\cite{ObS,Ma}, and we predict that the Betti 
numbers of $J(\c_{\vec\rr,\vec\ss})$ are the coefficients of
this polynomial. This example is from the table after
Theorem 13 in \cite{Pi} listed some cases where
the approach there is not applicable (including counting
the Euler number).
\smallskip

{\sf Khovanov polynomial:}
$$(3):\ D\!A\!H\!A'K\!h(C\!ab(19,3)T(3,2))=$$
\renewcommand{\baselinestretch}{0.5} 
{\small
\(
1+q^4 t^2+q^6 t^3+q^6 t^4+q^{10} t^5+q^8 t^6+q^{12} t^7
+q^{10} t^8+q^{12} t^8+2 q^{14} t^9+q^{12} t^{10}+q^{14} t^{10}
+2 q^{16} t^{11}+q^{18} t^{11}+3 q^{16} t^{12}+q^{20} t^{12}
+2 q^{18} t^{13}+2 q^{20} t^{13}+2 q^{18} t^{14}+q^{22} t^{14}
+4 q^{22} t^{15}+q^{20} t^{16}+2 q^{22} t^{16}+q^{24} t^{16}
+q^{26} t^{16}+4 q^{24} t^{17}+3 q^{24} t^{18}+q^{28} t^{18}
+3 q^{26} t^{19}+2 q^{28} t^{19}+q^{26} t^{20}+q^{28} t^{20}
+2 q^{30} t^{20}+4 q^{30} t^{21}+2 q^{30} t^{22}+2 q^{32} t^{22}
+3 q^{32} t^{23}+q^{32} t^{24}+q^{36} t^{24}+2 q^{36} t^{25}
+q^{34} t^{26}+q^{36} t^{26}+2 q^{38} t^{26}+q^{36} t^{27}
+q^{38} t^{27}+q^{40} t^{27}+q^{40} t^{29}+q^{42} t^{30}.
\)
}
\renewcommand{\baselinestretch}{1.2} 

This one is now different from the actual Khovanov polynomial:
\begin{align*}
&\tilde{K\!h}-D\!A\!H\!A'K\!h=
q^{34} t^{23} \left(1+2 t-t^3\right)\\
&+q^{36} t^{22} 
\left(1-t^2\right)+q^{40} t^{27} 
\left(1-t^2\right)+q^{42} t^{28} \left(1-t^2\right).
\end{align*}

There is no coincidence, but
the $t$\~degree corrections involve only top terms and 
all are in the expected direction. As it was discussed
after the Connection Conjecture, if the difference
$\tilde{K\!h}-D\!A\!H\!A'K\!h$ is a sum of the terms
$q^i t^j(1+t)$ and 
$q^i t^j(1-t^2)$ with positive coefficients,
then this is a confirmation of Part $(ii)$.
Here and below (where Khovanov polynomials
are discussed) $q,t$ are actually $q_{st},t_{st}$. 
This is what can be expected, taking into consideration the 
nature of the procedure from \cite{CJ} 
(and the example of $T(12,7)$ there). 
\medskip

{\mathversion{bold}$\h_{\{4,2\},\{3,1\}}$.}
The next example will be of the same deg$_a=5$.
Recall that deg${}_a\le (\ss_1\rr_2-1)$ 
in the uncolored case with $2$ iterations; conjecturally
they coincide for algebraic knots. Generally, this
must be multiplied by \,ord$(\la_b)$.
$$
(\mathbf 4):\ \vec\rr=\{4,2\},\, \vec\ss=\{3,1\},\,
\t= C\!ab(25,2)T(4,3);\ 
\h_{\,\vec\rr,\,\vec\ss}\,(\square\,;\,\,q,t,a)=
$$

\vfil
\renewcommand{\baselinestretch}{0.5} 
{\small
\(
1+q t+q^2 t+q^3 t+q^4 t+q^5 t+q^2 t^2+q^3 t^2+2 q^4 t^2
+2 q^5 t^2+3 q^6 t^2+q^7 t^2+q^8 t^2+q^3 t^3+q^4 t^3
+2 q^5 t^3+3 q^6 t^3+4 q^7 t^3+3 q^8 t^3+3 q^9 t^3
+q^4 t^4+q^5 t^4+2 q^6 t^4+3 q^7 t^4+5 q^8 t^4+4 q^9 t^4
+5 q^{10} t^4+q^{11} t^4+q^5 t^5+q^6 t^5+2 q^7 t^5
+3 q^8 t^5+5 q^9 t^5+5 q^{10} t^5+6 q^{11} t^5
+q^{12} t^5+q^6 t^6+q^7 t^6+2 q^8 t^6+3 q^9 t^6
+5 q^{10} t^6+5 q^{11} t^6+7 q^{12} t^6+q^{13} t^6
+q^7 t^7+q^8 t^7+2 q^9 t^7+3 q^{10} t^7+5 q^{11} t^7
+5 q^{12} t^7+6 q^{13} t^7+q^{14} t^7+q^8 t^8+q^9 t^8
+2 q^{10} t^8+3 q^{11} t^8+5 q^{12} t^8+5 q^{13} t^8
+5 q^{14} t^8+q^9 t^9+q^{10} t^9+2 q^{11} t^9+3 q^{12} t^9
+5 q^{13} t^9+4 q^{14} t^9+3 q^{15} t^9+q^{10} t^{10}
+q^{11} t^{10}+2 q^{12} t^{10}+3 q^{13} t^{10}
+5 q^{14} t^{10}+3 q^{15} t^{10}+q^{16} t^{10}
+q^{11} t^{11}+q^{12} t^{11}+2 q^{13} t^{11}
+3 q^{14} t^{11}+4 q^{15} t^{11}+q^{16} t^{11}
+q^{12} t^{12}+q^{13} t^{12}+2 q^{14} t^{12}
+3 q^{15} t^{12}+3 q^{16} t^{12}+q^{13} t^{13}
+q^{14} t^{13}+2 q^{15} t^{13}+2 q^{16} t^{13}
+q^{17} t^{13}+q^{14} t^{14}+q^{15} t^{14}
+2 q^{16} t^{14}+q^{17} t^{14}+q^{15} t^{15}
+q^{16} t^{15}+q^{17} t^{15}+q^{16} t^{16}
+q^{17} t^{16}+q^{17} t^{17}+q^{18} t^{18}
\)

\vfil
\noindent
\(+a^5 \bigl(q^{15}+q^{16} t+q^{17} t^2+q^{18} t^3\bigr)
\)

\vfil
\noindent
\(
+a^4 \bigl(q^{10}+q^{11}+q^{12}+q^{13}+q^{14}
+q^{11} t+2 q^{12} t+2 q^{13} t+2 q^{14} t
+2 q^{15} t+q^{12} t^2+2 q^{13} t^2+3 q^{14} t^2
+3 q^{15} t^2+2 q^{16} t^2+q^{13} t^3+2 q^{14} t^3
+3 q^{15} t^3+3 q^{16} t^3+2 q^{17} t^3+q^{14} t^4
+2 q^{15} t^4+3 q^{16} t^4+2 q^{17} t^4+q^{18} t^4
+q^{15} t^5+2 q^{16} t^5+2 q^{17} t^5+q^{18} t^5
+q^{16} t^6+2 q^{17} t^6+q^{18} t^6+q^{17} t^7
+q^{18} t^7+q^{18} t^8\bigr)
\)

\vfil
\noindent
\(
+a^3 \bigl(q^6+q^7
+2 q^8+2 q^9+2 q^{10}+q^{11}+q^{12}+q^7 t+2 q^8 t
+4 q^9 t+5 q^{10} t+6 q^{11} t+4 q^{12} t+3 q^{13} t
+q^{14} t+q^8 t^2+2 q^9 t^2+5 q^{10} t^2
+7 q^{11} t^2+10 q^{12} t^2+7 q^{13} t^2
+5 q^{14} t^2+q^{15} t^2+q^9 t^3+2 q^{10} t^3
+5 q^{11} t^3+8 q^{12} t^3+12 q^{13} t^3
+9 q^{14} t^3+6 q^{15} t^3+q^{16} t^3+q^{10} t^4
+2 q^{11} t^4+5 q^{12} t^4+8 q^{13} t^4+13 q^{14} t^4
+9 q^{15} t^4+5 q^{16} t^4+q^{17} t^4+q^{11} t^5
+2 q^{12} t^5+5 q^{13} t^5+8 q^{14} t^5+12 q^{15} t^5
+7 q^{16} t^5+3 q^{17} t^5+q^{12} t^6+2 q^{13} t^6
+5 q^{14} t^6+8 q^{15} t^6+10 q^{16} t^6+4 q^{17} t^6
+q^{18} t^6+q^{13} t^7+2 q^{14} t^7+5 q^{15} t^7
+7 q^{16} t^7+6 q^{17} t^7+q^{18} t^7+q^{14} t^8
+2 q^{15} t^8+5 q^{16} t^8+5 q^{17} t^8+2 q^{18} t^8
+q^{15} t^9+2 q^{16} t^9+4 q^{17} t^9+2 q^{18} t^9
+q^{16} t^{10}+2 q^{17} t^{10}+2 q^{18} t^{10}
+q^{17} t^{11}+q^{18} t^{11}+q^{18} t^{12}\bigr)
\)

\vfil
\noindent
\(
+a^2 \bigl(q^3+q^4+2 q^5+2 q^6+2 q^7+q^8+q^9+q^4 t
+2 q^5 t+4 q^6 t+6 q^7 t+7 q^8 t+6 q^9 t+5 q^{10} t
+2 q^{11} t+q^{12} t+q^5 t^2+2 q^6 t^2+5 q^7 t^2
+8 q^8 t^2+12 q^9 t^2+12 q^{10} t^2+11 q^{11} t^2
+5 q^{12} t^2+2 q^{13} t^2+q^6 t^3+2 q^7 t^3+5 q^8 t^3
+9 q^9 t^3+14 q^{10} t^3+17 q^{11} t^3+17 q^{12} t^3
+8 q^{13} t^3+3 q^{14} t^3+q^7 t^4+2 q^8 t^4+5 q^9 t^4
+9 q^{10} t^4+15 q^{11} t^4+19 q^{12} t^4+20 q^{13} t^4
+9 q^{14} t^4+3 q^{15} t^4+q^8 t^5+2 q^9 t^5+5 q^{10} t^5
+9 q^{11} t^5+15 q^{12} t^5+20 q^{13} t^5+20 q^{14} t^5
+8 q^{15} t^5+2 q^{16} t^5+q^9 t^6+2 q^{10} t^6
+5 q^{11} t^6+9 q^{12} t^6+15 q^{13} t^6+19 q^{14} t^6
+17 q^{15} t^6+5 q^{16} t^6+q^{17} t^6+q^{10} t^7
+2 q^{11} t^7+5 q^{12} t^7+9 q^{13} t^7+15 q^{14} t^7
+17 q^{15} t^7+11 q^{16} t^7+2 q^{17} t^7+q^{11} t^8
+2 q^{12} t^8+5 q^{13} t^8+9 q^{14} t^8+14 q^{15} t^8
+12 q^{16} t^8+5 q^{17} t^8+q^{12} t^9+2 q^{13} t^9
+5 q^{14} t^9+9 q^{15} t^9+12 q^{16} t^9+6 q^{17} t^9
+q^{18} t^9+q^{13} t^{10}+2 q^{14} t^{10}+5 q^{15} t^{10}
+8 q^{16} t^{10}+7 q^{17} t^{10}+q^{18} t^{10}
+q^{14} t^{11}+2 q^{15} t^{11}+5 q^{16} t^{11}
+6 q^{17} t^{11}+2 q^{18} t^{11}+q^{15} t^{12}
+2 q^{16} t^{12}+4 q^{17} t^{12}+2 q^{18} t^{12}
+q^{16} t^{13}+2 q^{17} t^{13}+2 q^{18} t^{13}
+q^{17} t^{14}+q^{18} t^{14}+q^{18} t^{15}\bigr)
\)

\vfil
\noindent
\(
+a \bigl(q+q^2+q^3+q^4+q^5+q^2 t+2 q^3 t+3 q^4 t
+4 q^5 t+5 q^6 t+3 q^7 t+2 q^8 t+q^9 t+q^3 t^2
+2 q^4 t^2+4 q^5 t^2+6 q^6 t^2+9 q^7 t^2+8 q^8 t^2
+7 q^9 t^2+3 q^{10} t^2+q^{11} t^2+q^4 t^3+2 q^5 t^3
+4 q^6 t^3+7 q^7 t^3+11 q^8 t^3+12 q^9 t^3
+13 q^{10} t^3+7 q^{11} t^3+2 q^{12} t^3
+q^5 t^4+2 q^6 t^4+4 q^7 t^4+7 q^8 t^4
+12 q^9 t^4+14 q^{10} t^4+17 q^{11} t^4
+10 q^{12} t^4+3 q^{13} t^4+q^6 t^5+2 q^7 t^5
+4 q^8 t^5+7 q^9 t^5+12 q^{10} t^5+15 q^{11} t^5
+19 q^{12} t^5+11 q^{13} t^5+3 q^{14} t^5
+q^7 t^6+2 q^8 t^6+4 q^9 t^6+7 q^{10} t^6
+12 q^{11} t^6+15 q^{12} t^6+19 q^{13} t^6
+10 q^{14} t^6+2 q^{15} t^6+q^8 t^7+2 q^9 t^7
+4 q^{10} t^7+7 q^{11} t^7+12 q^{12} t^7
+15 q^{13} t^7+17 q^{14} t^7+7 q^{15} t^7
+q^{16} t^7+q^9 t^8+2 q^{10} t^8+4 q^{11} t^8
+7 q^{12} t^8+12 q^{13} t^8+14 q^{14} t^8
+13 q^{15} t^8+3 q^{16} t^8+q^{10} t^9
+2 q^{11} t^9+4 q^{12} t^9+7 q^{13} t^9
+12 q^{14} t^9+12 q^{15} t^9+7 q^{16} t^9
+q^{17} t^9+q^{11} t^{10}+2 q^{12} t^{10}
+4 q^{13} t^{10}+7 q^{14} t^{10}+11 q^{15} t^{10}
+8 q^{16} t^{10}+2 q^{17} t^{10}+q^{12} t^{11}
+2 q^{13} t^{11}+4 q^{14} t^{11}+7 q^{15} t^{11}
+9 q^{16} t^{11}+3 q^{17} t^{11}+q^{13} t^{12}
+2 q^{14} t^{12}+4 q^{15} t^{12}+6 q^{16} t^{12}
+5 q^{17} t^{12}+q^{14} t^{13}+2 q^{15} t^{13}
+4 q^{16} t^{13}+4 q^{17} t^{13}+q^{18} t^{13}
+q^{15} t^{14}+2 q^{16} t^{14}+3 q^{17} t^{14}
+q^{18} t^{14}+q^{16} t^{15}+2 q^{17} t^{15}
+q^{18} t^{15}+q^{17} t^{16}+q^{18} t^{16}
+q^{18} t^{17}\bigr).
\)
}
\renewcommand{\baselinestretch}{1.2} 
\smallskip

{\sf Betti numbers} (conjecturally):
\begin{align} \label{bett4}
&\h(q\!=\!1,t,a\!=\!0)\!=\!
1\!+5 t\!+\!11 t^2\!+\!17 t^3\!+\!22 t^4\!+\!24 t^5\!+
\!25 t^6\!+\!24 t^7\\
+22 t^8\!&+\!19 t^9\!+\!16 t^{10}
\!+\!12 t^{11}\!+\!10 t^{12}
\!+\!7 t^{13}\!+\!5 t^{14}\!+\!3 t^{15}\!+\!2 t^{16}
\!+\!t^{17}\!+\!t^{18}.\notag
\end{align}

The corresponding (germ of) curve  in this case is given
by the ring $\r=\C[[t^6,t^{8}+t^9]]$ with 
the Puiseux exponents $(6,8,9)$ and the valuation
semigroup $\Gamma$ generated by $6,8,25$. Accordingly, 
%Generally,$\Ga=\lan 2d,2u,2(d-1)u+v\ran$ for the exponents 
%$(2d,2u,v)$ subject to $2d<2u<v$ and gcd$(v,2)=1$;
%see \cite{Pi} and references there.
$\de=|\N\setminus \Ga|=18$ in the considered case, 
which matches the top $t$\~degree in (\ref{bett4}). 

Furthermore, $\h(q\!=\!1,t=1,a\!=\!0)\!=\! 227$, which does
coincide with the specialization at $v=9$ of the general formula 
for the Euler number 
$e(J(\c_{\vec\rr,\vec\ss}))=229/2+25 v/2$ in the case of the 
Puiseux exponents $(6,8,v)$. This formula is from the 
Main Theorem in \cite{Pi}; its proof was involved technically and
was not extended to the Betti numbers. We 
conjecture that they are the $t$\~coefficients in (\ref{bett4})
in this case ($v=9$).
\smallskip

{\sf Khovanov polynomial:}
$$(4):\ D\!A\!H\!A'K\!h(C\!ab(25,2)T(4,3))=$$
\renewcommand{\baselinestretch}{0.5} 
{\small
\(
1+q^4 t^2+q^6 t^3+q^6 t^4+q^{10} t^5+q^8 t^6
+q^{12} t^7+q^{10} t^8+q^{12} t^8+2 q^{14} t^9
+q^{12} t^{10}+q^{14} t^{10}+2 q^{16} t^{11}
+q^{18} t^{11}+3 q^{16} t^{12}+q^{20} t^{12}
+2 q^{18} t^{13}+2 q^{20} t^{13}+q^{18} t^{14}
+q^{20} t^{14}+q^{22} t^{14}+4 q^{22} t^{15}
+q^{20} t^{16}+q^{22} t^{16}+q^{24} t^{16}
+q^{26} t^{16}+2 q^{24} t^{17}+q^{26} t^{17}
+2 q^{24} t^{18}+q^{28} t^{18}+2 q^{26} t^{19}
+q^{28} t^{19}+2 q^{28} t^{20}+q^{30} t^{20}
+3 q^{30} t^{21}+q^{32} t^{22}+q^{34} t^{23}
+q^{32} t^{24}+q^{36} t^{24}
+q^{34} t^{25}+q^{36} t^{26}+q^{38} t^{27}.
\)
}
\renewcommand{\baselinestretch}{1.2} 
\smallskip

This one is different from the corresponding
Khovanov polynomial as well as in the previous example.
Namely,
\begin{align*}
&\tilde{K\!h}-D\!A\!H\!A'K\!h=
q^{32} t^{22} \left(1-t^2\right)\\
&+q^{34} t^{23} 
\left(1-t^2\right)+q^{36} t^{24} 
\left(1-t^2\right)+q^{38} t^{25} \left(1-t^2\right).
\end{align*}
Once again, the $t$\~degrees
increase from $\tilde{K\!h}$ to $D\!A\!H\!A'K\!h$;
also, the number of $q,t$\~monomials
is the same in both, but this seems accidental. 
In spite of the nonzero difference, this 
confirms Conjecture \ref{CONCONJ}. 
\smallskip

{\sf Further Betti numbers.}
First of all, our DAHA-Betti polynomials (as far as we
calculated them) fully match
the table before Conjecture 23 in \cite{Pi}. For instance,
the coefficients of $\h_{\{5,2\},\{2,1\}}(q\!=\!1,t,a\!=\!0)$ 
for $C\!ab(21,2)T(5,2)$ 
constitute the line there for the Puiseux exponents
$(4,10,11)$. Without 
posting the (known) formulas for the superpolynomials, let us
provide below only the DAHA-Betti polynomials,
which conjecturally give the corresponding Betti numbers.
\smallskip

{\mathversion{bold}$\h_{\{4,2\},\{3,1\}}$.} In this example,
the knot, the ring of the singularity, 
the valuation semigroup and the arithmetic genus
are as follows:
$$ 
\t=C\!ab(20,3)T(3,2),\  \r=C[[t^6,t^9+t^{11}]],\
\Gamma=\lan 6,9,20\ran,\ \de=22. 
$$
The counterpart of the DAHA-Betti polynomial (\ref{bett3}) reads:
\begin{align} \label{bett3x}
&\h(q\!=\!1,t,a\!=\!0)=
1 + 5 t + 12 t^2 + 20 t^3 + 28 t^4 + 35 t^5 + 39 t^6\\
&+ 40 t^7 + 
39 t^8 +  37 t^9 + 33 t^{10} + 29 t^{11} + 25 t^{12} + 21 t^{13} 
+ 17 t^{14}\notag\\
&+ 13 t^{15} + 
 10 t^{16} + 7 t^{17} + 5 t^{18} + 3 t^{19} + 2 t^{20} + 
t^{21} + t^{22}.\notag
\end{align}
The arithmetic genus of $\c=\c_{\{3,3\},\{2,2\}}$
equals $\de=22$ and $e(J(\c))=423$ in this case; the 
Puiseux exponents are $(6,9,11)$, deg$_a$=5.
\medskip

{\mathversion{bold}$\h_{\{5,2\},\{3,2\}}$.} The Puiseux exponents
are $(6,15,17)$, deg$_a=5$,
$$
\t=C\!ab(32,3)T(5,2),\  \r=\C[[t^6,t^{15}+t^{17}]],\
\Gamma=\lan 6,15,32\ran,\ \de=37, 
$$
%and the counterpart of (\ref{bett3}) is
\begin{align} \label{bett3xx}
&\h(q\!=\!1,t,a\!=\!0)=
1+5 t+15 t^2+32 t^3+55 t^4+81 t^5+108 t^6\\
&+134 t^7+157 t^8+175 t^9+186 t^{10}+192 t^{11}+192 t^{12}
+189 t^{13}\notag\\
&+181 t^{14}+172 t^{15}+159 t^{16}+147 t^{17}+132 t^{18}
+120 t^{19}+105 t^{20}\notag\\
&+93 t^{21}+79 t^{22}+68 t^{23}+56 t^{24}+47 t^{25}+37 t^{26}
+30 t^{27}+23 t^{28}\notag\\
&+18 t^{29}+13 t^{30}+10 t^{31}+7 t^{32}
+5 t^{33}+3 t^{34}+2 t^{35}+t^{36}+t^{37}.\notag
\end{align}
The Euler number of $J(\c_{\{5,2\},\{3,2\}})$ 
is $3031$.
\medskip

{\mathversion{bold}$\h_{\{3,4\},\{2,1\}}$.} In this
case, deg$_a=7$ and the Puiseux exponents are $(8,12,13)$.
 One has: 
$$
\t=C\!ab(25,4)T(3,2),\  \r=\C[[t^8,t^{12}+t^{13}]],\
\Gamma=\lan 8,12,25\ran,\ \de=40, 
$$
and the DAHA-Bette polynomial is
\begin{align} \label{bett3xxx}
&\h(q\!=\!1,t,a\!=\!0)=
1+7 t+24 t^2+56 t^3+104 t^4+165 t^5+232 t^6\\
&+297 t^7+355 t^8
+402 t^9+435 t^{10}+454 t^{11}+461 t^{12}+456 t^{13}\notag\\
&+442 t^{14}+420 t^{15}+394 t^{16}
+362 t^{17}+330 t^{18}+295 t^{19}\notag\\
&+262 t^{20}+229 t^{21}+199 t^{22}
+168 t^{23}+143 t^{24}+118 t^{25}\notag\\
&+97 t^{26}+78 t^{27}+63 t^{28}+48 t^{29}
+38 t^{30}+28 t^{31}+21 t^{32}\notag\\
&+15 t^{33}+11 t^{34}
+7 t^{35}+5 t^{36}+3 t^{37}+2 t^{38}+t^{39}+t^{40}.\notag
\end{align}
The Euler number of $J(\c_{\{3,4\},\{2,1\}})$ 
is $7229$. The latter follows from \cite{ObS,Ma};
the Betti numbers are conjectured.

\smallskip
{\mathversion{bold} {\em Double iteration\,}: $\vec\rr=\{3,2,2\},
\vec\ss=\{2,1,1\}$}. The corresponding knot is
$\t=C\!ab(53,2)C\!ab(13,2)T(3,2)$ and deg$_a=7$.
In this case, 
$$
%\t=C\!ab(53,2)C\!ab(13,2)T(3,2),\  
\r=\C[[t^8,t^{12}+t^{14}+t^{15}]],\
\Gamma=\lan 8,12,26,53\ran,\ \de=42, 
$$
and the DAHA-Bette polynomial is
\begin{align} \label{bett3xxxx}
&\h(q\!=\!1,t,a\!=\!0)=
1+7 t+24 t^2+56 t^3+104 t^4+166 t^5\\
&+236 t^6+306 t^7
+370 t^8+424 t^9+465 t^{10}+492 t^{11}+507 t^{12}\notag\\
&+510 t^{13}+504 t^{14}+488 t^{15}+466 t^{16}+437 t^{17}
+406 t^{18}+370 t^{19}\notag\\
&+335 t^{20}+298 t^{21}+264 t^{22}
+230 t^{23}+199 t^{24}+168 t^{25}+143 t^{26}\notag\\
&+118 t^{27}
+97 t^{28}+78 t^{29}+63 t^{30}+48 t^{31}+38 t^{32}+28 t^{33}\notag\\
&+21 t^{34}+15 t^{35}+11 t^{36}
+7 t^{37}+5 t^{38}+3 t^{39}+2 t^{40}+t^{41}+t^{42}.\notag
\end{align}
\medskip
The Euler number
of $J(\c_{\{3,2,2\},\{2,1,1\}})$ is $8512$. This case
is listed in \cite{Pi} (right before Section 4) as an example
beyond the technique there.
\smallskip 

We omit the complete DAHA superpolynomials and the DAHA-Khovanov
polynomials in this and the previous examples (they are all
known).
Quite a few interesting relations for the (conjectural) 
Bette numbers can be seen from these and other examples 
of the DAHA-Bette polynomials. For instance, we expect 
that {\em $b_2$ of the Jacobian factors 
coincide with deg$_a$ of the corresponding DAHA 
superpolynomials\,} (though do not see algebraic reasons 
for this).  

\subsection{\bf The case of two boxes}
We will consider now the iterated knot from the first
example with the weight $\,2\om_1$. We calculated
the one for $\om_2$ too; it is super-dual to the one we 
provide and will be omitted. Also, we will
omit the corresponding DAHA- Khovanov 
polynomial (calculated using the algorithm from
\cite{CJ}), since the Khovanov polynomial
are uncolored.

$$
(\mathbf 5):\ \vec\rr=\{3,2\},\,
\vec\ss\!=\!\{2,1\},\,
\t\!=\! C\!ab(13,2)T(3,2);\, 
\h_{\vec\rr,\vec\ss}(\square\square;q,t,a)=
$$
\renewcommand{\baselinestretch}{0.5} 
{\small
\(
1+t^{16} q^{32}+t^{15} q^{31}+t^{14} q^{31}
+t^{13} q^{31}
+t^{15} q^{30}+2 t^{14} q^{30}+2 t^{13} q^{30}
+2 t^{12} q^{30}
+t^{14} q^{29}+3 t^{13} q^{29}+4 t^{12} q^{29}
+3 t^{11} q^{29}
+t^{14} q^{28}+2 t^{13} q^{28}+5 t^{12} q^{28}
+6 t^{11} q^{28}
+4 t^{10} q^{28}+t^{13} q^{27}+3 t^{12} q^{27}
+7 t^{11} q^{27}
+8 t^{10} q^{27}+4 t^9 q^{27}+t^{13} q^{26}
+2 t^{12} q^{26}
+5 t^{11} q^{26}+10 t^{10} q^{26}+9 t^9 q^{26}
+2 t^8 q^{26}
+t^{12} q^{25}+3 t^{11} q^{25}+7 t^{10} q^{25}
+12 t^9 q^{25}
+7 t^8 q^{25}+t^7 q^{25}+t^{12} q^{24}+2 t^{11} q^{24}
+5 t^{10} q^{24}+10 t^9 q^{24}+14 t^8 q^{24}
+4 t^7 q^{24}
+t^{11} q^{23}+3 t^{10} q^{23}+7 t^9 q^{23}
+12 t^8 q^{23}
+11 t^7 q^{23}+t^6 q^{23}+t^{11} q^{22}+2 t^{10} q^{22}
+5 t^9 q^{22}+10 t^8 q^{22}+13 t^7 q^{22}
+6 t^6 q^{22}
+t^{10} q^{21}+3 t^9 q^{21}+7 t^8 q^{21}
+12 t^7 q^{21}
+10 t^6 q^{21}+t^5 q^{21}+t^{10} q^{20}+2 t^9 q^{20}
+5 t^8 q^{20}+10 t^7 q^{20}+13 t^6 q^{20}+3 t^5 q^{20}
+t^9 q^{19}+3 t^8 q^{19}+7 t^7 q^{19}+11 t^6 q^{19}
+9 t^5 q^{19}+t^9 q^{18}+2 t^8 q^{18}+5 t^7 q^{18}
+10 t^6 q^{18}+11 t^5 q^{18}+2 t^4 q^{18}+t^8 q^{17}
+3 t^7 q^{17}+7 t^6 q^{17}+10 t^5 q^{17}+4 t^4 q^{17}
+t^8 q^{16}+2 t^7 q^{16}+5 t^6 q^{16}+9 t^5 q^{16}
+8 t^4 q^{16}+t^7 q^{15}+3 t^6 q^{15}+7 t^5 q^{15}
+8 t^4 q^{15}+2 t^3 q^{15}+t^7 q^{14}+2 t^6 q^{14}
+5 t^5 q^{14}+8 t^4 q^{14}+3 t^3 q^{14}+t^6 q^{13}
+3 t^5 q^{13}+6 t^4 q^{13}+5 t^3 q^{13}+t^6 q^{12}
+2 t^5 q^{12}+5 t^4 q^{12}+6 t^3 q^{12}+t^2 q^{12}
+t^5 q^{11}+3 t^4 q^{11}+5 t^3 q^{11}+t^2 q^{11}
+t^5 q^{10}+2 t^4 q^{10}+4 t^3 q^{10}+3 t^2 q^{10}
+t^4 q^9+3 t^3 q^9+3 t^2 q^9+t^4 q^8+2 t^3 q^8
+3 t^2 q^8+t^3 q^7+2 t^2 q^7+t q^7+t^3 q^6+2 t^2 q^6
+t q^6+t^2 q^5+t q^5+t^2 q^4+t q^4+t q^3+t q^2
\)

\vfil
\noindent
\(
+a^6 \bigl(t^4 q^{35}+t^3 q^{34}+t^3 q^{33}+t^2 q^{33}
+t^2 q^{32}+t^2 q^{31}+t q^{30}+t q^{29}+q^{27}\bigr)
\)

\vfil
\noindent
\(
+a^5 \bigl(t^7 q^{35}+t^6 q^{35}+t^5 q^{35}+t^7 q^{34}
+2 t^6 q^{34}+2 t^5 q^{34}+t^4 q^{34}+2 t^6 q^{33}
+4 t^5 q^{33}+4 t^4 q^{33}+t^3 q^{33}+t^6 q^{32}
+4 t^5 q^{32}+5 t^4 q^{32}+3 t^3 q^{32}+2 t^5 q^{31}
+5 t^4 q^{31}+5 t^3 q^{31}+t^2 q^{31}+t^5 q^{30}
+4 t^4 q^{30}+6 t^3 q^{30}+3 t^2 q^{30}+2 t^4 q^{29}
+6 t^3 q^{29}+5 t^2 q^{29}+t^4 q^{28}+4 t^3 q^{28}
+5 t^2 q^{28}+t q^{28}+2 t^3 q^{27}+5 t^2 q^{27}
+3 t q^{27}+t^3 q^{26}+4 t^2 q^{26}+3 t q^{26}
+2 t^2 q^{25}+3 t q^{25}+q^{25}+t^2 q^{24}+3 t q^{24}
+q^{24}+2 t q^{23}+q^{23}+t q^{22}+q^{22}+q^{21}
+q^{20}\bigr)
\)

\vfil
\noindent
\(
+a^4 \bigl(t^9 q^{35}+t^8 q^{35}
+t^7 q^{35}+t^{10} q^{34}+2 t^9 q^{34}+4 t^8 q^{34}
+3 t^7 q^{34}+2 t^6 q^{34}+2 t^9 q^{33}+5 t^8 q^{33}
+8 t^7 q^{33}+6 t^6 q^{33}+3 t^5 q^{33}+t^9 q^{32}
+5 t^8 q^{32}+11 t^7 q^{32}+13 t^6 q^{32}+8 t^5 q^{32}
+2 t^4 q^{32}+2 t^8 q^{31}+9 t^7 q^{31}+16 t^6 q^{31}
+15 t^5 q^{31}+6 t^4 q^{31}+t^3 q^{31}+t^8 q^{30}
+5 t^7 q^{30}+16 t^6 q^{30}+21 t^5 q^{30}+14 t^4 q^{30}
+3 t^3 q^{30}+2 t^7 q^{29}+9 t^6 q^{29}+22 t^5 q^{29}
+20 t^4 q^{29}+8 t^3 q^{29}+t^7 q^{28}+5 t^6 q^{28}
+16 t^5 q^{28}+25 t^4 q^{28}+14 t^3 q^{28}+2 t^2 q^{28}
+2 t^6 q^{27}+9 t^5 q^{27}+21 t^4 q^{27}+20 t^3 q^{27}
+5 t^2 q^{27}+t^6 q^{26}+5 t^5 q^{26}+16 t^4 q^{26}
+22 t^3 q^{26}+10 t^2 q^{26}+2 t^5 q^{25}+9 t^4 q^{25}
+20 t^3 q^{25}+13 t^2 q^{25}+2 t q^{25}+t^5 q^{24}
+5 t^4 q^{24}+15 t^3 q^{24}+17 t^2 q^{24}+4 t q^{24}
+2 t^4 q^{23}+9 t^3 q^{23}+15 t^2 q^{23}+6 t q^{23}
+t^4 q^{22}+5 t^3 q^{22}+13 t^2 q^{22}+8 t q^{22}
+q^{22}+2 t^3 q^{21}+8 t^2 q^{21}+9 t q^{21}+q^{21}
+t^3 q^{20}+5 t^2 q^{20}+8 t q^{20}+2 q^{20}+2 t^2 q^{19}
+6 t q^{19}+2 q^{19}+t^2 q^{18}+4 t q^{18}+3 q^{18}
+2 t q^{17}+2 q^{17}+t q^{16}+2 q^{16}+q^{15}+q^{14}\bigr)
\)

\vfil
\noindent
\(
+a^3 \bigl(t^{10} q^{35}+t^{12} q^{34}+2 t^{11} q^{34}
+3 t^{10} q^{34}+3 t^9 q^{34}+t^8 q^{34}+t^{12} q^{33}
+3 t^{11} q^{33}+6 t^{10} q^{33}+8 t^9 q^{33}
+6 t^8 q^{33}+2 t^7 q^{33}+2 t^{11} q^{32}
+8 t^{10} q^{32}+14 t^9 q^{32}+16 t^8 q^{32}
+10 t^7 q^{32}+3 t^6 q^{32}+t^{11} q^{31}
+5 t^{10} q^{31}+15 t^9 q^{31}+25 t^8 q^{31}
+22 t^7 q^{31}+11 t^6 q^{31}+2 t^5 q^{31}
+2 t^{10} q^{30}+10 t^9 q^{30}+26 t^8 q^{30}
+36 t^7 q^{30}+26 t^6 q^{30}+9 t^5 q^{30}
+t^4 q^{30}+t^{10} q^{29}+5 t^9 q^{29}
+18 t^8 q^{29}+39 t^7 q^{29}+42 t^6 q^{29}
+21 t^5 q^{29}+4 t^4 q^{29}+2 t^9 q^{28}
+10 t^8 q^{28}+29 t^7 q^{28}+50 t^6 q^{28}
+39 t^5 q^{28}+12 t^4 q^{28}+t^3 q^{28}
+t^9 q^{27}+5 t^8 q^{27}+18 t^7 q^{27}
+43 t^6 q^{27}+52 t^5 q^{27}+25 t^4 q^{27}
+4 t^3 q^{27}+2 t^8 q^{26}+10 t^7 q^{26}
+29 t^6 q^{26}+53 t^5 q^{26}+41 t^4 q^{26}
+10 t^3 q^{26}+t^8 q^{25}+5 t^7 q^{25}
+18 t^6 q^{25}+43 t^5 q^{25}+50 t^4 q^{25}
+20 t^3 q^{25}+t^2 q^{25}+2 t^7 q^{24}
+10 t^6 q^{24}+29 t^5 q^{24}+50 t^4 q^{24}
+32 t^3 q^{24}+5 t^2 q^{24}+t^7 q^{23}
+5 t^6 q^{23}+18 t^5 q^{23}+41 t^4 q^{23}
+39 t^3 q^{23}+10 t^2 q^{23}+2 t^6 q^{22}
+10 t^5 q^{22}+28 t^4 q^{22}+41 t^3 q^{22}
+18 t^2 q^{22}+t q^{22}+t^6 q^{21}+5 t^5 q^{21}
+18 t^4 q^{21}+36 t^3 q^{21}+24 t^2 q^{21}
+3 t q^{21}+2 t^5 q^{20}+10 t^4 q^{20}+26 t^3 q^{20}
+27 t^2 q^{20}+6 t q^{20}+t^5 q^{19}
+5 t^4 q^{19}+17 t^3 q^{19}+26 t^2 q^{19}
+9 t q^{19}+2 t^4 q^{18}+10 t^3 q^{18}
+21 t^2 q^{18}+12 t q^{18}+q^{18}+t^4 q^{17}
+5 t^3 q^{17}+15 t^2 q^{17}+13 t q^{17}
+q^{17}+2 t^3 q^{16}+9 t^2 q^{16}+12 t q^{16}
+2 q^{16}+t^3 q^{15}+5 t^2 q^{15}+10 t q^{15}
+3 q^{15}+2 t^2 q^{14}+7 t q^{14}+3 q^{14}
+t^2 q^{13}+4 t q^{13}+3 q^{13}+2 t q^{12}
+3 q^{12}+t q^{11}+2 q^{11}+q^{10}+q^9\bigr)
\)

\vfil
\noindent
\(
+a^2 \bigl(t^{13} q^{34}+t^{12} q^{34}
+t^{11} q^{34}+t^{14} q^{33}+2 t^{13} q^{33}
+4 t^{12} q^{33}+4 t^{11} q^{33}+3 t^{10} q^{33}
+2 t^{13} q^{32}+5 t^{12} q^{32}+10 t^{11} q^{32}
+10 t^{10} q^{32}+6 t^9 q^{32}+t^8 q^{32}
+t^{13} q^{31}+5 t^{12} q^{31}+12 t^{11} q^{31}
+20 t^{10} q^{31}+18 t^9 q^{31}+8 t^8 q^{31}
+t^7 q^{31}+2 t^{12} q^{30}+9 t^{11} q^{30}
+21 t^{10} q^{30}+32 t^9 q^{30}+24 t^8 q^{30}
+9 t^7 q^{30}+t^6 q^{30}+t^{12} q^{29}
+5 t^{11} q^{29}+17 t^{10} q^{29}+35 t^9 q^{29}
+45 t^8 q^{29}+25 t^7 q^{29}+7 t^6 q^{29}
+2 t^{11} q^{28}+9 t^{10} q^{28}+27 t^9 q^{28}
+49 t^8 q^{28}+49 t^7 q^{28}+20 t^6 q^{28}
+3 t^5 q^{28}+t^{11} q^{27}+5 t^{10} q^{27}
+17 t^9 q^{27}+42 t^8 q^{27}+62 t^7 q^{27}
+44 t^6 q^{27}+11 t^5 q^{27}+t^4 q^{27}
+2 t^{10} q^{26}+9 t^9 q^{26}+27 t^8 q^{26}
+57 t^7 q^{26}+63 t^6 q^{26}+28 t^5 q^{26}
+3 t^4 q^{26}+t^{10} q^{25}+5 t^9 q^{25}
+17 t^8 q^{25}+42 t^7 q^{25}+70 t^6 q^{25}
+48 t^5 q^{25}+11 t^4 q^{25}+2 t^9 q^{24}
+9 t^8 q^{24}+27 t^7 q^{24}+56 t^6 q^{24}
+66 t^5 q^{24}+24 t^4 q^{24}+2 t^3 q^{24}
+t^9 q^{23}+5 t^8 q^{23}+17 t^7 q^{23}
+42 t^6 q^{23}+67 t^5 q^{23}+42 t^4 q^{23}
+6 t^3 q^{23}+2 t^8 q^{22}+9 t^7 q^{22}
+27 t^6 q^{22}+55 t^5 q^{22}+54 t^4 q^{22}
+15 t^3 q^{22}+t^8 q^{21}+5 t^7 q^{21}
+17 t^6 q^{21}+41 t^5 q^{21}+60 t^4 q^{21}
+26 t^3 q^{21}+2 t^2 q^{21}+2 t^7 q^{20}
+9 t^6 q^{20}+27 t^5 q^{20}+50 t^4 q^{20}
+38 t^3 q^{20}+5 t^2 q^{20}+t^7 q^{19}
+5 t^6 q^{19}+17 t^5 q^{19}+39 t^4 q^{19}
+43 t^3 q^{19}+12 t^2 q^{19}+2 t^6 q^{18}
+9 t^5 q^{18}+26 t^4 q^{18}+41 t^3 q^{18}
+18 t^2 q^{18}+t q^{18}+t^6 q^{17}
+5 t^5 q^{17}+17 t^4 q^{17}+34 t^3 q^{17}
+25 t^2 q^{17}+2 t q^{17}+2 t^5 q^{16}
+9 t^4 q^{16}+24 t^3 q^{16}+25 t^2 q^{16}
+5 t q^{16}+t^5 q^{15}+5 t^4 q^{15}
+16 t^3 q^{15}+25 t^2 q^{15}+8 t q^{15}
+2 t^4 q^{14}+9 t^3 q^{14}+19 t^2 q^{14}
+10 t q^{14}+t^4 q^{13}+5 t^3 q^{13}
+14 t^2 q^{13}+11 t q^{13}+q^{13}
+2 t^3 q^{12}+8 t^2 q^{12}+11 t q^{12}
+q^{12}+t^3 q^{11}+5 t^2 q^{11}+9 t q^{11}
+2 q^{11}+2 t^2 q^{10}+6 t q^{10}+2 q^{10}
+t^2 q^9+4 t q^9+3 q^9+2 t q^8+2 q^8+t q^7
+2 q^7+q^6+q^5\bigr)
\)

\vfil
\noindent
\(
+a \bigl(t^{15} q^{33}
+t^{14} q^{33}+t^{13} q^{33}+t^{15} q^{32}
+2 t^{14} q^{32}+3 t^{13} q^{32}+3 t^{12} q^{32}
+t^{11} q^{32}+2 t^{14} q^{31}+5 t^{13} q^{31}
+8 t^{12} q^{31}+7 t^{11} q^{31}+2 t^{10} q^{31}
+t^{14} q^{30}+4 t^{13} q^{30}+9 t^{12} q^{30}
+14 t^{11} q^{30}+11 t^{10} q^{30}+3 t^9 q^{30}
+2 t^{13} q^{29}+7 t^{12} q^{29}+16 t^{11} q^{29}
+22 t^{10} q^{29}+15 t^9 q^{29}+3 t^8 q^{29}
+t^{13} q^{28}+4 t^{12} q^{28}+12 t^{11} q^{28}
+24 t^{10} q^{28}+29 t^9 q^{28}+14 t^8 q^{28}
+2 t^7 q^{28}+2 t^{12} q^{27}+7 t^{11} q^{27}
+19 t^{10} q^{27}+34 t^9 q^{27}+32 t^8 q^{27}
+11 t^7 q^{27}+t^6 q^{27}+t^{12} q^{26}
+4 t^{11} q^{26}+12 t^{10} q^{26}+28 t^9 q^{26}
+42 t^8 q^{26}+27 t^7 q^{26}+5 t^6 q^{26}
+2 t^{11} q^{25}+7 t^{10} q^{25}+19 t^9 q^{25}
+38 t^8 q^{25}+43 t^7 q^{25}+17 t^6 q^{25}
+t^5 q^{25}+t^{11} q^{24}+4 t^{10} q^{24}
+12 t^9 q^{24}+28 t^8 q^{24}+46 t^7 q^{24}
+34 t^6 q^{24}+6 t^5 q^{24}+2 t^{10} q^{23}
+7 t^9 q^{23}+19 t^8 q^{23}+38 t^7 q^{23}
+46 t^6 q^{23}+17 t^5 q^{23}+t^4 q^{23}
+t^{10} q^{22}+4 t^9 q^{22}+12 t^8 q^{22}
+28 t^7 q^{22}+45 t^6 q^{22}+31 t^5 q^{22}
+4 t^4 q^{22}+2 t^9 q^{21}+7 t^8 q^{21}
+19 t^7 q^{21}+37 t^6 q^{21}+42 t^5 q^{21}
+11 t^4 q^{21}+t^9 q^{20}+4 t^8 q^{20}
+12 t^7 q^{20}+28 t^6 q^{20}+42 t^5 q^{20}
+21 t^4 q^{20}+t^3 q^{20}+2 t^8 q^{19}
+7 t^7 q^{19}+19 t^6 q^{19}+35 t^5 q^{19}
+32 t^4 q^{19}+5 t^3 q^{19}+t^8 q^{18}
+4 t^7 q^{18}+12 t^6 q^{18}+27 t^5 q^{18}
+35 t^4 q^{18}+11 t^3 q^{18}+2 t^7 q^{17}
+7 t^6 q^{17}+19 t^5 q^{17}+31 t^4 q^{17}
+18 t^3 q^{17}+t^2 q^{17}+t^7 q^{16}
+4 t^6 q^{16}+12 t^5 q^{16}+25 t^4 q^{16}
+23 t^3 q^{16}+3 t^2 q^{16}+2 t^6 q^{15}
+7 t^5 q^{15}+18 t^4 q^{15}+24 t^3 q^{15}
+7 t^2 q^{15}+t^6 q^{14}+4 t^5 q^{14}
+12 t^4 q^{14}+21 t^3 q^{14}+10 t^2 q^{14}
+2 t^5 q^{13}+7 t^4 q^{13}+16 t^3 q^{13}
+13 t^2 q^{13}+t q^{13}+t^5 q^{12}
+4 t^4 q^{12}+11 t^3 q^{12}+14 t^2 q^{12}
+2 t q^{12}+2 t^4 q^{11}+7 t^3 q^{11}
+12 t^2 q^{11}+3 t q^{11}+t^4 q^{10}
+4 t^3 q^{10}+9 t^2 q^{10}+5 t q^{10}
+2 t^3 q^9+6 t^2 q^9+6 t q^9+t^3 q^8
+4 t^2 q^8+5 t q^8+2 t^2 q^7+4 t q^7+q^7
+t^2 q^6+3 t q^6+q^6+2 t q^5+q^5+t q^4
+q^4+q^3+q^2\bigr).
\)
}
\renewcommand{\baselinestretch}{1.2} 
\smallskip

{\sf Betti numbers} (presumably):
\begin{align} \label{bett5}
&\h(q\!=\!1,t,a\!=\!0)\!=\!
1+6 t+17 t^2+32 t^3+48 t^4+62 t^5+70 t^6\\
+70& t^7\!+\!64 t^8\!+\!54 t^9\!+\!41 t^{10}\!+\!28 t^{11}\!+
\!18 t^{12}
\!+\!10 t^{13}\!+\!5 t^{14}\!+\!2 t^{15}\!+\!t^{16}\notag\\
&=\bigl(1+3 t+4 t^2+4 t^3+4 t^4+3 t^5+2 t^6+t^7+t^8\bigr)^2.
\notag
\end{align}

We put here (and below) ``presumably", since at the moment
it is not clear which Betti numbers this polynomial calculates
(even conjecturally).
The last relation is a particular case of Part $(iv)$ 
of Theorem \ref{STABILIZ} (the specialization at $q=1$). 
The full specialization formula in this case reads
as follows:

\begin{align*}
&\h_{\{3,2\},\{2,1\}}\,(2\om_1;q\!=\!1,t,a)
\!=\!\bigl(1\!+\!3 a\!+\!3 a^2\!+\!a^3\!+\!3 t\!+\!7 a t\!+\!
5 a^2 t\!+\!a^3 t\!+\!4 t^2\\
&+8 a t^2
+5 a^2 t^2+a^3 t^2+4 t^3+8 a t^3+4 a^2 t^3+4 t^4+6 a t^4+2 a^2 t^4
+3 t^5+\\
&4 a t^5+a^2 t^5+2 t^6+2 a t^6+t^7+a t^7+t^8\bigr)^2
=\h_{\{3,2\},\{2,1\}}\,(\om_1\,;\, q\!=\!1,t,a)^2.
\end{align*}
\smallskip

{\mathversion{bold}$\h_{\{3,2\},\{2,1\}}(\om_2)$.}
The superpolynomial for the same iterated knot and
the weight $\om_2$ is connected with the one for
$2\om_1$ by the Super-Duality, so we will not list it.
However, let us provide its value at $a=0,q=1$ (presumably
a sequence of  Betti numbers):
\begin{align} \label{bett55}
&\h_{\{3,2\},\{2,1\}}\,(\om_2;q\!=\!1,t,a\!=\!0)\!=\!
1\! +\! 3 t\! + 7 t^2\! + 11 t^3\! + 18 t^4\! + 23 t^5\\
+ 29 &t^6+ 31 t^7 + 
36 t^8 + 35 t^9 + 37 t^{10} + 34 t^{11}+ 34 t^{12} 
+31 t^{13} + 31 t^{14}\notag\\
+ 25 &t^{15}
+ 25 t^{16} + 21 t^{17}+ 19 t^{18} + 15 t^{19} + 15 t^{20} + 
10 t^{21} + 10 t^{22}\notag\\ 
+ 7 t&^{23}+ 6 t^{24}+ 4 t^{25}+ 4 t^{26} + 2 t^{27} + 
 2 t^{28} + t^{29} + t^{30} + t^{32}.\notag
\end{align}

A geometric meaning of formulas
(\ref{bett5}) and (\ref{bett55}) 
in terms of  certain Betti numbers
remains to be found. Generally,
we expect that a variant of the construction
from \cite{Ma,DHS} can be used here.
\medskip

\subsection{\bf Pseudo-algebraic knots}
This section is a theoretical challenge.
It indicates that some non-algebraic iterated
knots are very similar to the
algebraic ones. The positivity
of their DAHA-superpolynomials is the key 
``algebraic" feature.
Recall that we call such knots {\em pseudo-algebraic\,},
which conjecturally includes all algebraic knots
and their mirror images for any rectangle Young diagrams. 

The two uncolored knots posses such a positivity and all 
other DAHA properties of algebraic knots, for instance
the existence of the tilde- normalization from
Theorem \ref{JONITER} for any $A_N$ 
and the strict equality for deg$_a$ in (\ref{tilde-sup}).
We give below some other examples of pseudo-algebraic 
knots that do not have these two features; see
(\ref{nac7},\ref{nac-3}).
 
In spite of the presence of negative $\ss_i$ there, we think that
the coefficients of the DAHA-Bette polynomials in 
(\ref{bett6},\ref{bett7}) can be expected to be
the Bette numbers of certain Jacobian factors for the corresponding
curves (\ref{yxcurve}), however this is unknown at the moment.

$$
(\mathbf 6):\ \vec\rr=\{3 ,2\},\,
\vec\ss=\{2,-3\},\,
\t= C\!ab(9,2)T(3,2);\ 
\h_{\,\vec\rr,\,\vec\ss}\,(\square\,;\,\,q,t,a)=
$$
\renewcommand{\baselinestretch}{0.5} 
{\small
\(
1+a^3 q^6+q t+q^2 t+q^3 t+q^2 t^2+q^3 t^2+2 q^4 t^2
+q^3 t^3+q^4 t^3+q^5 t^3+q^4 t^4+q^5 t^4+q^5 t^5
+q^6 t^6 +a^2 \bigl(q^3+q^4+q^5+q^4 t+2 q^5 t
+q^6 t+q^5 t^2+q^6 t^2+q^6 t^3\bigr)
\)

\noindent
\(+a \bigl(q+q^2+q^3+q^2 t+2 q^3 t+3 q^4 t
+q^5 t+q^3 t^2+2 q^4 t^2+3 q^5 t^2+q^4 t^3
+2 q^5 t^3+q^6 t^3+q^5 t^4+q^6 t^4+q^6 t^5\bigr).
\)
}
\renewcommand{\baselinestretch}{1.2} 
\smallskip

{\sf Betti numbers} (presumably):
\begin{align} \label{bett6}
\h(q\!=\!1,t,a\!=\!0)=
1+3 t+4 t^2+3 t^3+2 t^4+t^5+t^6.
\end{align}

This sequence of coefficients and the one in the next example
do not come from algebraic knots. However they fit the table  
of Betti numbers from \cite{Pi} for the Puiseux
exponents $(4,6,v)$. The first two entries there
(our examples $(1),(2)$) are for $(4,6,7),(4,6,9)$;  
the knots from $(6),(7)$  {\em formally\,}
correspond to non-algebraic $(4,6,3),(4,6,5)$. 
\smallskip

{\sf Khovanov polynomial:}
$$(6):\ D\!A\!H\!A'K\!h(C\!ab(9,2)T(3,2))=$$
\renewcommand{\baselinestretch}{0.5} 
{\small
\(
 1 + q^4t^2 + q^6t^3 + q^6t^4 + q^{10}t^5 + q^8t^6 + q^{12}t^7 + 
  2q^{12}t^8 + 2q^{14}t^9 + q^{18}t^{11} + q^{20}t^{12}.
\)
}
\renewcommand{\baselinestretch}{1.2}
 
This polynomial, as well as the next $D\!A\!H\!A'K\!h$\~polynomial,
coincide with the corresponding (reduced) Khovanov polynomials.
This confirms the Connection Conjecture $(ii)$ and also demonstrates
that the differentials necessary to extract the Khovanov polynomials
from the stable ones satisfy the assumptions from Section 3.6 
in \cite{CJ}. 
\medskip

{\mathversion{bold}$\h_{\{3,2\},\{2,-1\}}$.}
This knot is pseudo-algebraic but non-algebraic as well. In a
sense, it is the greatest non-algebraic, since 
$C\!ab(13,2)T(3,2)$ from $(1)$ above is the first algebraic one
in this series.

$$
(\mathbf 7):\ \vec\rr=\{3, 2\},\,
\vec\ss=\{2,-1\},\,
\t= C\!ab(11,2)T(3,2);\ 
\h_{\,\vec\rr,\,\vec\ss}\,(\square\,;\,\,q,t,a)=
$$
\renewcommand{\baselinestretch}{0.5} 
{\small
\(
1+q t+q^2 t+q^3 t+q^2 t^2+q^3 t^2+2 q^4 t^2
+q^3 t^3+q^4 t^3+2 q^5 t^3+q^4 t^4+q^5 t^4
+q^6 t^4+q^5 t^5+q^6 t^5+q^6 t^6+q^7 t^7
+a^3 \bigl(q^6+q^7 t\bigr)
\)

\noindent
\(
+a^2 \bigl(q^3+q^4+q^5+q^4 t+2 q^5 t
+2 q^6 t+q^5 t^2+2 q^6 t^2+q^7 t^2
+q^6 t^3+q^7 t^3+q^7 t^4\bigr)
\)

\noindent
\(
+a \bigl(q+q^2+q^3+q^2 t+2 q^3 t+3 q^4 t
+q^5 t+q^3 t^2+2 q^4 t^2+4 q^5 t^2
+q^6 t^2+q^4 t^3+2 q^5 t^3+3 q^6 t^3
+q^5 t^4+2 q^6 t^4+q^7 t^4+q^6 t^5
+q^7 t^5+q^7 t^6\bigr).
\)
}
\renewcommand{\baselinestretch}{1.2} 
\smallskip

{\sf Betti numbers} (presumably):
\begin{align} \label{bett7}
\h(q\!=\!1,t,a\!=\!0)=
1+3 t+4 t^2+4 t^3+3 t^4+2 t^5+t^6+t^7.
\end{align}

The general formula is as follows:
\begin{align}\label{bet-gen-pos}
&\h_{\{3,2\},\{2,2n-15\}}(q\!=\!1,t,a\!=\!0)=
1+3 t+4\bigl(t^2\!+\!\ldots\!+\!t^{n-4}\bigr)\\
+ 3 &t^{n-3}\!+2 t^{n-2}+t^{n-1}+t^n 
\for C\!ab(2n\!-\!3,2)T(3,2),\ n\!\ge\! 6.\notag
\end{align}

\smallskip
{\sf Khovanov polynomial:}
$$(7):\ D\!A\!H\!A'K\!h(C\!ab(11,2)T(3,2))= 
\hbox{\small $1 + q^4t^2 + q^6t^3$} $$
\renewcommand{\baselinestretch}{0.5} 
{\small
\(+ q^6t^4 + q^{10}t^5 + q^8t^6 + q^{12}t^7 + 
 2q^{12}t^8 + 2q^{14}t^9 + q^{16}t^{10} + 2q^{18}t^{11} + 
q^{20}t^{12}.
\)
}
\renewcommand{\baselinestretch}{1.2} 
\smallskip

This polynomial and the previous one coincide with the
corresponding Khovanov polynomials (under the tilde-normalization);
recall that we use the {\em standard parameters\,} in all formulas
involving Khovanov polynomials. Thus $q,t$ above are actually
$q_{st},t_{st}$.
\medskip

\comment{
%\vskip 1.5in
\begin{figure}[htbp]
\begin{center}
\vskip -0.0in
\hskip -0.5in
\includegraphics[scale=0.8]{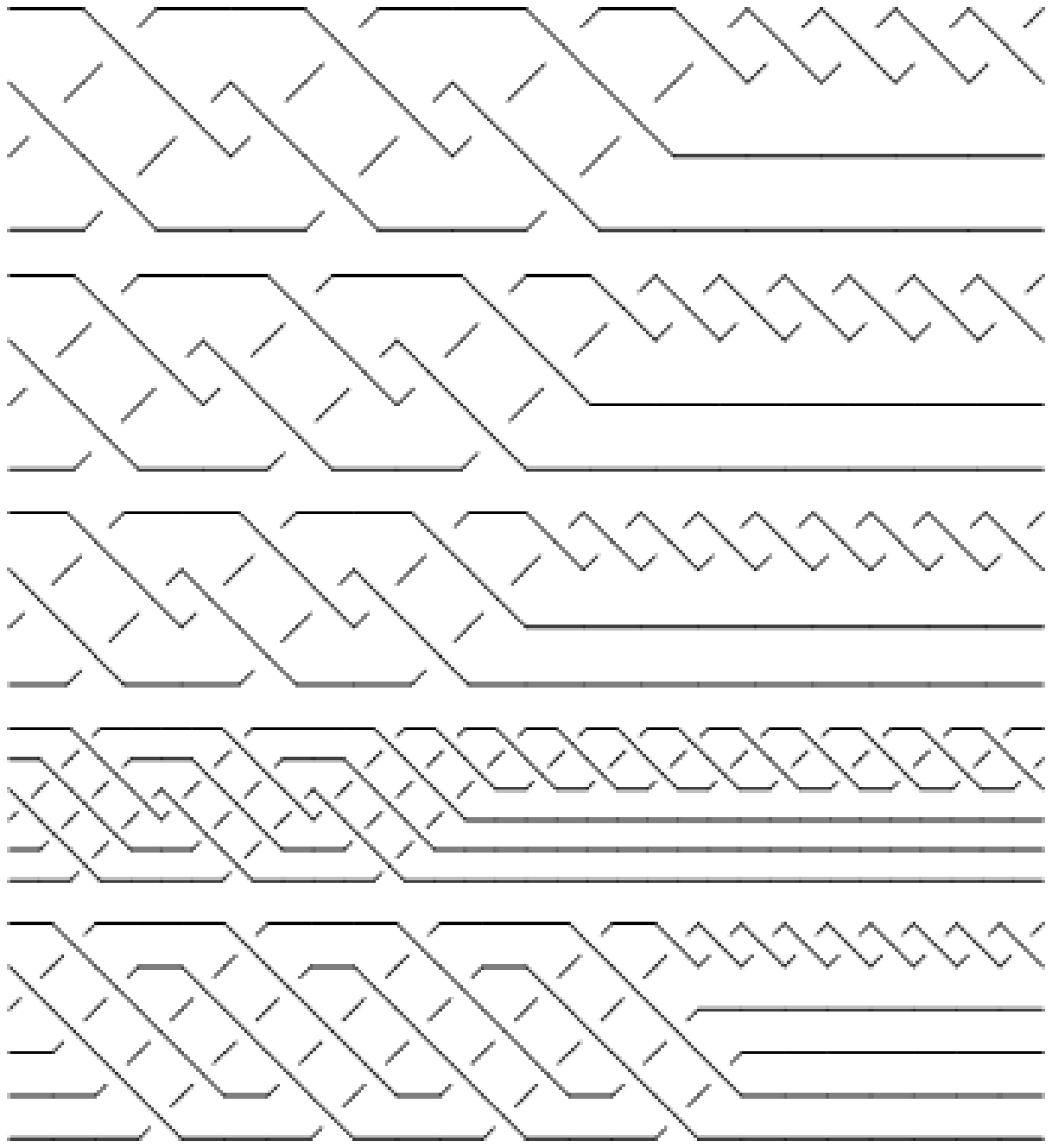}
\vskip +0.in
\caption{
{Cables (11,2),(13,2),(15,2) of T(3,2),}
\newline 
{Cab(19,3)T(3,2) and Cab(25,2)T(4,3)}}
\label{bngt}
\end{center}
\end{figure}
\medskip
}

\section{\sc Further aspects}
\subsection{\bf Non-algebraic knots}
Continuing the previous section, let
us discuss non-algebraic knots that 
are not {\em pseudo-algebraic\,}. 

The first knot we will consider is actually transitional.
It is {\em pseudo-algebraic\,}, but really
different from knots $({\bf 6,7})$
(see (\ref{bett6},\ref{bett7})).
For instance, the corresponding 
\,deg$_a$ is smaller than $(\ss_1\rr_2-1)=3$; 
the coincidence holds for  $({\bf 6,7})$
(and is conjectured for any algebraic ones). 
%The (formal) DAHA-Bette polynomials are
%``chaotic"  for the next knots in this family until 
%$C\!ab(-3,2)T(3,2)$.

$$
({\bf a})\ \vec\rr=\{3,-2\},\,
\vec\ss\!=\!\{2,5\},\,
\t\!=\! C\!ab(7,2)T(3,2);\, 
\h_{\vec\rr,\vec\ss}(\square;q,t,a)=
$$
\renewcommand{\baselinestretch}{0.5} 
{\small
\(
1+q t+q^2 t+q^3 t+q^2 t^2+q^3 t^2+q^4 t^2+q^3 t^3+q^4 t^3+q^4 t^4
+q^5 t^5
+a^2 (q^3+q^4+q^4 t+q^5 t+q^5 t^2)
+a (q+q^2+q^3+q^2 t +2 q^3 t+2 q^4 t
+q^3 t^2+2 q^4 t^2+q^5 t^2 +q^4 t^3+q^5 t^3+q^5 t^4).
\)
}
\renewcommand{\baselinestretch}{1.2} 
\smallskip
\begin{align}\label{nac7}
&\h(q\!=\!1,t,a\!=\!0)\!=\!
1+3 t+3 t^2+2 t^3+ t^4+t^5,\\
&D\!A\!H\!A'K\!h(C\!ab(7,2)T(3,2))\!=\!
1+q^4 t^2+q^6 t^3+q^6 t^4+q^{10} t^5\notag\\
&+q^8 t^6+q^{12} t^7+q^{12} t^8
+q^{14} t^9+q^{18} t^{11}+q^{20} t^{12}.
\notag
\end{align}

The first formula matches (\ref{bet-gen-pos}).
The DAHA-Khovanov polynomial here coincides with
the actual Khovanov one in this case (reduced and
under the tilde normalization); recall
that we use the standard parameters in the
formulas for DAHA-Khovanov polynomials.
\smallskip

$$
({\bf b})\ \vec\rr=\{3,-2\},\,
\vec\ss\!=\!\{2,7\},\,
\t\!=\! C\!ab(5,2)T(3,2);\, 
\h_{\vec\rr,\vec\ss}(\square;q,t,a)=
$$
\renewcommand{\baselinestretch}{0.5} 
{\small
\(
1-\frac{a^3 q^5}{t}+q t+q^2 t+q^2 t^2+q^3 t^2+q^3 t^3
+q^4 t^4+a^2 \left(q^3-q^5-\frac{q^4}{t}+q^4 t\right)\\
+a \left(q+q^2-q^4+q^2 t+2 q^3 t+q^3 t^2+q^4 t^2+q^4 t^3\right).
\)
}
\renewcommand{\baselinestretch}{1.2} 
\smallskip
\begin{align}\label{nac5}
&\h(q\!=\!1,t,a\!=\!0)\!=\!
1+2 t+2 t^2+t^3+t^4,\\
&D\!A\!H\!A'K\!h(C\!ab(5,2)T(3,2))\!=\!
1+q^4 t^2+q^6 t^3+q^6 t^4+q^{10} t^5\notag\\
&+q^{10} t^6+q^{12} t^7
-q^{12} t^9-q^{14} t^{10}-q^{18} t^{12}-q^{20} t^{13}.
\notag
\end{align}

The corresponding DAHA-Jones polynomial for $A_1$ is
\begin{align}\label{jon5-2}
&\tilde{J\!D}_{\{3,-2\},\,\{2,7\}}^{A_1}(\om_1\,;\,q,t)=
\h(q,t,a\!=\!-t^2) =\\
1+q t+q^2 t-&q t^2+q^3 t^2+q^4 t^2
-q^2 t^3-q^3 t^3-q^4 t^3-q^5 t^4+q^5 t^5.\notag
\end{align}

Its value at $t=q$ is the (reduced) {\em Jones 
polynomial\,} $\tilde{J}=1+q^2-q^7-q^9+q^{10}$ 
(for our $q$ and under the tilde normalization); the 
(reduced) one 
produced by the KnotTheory package
is $q^4+q^6-q^{11}-q^{13}+q^{14}$. 

We note that
(\ref{jon5-2}) has little to do with 
the polynomial obtained in \cite{Sam}, Section 5.2
(for the same knot).
We cannot comment on this discrepancy and the general approach 
used in \cite{Sam}. The latter stimulated our one, but the
construction there remains unclear to us. 
\smallskip

The knot $C\!ab(5,2)T(3,2)$ is the first with non-positive 
$\h_{\vec\rr,\vec\ss}(\square;q,t,a)$ in the
considered family. In this and 
further non-positive examples below the actual Khovanov 
polynomials can be obtained from DAHA-Khovanov by the 
substitution $-1\mapsto 1/t$, though the tilde-normalization
creates here some problems and $D\!A\!H\!A'K\!h(C\!ab(1,2)T(3,2))$ 
contains  the term $q^{10}(t^4 - t^6)$, which has no counterpart
in the corresponding $\tilde{K\!h}$, presumably due to the 
(topological) differential. Examples show that such a simple 
substitution is not generally sufficient.
%\smallskip

$$
({\bf c})\ \vec\rr=\{3,-2\},\,
\vec\ss\!=\!\{2,9\},\,
\t\!=\! C\!ab(3,2)T(3,2);\, 
\h_{\vec\rr,\vec\ss}(\square;q,t,a)=
$$
\renewcommand{\baselinestretch}{0.5} 
{\small
\(
1-q^2+q t+q^2 t
-q^3 t+q^2 t^2+q^3 t^3+a^3 \left(-\frac{q^4}{t^2}
-\frac{q^5}{t}\right)+
a^2 \bigl(q^3-q^4-q^5-\frac{q^3}{t^2}-\frac{q^3}{t}
-\frac{2 q^4}{t}\bigr)
+a \left(q+q^2-2 q^3-q^4-\frac{q^2}{t}
-\frac{q^3}{t}+q^2 t+q^3 t-q^4 t+q^3 t^2\right).
\)
}
\renewcommand{\baselinestretch}{1.2} 
\smallskip
\begin{align}\label{nac3}
&\h(q\!=\!1,t,a\!=\!0)\!=\!
t+t^2+t^3,\\
&D\!A\!H\!A'K\!h(C\!ab(3,2)T(3,2))\!=\!
1+q^4 t^2+q^6 t^3-q^4 t^4+q^6 t^4\notag\\
&-q^6 t^5+q^{10} t^5-q^{10} t^7
-q^{12} t^8-q^{12} t^9-q^{14} t^{10}-q^{18} t^{12}-q^{20} t^{13}.
\notag
\end{align}
%\smallskip
$$
({\bf d})\ \vec\rr=\{3,-2\},\,
\vec\ss\!=\!\{2,11\},\,
\t\!=\! C\!ab(1,2)T(3,2);\, 
\h_{\vec\rr,\vec\ss}(\square;q,t,a)=
$$
\renewcommand{\baselinestretch}{0.5} 
{\small
\(
-q+t-2 q^2 t+q t^2-q^3 t^2+q^2 t^3+a^3 \bigl(-q^5-\frac{q^3}{t^2}
-\frac{q^4}{t}\bigr)+a^2 \bigl(-2 q^3-2 q^4
-\frac{q^2}{t^2}-\frac{q^2}{t}-\frac{2 q^3}{t}-q^4 t-q^5 t\bigr)
+a \bigl(-3 q^2-q^3-\frac{q}{t}
-\frac{q^2}{t}+q t-3 q^3 t-q^4 t+q^2 t^2-q^4 t^2\bigr).
\)
}
\renewcommand{\baselinestretch}{1.2} 

%\smallskip
\begin{align}\label{nac1}
&\h(q\!=\!1,t,a\!=\!0)\!=\!
-1-t+t^3,\\
&D\!A\!H\!A'K\!h(C\!ab(1,2)T(3,2))\!=\!
q^2-q^2 t^2+q^6 t^2-q^4 t^3\notag\\
&+q^8 t^3-q^6 t^4
+q^{10} t^4-2 q^8 t^5+q^{12} t^5-q^{10} t^6-q^{12} t^7\notag\\
&-q^{14} t^8
-q^{14} t^9-q^{16} t^{10}-q^{20} t^{12}-q^{22} t^{13}.
\notag
\end{align}
%\smallskip
$$
({\bf e})\ \vec\rr=\{3,-2\},\,
\vec\ss\!=\!\{2,13\},\,
\t\!=\! C\!ab(-1,2)T(3,2);\, 
\h_{\vec\rr,\vec\ss}(\square;q,t,a)=
$$
\renewcommand{\baselinestretch}{0.5} 
{\small
\(
1+2 q t-t^2+q t^2+2 q^2 t^2-q t^3+q^3 t^3
+a^3 \left(q^4+\frac{q^2}{t^2}+\frac{q^3}{t}+q^5 t\right)\\
+a^2 \left(2 q^2+2 q^3+\frac{q}{t^2}+\frac{q}{t}
+\frac{2 q^2}{t}+q^2 t+2 q^3 t+2 q^4 t+q^4 t^2+q^5 t^2\right)\\
+a \left(3 q+q^2+\frac{1}{t}+\frac{q}{t}+q t+4 q^2 t
+q^3 t-q t^2+q^2 t^2+3 q^3 t^2+q^4 t^2+q^4 t^3\right).
\)
}
\renewcommand{\baselinestretch}{1.2} 
\smallskip

Note that this is the first case of nontrivial
$\h(q=0,t,a=0)=1+a/t$; it was $1$ (as for algebraic
knots) for $C\!ab(2\mm+1\ge 3,2)T(3,2)$ and then $t$ before 
this example.
For all further negative $\mm$, 
$\h(q=0,t,a=0)=1+a/t$  in the family $C\!ab(2\mm+1,2)T(3,2)$.

\begin{align}\label{nac-1}
&\h(q\!=\!1,t,a\!=\!0)\!=\!
1+2 t+2 t^2,\\
&D\!A\!H\!A'K\!h(C\!ab(-1,2)T(3,2))\!=\!
1 - q^4 + q^2 t + q^4 t^2 \notag\\
&- q^8 t^2+ q^6 t^3 - q^{10} t^3 + q^8 t^4 + 
2 q^{10} t^5 + q^{14} t^6 + q^{14} t^7 \notag\\
&+ q^{16} t^8 + q^{16} t^9 
+ q^{18} t^{10} + 
q^{22} t^{12} + q^{24} t^{13}.
\notag
\end{align}
%\smallskip
$$
({\bf f})\ \vec\rr=\{3,-2\},\,
\vec\ss\!=\!\{2,15\},\,
\t\!=\! C\!ab(-3,2)T(3,2);\, 
\h_{\vec\rr,\vec\ss}(\square;q,t,a)=
$$
\renewcommand{\baselinestretch}{0.5} 
{\small
\(
1+2 q t+q t^2+2 q^2 t^2-q t^3+q^2 t^3+2 q^3 t^3+q^4 t^4
+a^3 \left(q^4+\frac{q^2}{t^2}+\frac{q^3}{t}+q^5 t+q^6 t^2\right)
\\
+a^2 \bigl(2 q^2+2 q^3+\frac{q}{t^2}+\frac{q}{t}
+\frac{2 q^2}{t}+q^2 t+2 q^3 t+2 q^4 t+q^3 t^2
+2 q^4 t^2+2 q^5 t^2\\
+q^5 t^3+q^6 t^3\bigr)
+a \bigl(3 q+q^2+\frac{1}{t}+\frac{q}{t}+q t
+4 q^2 t
+q^3 t+2 q^2 t^2+4 q^3 t^2
+q^4 t^2\\+q^3 t^3+3 q^4 t^3+q^5 t^3+q^5 t^4\bigr).
\)
}
\renewcommand{\baselinestretch}{1.2} 

%\smallskip
\begin{align}\label{nac-3}
&\h(q\!=\!1,t,a\!=\!0)\!=\!
1+2 t+3 t^2+2 t^3+t^4,\\
&D\!A\!H\!A'K\!h(C\!ab(-3,2)T(3,2))\!=\!
1+q^2 t+q^4 t^2-q^8 t^2+q^6 t^3\notag\\
&+q^8 t^4+q^{10} t^5
+q^{12} t^6+q^{14} t^6+q^{12} t^7+2 q^{14} t^7
+q^{18} t^8\notag\\
&+q^{18} t^9+q^{20} t^{10}
+q^{20} t^{11}+q^{22} t^{12}+q^{26} t^{14}+q^{28} t^{15}.
\notag
\end{align}
%\smallskip
$$
({\bf g})\ \vec\rr=\{3,-2\},\,
\vec\ss\!=\!\{2,17\},\,
\t\!=\! C\!ab(-5,2)T(3,2);\, 
\h_{\vec\rr,\vec\ss}(\square;q,t,a)=
$$
\renewcommand{\baselinestretch}{0.5} 
{\small
\(
1+2 q t+q t^2+2 q^2 t^2+q^2 t^3+2 q^3 t^3+q^3 t^4
+2 q^4 t^4+q^5 t^5+a^3 \bigl(q^4+\frac{q^2}{t^2}
+\frac{q^3}{t}+q^5 t+q^6 t^2+q^7 t^3\bigr)
+a^2 \bigl(2 q^2+2 q^3+\frac{q}{t^2}
+\frac{q}{t}+\frac{2 q^2}{t}+q^2 t
+2 q^3 t+2 q^4 t+q^3 t^2+2 q^4 t^2
+2 q^5 t^2+q^4 t^3+2 q^5 t^3+2 q^6 t^3
+q^6 t^4+q^7 t^4\bigr)+a \bigl(3 q+q^2
+\frac{1}{t}+\frac{q}{t}+q t+4 q^2 t+q^3 t
+2 q^2 t^2+4 q^3 t^2+q^4 t^2+q^2 t^3
+2 q^3 t^3+4 q^4 t^3+q^5 t^3
+q^4 t^4+3 q^5 t^4+q^6 t^4+q^6 t^5\bigr).
\)
}
\renewcommand{\baselinestretch}{1.2} 

%\smallskip
\begin{align}\label{nac-5}
&\h(q\!=\!1,t,a\!=\!0)\!=\!
1 + 2 t + 3 t^2 + 3 t^3 + 3 t^4 + t^5,\\
&D\!A\!H\!A'K\!h(C\!ab(-5,2)T(3,2))\!=\!
1+q^2 t+q^4 t^2+q^6 t^3+q^8 t^4\notag\\
&+q^{10} t^5+q^{12} t^6
+q^{12} t^7+q^{14} t^7+q^{16} t^8
+q^{18} t^8+q^{16} t^9+2 q^{18} t^9
\notag\\
&+q^{22} t^{10}+q^{22} t^{11}+q^{24} t^{12}
+q^{24} t^{13}+q^{26} t^{14}+q^{30} t^{16}+q^{32} t^{17}.
\notag
\end{align}

One has that $\tilde{K\!h}-D\!A\!H\!A'K\!h=
q^{12}(t^5 - t^7) + q^{16}(t^7 - t^9),$ as can be
expected, though the first difference is 
obviously not of ``top" $q$\~degree.

This is a pseudo-algebraic knot, as well as all
further knots in the considered family, however 
numerically quite different from those in 
(\ref{bett6},\ref{bett7}). The polynomial 
$\h(q\!=\!1,t,a\!=\!0)$  becomes
$1+2 t+3 t^2+3 t^3+4 t^4+3 t^5+t^6$ for $C\!ab(-7,2)T(3,2)$
and 
\begin{align}\label{bet-gen-neg}
&\h_{\{3,-2\},\{2,2n+7\}}(q\!=\!1,t,a\!=\!0)=
1+2 t+3 t^2+3 t^3\\
+4\bigl(t^4\!+&\!\ldots\!+\!t^{n-2}\bigr)+
3 t^{n-1}\!+\!t^n \for C\!ab(-2n\!+\!5,2)T(3,2),\, n\!\ge\! 5.\notag
\end{align}
\medskip

Let us consider a more involved
{\em pseudo-algebraic\,} knot:

$$
({\bf h})\ \vec\rr=\{4,-2\},\,
\vec\ss\!=\!\{3,7\},\,
\t\!=\! C\!ab(17,2)T(4,3);\, 
\h_{\vec\rr,\vec\ss}(\square;q,t,a)=
$$
\renewcommand{\baselinestretch}{0.5} 
{\small
\(
1+q t+q^2 t+q^3 t+q^4 t+q^5 t+q^2 t^2+q^3 t^2+2 q^4 t^2
+2 q^5 t^2+3 q^6 t^2+q^7 t^2+q^8 t^2+q^3 t^3+q^4 t^3+2 q^5 t^3
+3 q^6 t^3+4 q^7 t^3+3 q^8 t^3+2 q^9 t^3+q^4 t^4+q^5 t^4+2 q^6 t^4
+3 q^7 t^4+5 q^8 t^4+4 q^9 t^4+3 q^{10} t^4+q^5 t^5+q^6 t^5
+2 q^7 t^5+3 q^8 t^5+5 q^9 t^5+4 q^{10} t^5+2 q^{11} t^5+q^6 t^6
+q^7 t^6+2 q^8 t^6+3 q^9 t^6+5 q^{10} t^6+3 q^{11} t^6+q^{12} t^6
+q^7 t^7+q^8 t^7+2 q^9 t^7+3 q^{10} t^7+4 q^{11} t^7+q^{12} t^7
+q^8 t^8+q^9 t^8+2 q^{10} t^8+3 q^{11} t^8+3 q^{12} t^8+q^9 t^9
+q^{10} t^9+2 q^{11} t^9+2 q^{12} t^9+q^{13} t^9+q^{10} t^{10}
+q^{11} t^{10}+2 q^{12} t^{10}+q^{13} t^{10}+q^{11} t^{11}
+q^{12} t^{11}+q^{13} t^{11}+q^{12} t^{12}+q^{13} t^{12}
+q^{13} t^{13}+q^{14} t^{14}+a^4 \bigl(q^{10}+q^{11}+q^{12}
+q^{11} t+2 q^{12} t+q^{13} t+q^{12} t^2+2 q^{13} t^2+q^{14} t^2
+q^{13} t^3+q^{14} t^3+q^{14} t^4\bigr)+a^3 \bigl(q^6+q^7+2 q^8
+2 q^9+2 q^{10}+q^{11}+q^7 t+2 q^8 t+4 q^9 t+5 q^{10} t+5 q^{11} t
+2 q^{12} t+q^8 t^2+2 q^9 t^2+5 q^{10} t^2+7 q^{11} t^2
+7 q^{12} t^2+2 q^{13} t^2+q^9 t^3+2 q^{10} t^3+5 q^{11} t^3
+7 q^{12} t^3+5 q^{13} t^3+q^{14} t^3+q^{10} t^4+2 q^{11} t^4
+5 q^{12} t^4+5 q^{13} t^4+2 q^{14} t^4+q^{11} t^5+2 q^{12} t^5
+4 q^{13} t^5+2 q^{14} t^5+q^{12} t^6+2 q^{13} t^6+2 q^{14} t^6
+q^{13} t^7+q^{14} t^7+q^{14} t^8\bigr)+a^2 \bigl(q^3+q^4+2 q^5
+2 q^6+2 q^7+q^8+q^9+q^4 t+2 q^5 t+4 q^6 t+6 q^7 t+7 q^8 t+6 q^9 t
+4 q^{10} t+q^{11} t+q^5 t^2+2 q^6 t^2+5 q^7 t^2+8 q^8 t^2
+12 q^9 t^2+11 q^{10} t^2+7 q^{11} t^2+q^{12} t^2+q^6 t^3+2 q^7 t^3
+5 q^8 t^3+9 q^9 t^3+14 q^{10} t^3+14 q^{11} t^3+7 q^{12} t^3
+q^{13} t^3+q^7 t^4+2 q^8 t^4+5 q^9 t^4+9 q^{10} t^4+14 q^{11} t^4
+11 q^{12} t^4+4 q^{13} t^4+q^8 t^5+2 q^9 t^5+5 q^{10} t^5
+9 q^{11} t^5+12 q^{12} t^5+6 q^{13} t^5+q^{14} t^5+q^9 t^6
+2 q^{10} t^6+5 q^{11} t^6+8 q^{12} t^6+7 q^{13} t^6+q^{14} t^6
+q^{10} t^7+2 q^{11} t^7+5 q^{12} t^7+6 q^{13} t^7+2 q^{14} t^7
+q^{11} t^8+2 q^{12} t^8+4 q^{13} t^8+2 q^{14} t^8+q^{12} t^9
+2 q^{13} t^9+2 q^{14} t^9+q^{13} t^{10}+q^{14} t^{10}
+q^{14} t^{11}\bigr)+a \bigl(q+q^2+q^3+q^4+q^5+q^2 t+2 q^3 t
+3 q^4 t+4 q^5 t+5 q^6 t+3 q^7 t+2 q^8 t+q^9 t+q^3 t^2+2 q^4 t^2
+4 q^5 t^2+6 q^6 t^2+9 q^7 t^2+8 q^8 t^2+6 q^9 t^2+2 q^{10} t^2
+q^4 t^3+2 q^5 t^3+4 q^6 t^3+7 q^7 t^3+11 q^8 t^3+12 q^9 t^3
+10 q^{10} t^3+3 q^{11} t^3+q^5 t^4+2 q^6 t^4+4 q^7 t^4+7 q^8 t^4
+12 q^9 t^4+13 q^{10} t^4+10 q^{11} t^4+2 q^{12} t^4+q^6 t^5
+2 q^7 t^5+4 q^8 t^5+7 q^9 t^5+12 q^{10} t^5+12 q^{11} t^5
+6 q^{12} t^5+q^{13} t^5+q^7 t^6+2 q^8 t^6+4 q^9 t^6+7 q^{10} t^6
+11 q^{11} t^6+8 q^{12} t^6+2 q^{13} t^6+q^8 t^7+2 q^9 t^7
+4 q^{10} t^7+7 q^{11} t^7+9 q^{12} t^7+3 q^{13} t^7+q^9 t^8
+2 q^{10} t^8+4 q^{11} t^8+6 q^{12} t^8+5 q^{13} t^8+q^{10} t^9
+2 q^{11} t^9+4 q^{12} t^9+4 q^{13} t^9+q^{14} t^9+q^{11} t^{10}
+2 q^{12} t^{10}+3 q^{13} t^{10}+q^{14} t^{10}+q^{12} t^{11}
+2 q^{13} t^{11}
+q^{14} t^{11}+q^{13} t^{12}+q^{14} t^{12}+q^{14} t^{13}\bigr).
\)
}
\renewcommand{\baselinestretch}{1.2} 
\smallskip

The corresponding $D\!A\!H\!A'K\!h(C\!ab(17,2)T(4,3))$
equals
\renewcommand{\baselinestretch}{0.5} 
{\small
\(
1+q^4 t^2+q^6 t^3+q^6 t^4+q^{10} t^5+q^8 t^6
+q^{12} t^7+q^{10} t^8+q^{12} t^8+2 q^{14} t^9+q^{12} t^{10}
+q^{14} t^{10}+2 q^{16} t^{11}+q^{18} t^{11}+3 q^{16} t^{12}
+q^{20} t^{12}+2 q^{18} t^{13}+2 q^{20} t^{13}+q^{18} t^{14}
+q^{20} t^{14}+q^{22} t^{14}+4 q^{22} t^{15}+q^{20} t^{16}
+q^{22} t^{16}+q^{24} t^{16}+q^{26} t^{16}+2 q^{24} t^{17}
+q^{26} t^{17}+q^{24} t^{18}+q^{28} t^{18}+q^{26} t^{19}
+q^{28} t^{19}+q^{28} t^{20}+q^{30} t^{20}
+2 q^{30} t^{21}+q^{32} t^{22}+q^{36} t^{24}+q^{34} t^{25}.
\)
}
\renewcommand{\baselinestretch}{1.2} 
%\smallskip

One has:
$\tilde{K\!h}-D\!A\!H\!A'K\!h=q^{34}(t^{23} - t^{25})$.
The corresponding (reduced and tilde-normalized) 
HOMFLYPT polynomial coincides with   
$\h(\square;q,q,-a)$ as well as in all examples
we considered.

\comment{
$$
({\bf i})\ \vec\rr=\{4,-2\},\,
\vec\ss\!=\!\{3,7\},\,
\t\!=\! C\!ab(7,2)T(4,3);\, 
\h_{\vec\rr,\vec\ss}(\square;q,t,a)=
$$
\renewcommand{\baselinestretch}{0.5} 
{\small
\(
-q^4-2 q^5 t-q^6 t+t^2-q^5 t^2-4 q^6 t^2-q^7 t^2+q t^3+q^2 t^3
+q^3 t^3+q^4 t^3+q^5 t^3-2 q^6 t^3-6 q^7 t^3-q^8 t^3+q^2 t^4
+q^3 t^4+2 q^4 t^4+2 q^5 t^4+2 q^6 t^4-3 q^7 t^4-6 q^8 t^4
-q^9 t^4+q^3 t^5+q^4 t^5+2 q^5 t^5+3 q^6 t^5+2 q^7 t^5-2 q^8 t^5
-4 q^9 t^5-q^{10} t^5+q^4 t^6+q^5 t^6+2 q^6 t^6+2 q^7 t^6+q^8 t^6
-q^9 t^6-2 q^{10} t^6+q^5 t^7+q^6 t^7+2 q^7 t^7+q^8 t^7-q^{11} t^7
+q^6 t^8+q^7 t^8+q^8 t^8+q^7 t^9+q^8 t^9+q^8 t^{10}+q^9 t^{11}
+a^5 \bigl(-q^{13}-\frac{q^{10}}{t^3}-\frac{q^{11}}{t^2}
-\frac{q^{12}}{t}-q^{14} t\bigr)+a^4 \bigl(-q^9-3 q^{10}-3 q^{11}
-2 q^{12}-\frac{q^8}{t^3}-\frac{q^9}{t^3}-\frac{q^8}{t^2}
-\frac{2 q^9}{t^2}-\frac{2 q^{10}}{t^2}-\frac{2 q^9}{t}
-\frac{3 q^{10}}{t}-\frac{2 q^{11}}{t}-2 q^{10} t-3 q^{11} t
-3 q^{12} t-2 q^{13} t-q^{11} t^2-2 q^{12} t^2-2 q^{13} t^2
-q^{14} t^2-q^{13} t^3-q^{14} t^3\bigr)+a^3 \bigl(-q^7-7 q^8
-8 q^9-6 q^{10}-q^{11}-\frac{q^7}{t^3}-\frac{q^6}{t^2}
-\frac{2 q^7}{t^2}-\frac{3 q^8}{t^2}-\frac{q^9}{t^2}
-\frac{3 q^7}{t}-\frac{5 q^8}{t}-\frac{5 q^9}{t}-\frac{q^{10}}{t}
-3 q^8 t-11 q^9 t-9 q^{10} t-6 q^{11} t-q^{12} t+q^6 t^2+q^7 t^2
+q^8 t^2-4 q^9 t^2-11 q^{10} t^2-8 q^{11} t^2-5 q^{12} t^2
-q^{13} t^2+q^7 t^3+2 q^8 t^3+q^9 t^3-3 q^{10} t^3-7 q^{11} t^3
-5 q^{12} t^3-3 q^{13} t^3+q^8 t^4+q^9 t^4-q^{11} t^4-3 q^{12} t^4
-2 q^{13} t^4-q^{14} t^4+q^9 t^5-q^{13} t^5\bigr)
+a^2 \bigl(-3 q^6-10 q^7-7 q^8-3 q^9-\frac{q^5}{t^2}
-\frac{q^6}{t^2}-\frac{q^7}{t^2}-\frac{q^5}{t}-\frac{4 q^6}{t}
-\frac{4 q^7}{t}-\frac{2 q^8}{t}-q^6 t-8 q^7 t-16 q^8 t-9 q^9 t
-3 q^{10} t+q^3 t^2+q^4 t^2+2 q^5 t^2+2 q^6 t^2-q^7 t^2-13 q^8 t^2
-19 q^9 t^2-9 q^{10} t^2-3 q^{11} t^2+q^4 t^3+2 q^5 t^3+4 q^6 t^3
+5 q^7 t^3-q^8 t^3-13 q^9 t^3-16 q^{10} t^3-7 q^{11} t^3
-2 q^{12} t^3+q^5 t^4+2 q^6 t^4+5 q^7 t^4+5 q^8 t^4-q^9 t^4
-8 q^{10} t^4-10 q^{11} t^4-4 q^{12} t^4-q^{13} t^4+q^6 t^5
+2 q^7 t^5+4 q^8 t^5+2 q^9 t^5-q^{10} t^5-3 q^{11} t^5-4 q^{12} t^5
-q^{13} t^5+q^7 t^6+2 q^8 t^6+2 q^9 t^6-q^{12} t^6-q^{13} t^6
+q^8 t^7+q^9 t^7+q^9 t^8\bigr)+a \bigl(-3 q^5-4 q^6-2 q^7
-\frac{q^4}{t}-\frac{q^5}{t}-\frac{q^6}{t}-q^5 t-7 q^6 t-8 q^7 t
-3 q^8 t+q t^2+q^2 t^2+q^3 t^2+q^4 t^2+q^5 t^2-3 q^6 t^2-13 q^7 t^2
-10 q^8 t^2-3 q^9 t^2+q^2 t^3+2 q^3 t^3+3 q^4 t^3+4 q^5 t^3
+4 q^6 t^3-4 q^7 t^3-16 q^8 t^3-10 q^9 t^3-3 q^{10} t^3+q^3 t^4
+2 q^4 t^4+4 q^5 t^4+6 q^6 t^4+6 q^7 t^4-4 q^8 t^4-13 q^9 t^4
-8 q^{10} t^4-2 q^{11} t^4+q^4 t^5+2 q^5 t^5+4 q^6 t^5+6 q^7 t^5
+4 q^8 t^5-3 q^9 t^5-7 q^{10} t^5-4 q^{11} t^5-q^{12} t^5+q^5 t^6
+2 q^6 t^6+4 q^7 t^6+4 q^8 t^6+q^9 t^6-q^{10} t^6-3 q^{11} t^6
-q^{12} t^6+q^6 t^7+2 q^7 t^7+3 q^8 t^7+q^9 t^7-q^{12} t^7
+q^7 t^8+2 q^8 t^8+q^9 t^8+q^8 t^9+q^9 t^9+q^9 t^{10}\bigr)
\)
}
\renewcommand{\baselinestretch}{1.2} 
\smallskip
}

%%%%%%%%%%%%%

\subsection{\bf The case of 
\texorpdfstring{{\mathversion{bold}$A_1$}}{A1}}
We are going to prove Part $(i)$ of the Connection Conjecture
in the $A_1$\~case, i.e. the equality
\begin{align}\label{homjon2x}
\h(b\om_1;q,q,-q^2)\!=\!
\tilde{\hbox{\small\em H\!O\!M}}(b\om_1;q,q^2)
\!=\!\tilde{J\!D}^{A_1}(b\om_1;q,q), \ b\in \Z_+.
\end{align}

We need the simplest case of the theory of the
difference {\em shift operators\,}; see  
Theorem 2.4 of \cite{C4} and \cite{CJ}. 
Let $k=\{k_\nu=0,1\}$. Generally, we define
\begin{align}\label{shiftx}
\x_k  = \prod_{\al\in R_+}\prod_{j=0}^{k_\al-1}
\bigl((q_\al^{j} X_{\al})^{1/2}-
(q_\al^{j} X_{\al})^{-1/2})\bigr).
\end{align}
Recall that $q_\al=q^{\nu_{\al}}$,
where $\nu_\al\equal (\al,\al)/2$.

We put $H=H^{(k)}$ for elements $H$ from the double
affine Hecke algebra
$\HH^{(k)}$ with the structural parameters
$q$ and $t_\nu=q_\nu^{k_\nu}$ and their images in the
corresponding polynomial representation $\v^{(k)}$
Accordingly, $\tau_{\pm}^{(k)}$ will be automorphisms of
$\HH^{(k)}$.

We set $\circ=(\{k_\nu=0\})$. Then the operator 
$H^{\circ}$ is obtained from $H$ by
replacing every $Y_b (b\in P)$ by the difference
operators $b^{-1}$ and $T_w (w\in W)$ by $w$. Accordingly,
\begin{align}\label{tauzero}
&\tau_{+}^\circ (Y_b)=q^{-(b,b)/2}X_b Y_b,\ \,
\tau_{-}^\circ (X_b)=q^{(b,b)/2}Y_b X_b,\ b\in P.
%&\hbox{where\ \ }\tau_{\pm}^\circ\, \equal\, \tau_{\pm}^{(0)} \for
%b\in P,\, 0=\{k_\nu=0\}.\notag
\end{align}

\begin{lemma}\label{SHIFTAU}
Let $H$ be an algebraic expression with complex coefficients
in terms of the standard $W$\~symmetrizations
of monomials with respect to
$\{X_b\}$ and those with respect to $\{Y_b\}$. The restriction
of the corresponding operator to 
the subsubspace $(\v^{(k)})^W$ of $W$\~invariant elements 
in $\v^{(k)}$ will be denoted by $H^{(k)}_{sym}$. Then
\begin{align}%\label{xh0x}
&H^{(k)}_{sym}\!=\! \x_k^{-1}\, H_{sym}^{\circ}\, \x_k,\ 
\bigr(\tau_{\pm}^{(k)}(H^{(k)})\bigl)_{sym} =\! 
\x_k^{-1}\, \bigl( \tau_{\pm}^\circ(
H^{\circ})\bigr)_{sym}\, \x_k.\notag%\label{xhtaux}
\end{align}
\vskip -1.4cm 
\sq 
\end{lemma}
%\medskip

We will also treat $\tau_{\pm}^\circ$ and
any elements $\ga\in SL(2,\Z)$ as automorphisms of
the linear span $\tilde{\v}$ of $X_\la q^{M x^2/2}$
for $\la, M$ from sufficiently general subsets of $\C$.
Let us formally set:
$$
X_\la=q^{(\la,x)},\ 
x^2=(x,x)=\sum_{i=1}^n (x,\om_i)(x,\al_i^\vee).
$$
The action of $\hat{W}$ in $\v$ is naturally extended
to $\tilde{\v}$ through its action on $\{x_\la\}$; note
that $q^{M x^2/2}$ are symmetric ($W$\~invariant).
\smallskip

The main formulas we need are as follows:
\begin{align}\label{sl2Mq}
&\si^\circ(X_\la\,q^{-M x^2/2})\,=\,
\frac{q^{\la^2/(2M)}}{M^{1/2}}\,
X_{\la/M}\,q^{+x^2/(2M)} \for M\neq 0,\\
&(\tau_+^\circ)^N\,(X_\la\,q^{-M x^2/2})\ =
\ X_\la\,q^{(N-M)x^2/2}\for N\in \Z
\and \notag\\ 
&(\tau_-^\circ)^N\,(X_\la\,q^{-M x^2/2})\ =\ 
((\si^\circ)^{-1}\,(\tau_+^\circ)^{-N}\,
\si^\circ)(X_\la\,e^{-M x^2/2})
\notag\\
&=\ \frac{1}{(1-MN)^{1/2}}\,
q^{\frac{\la^2\,N}{2(1-MN)}}\
X_{\la/(1-MN)}\, q^{-x^2\frac{M}{2(1-MN)}}.\notag
\end{align}

Here $N\in \Z$ and the parameters $M$ are sufficiently
general complex numbers; for instance $MN\neq 1$ 
in the last formula. We also assume that $(MM')^{1/2}$
$= M^{1/2}(M')^{1/2}$ for $M,M'\in \C$. 
\smallskip

These formulas can be readily extended to 
the following
{\em free action\,} of the whole $SL(2,\Z)$ in $\tilde{\v}$:
\begin{align}\label{sl2zq}
&\ga^\circ\,(X_\la\,q^{z x^2/2})\,=\,
\frac{1}{(cz+d)^{1/2}}\ q^{-\frac{\la^2\, c}{2(cz+d)}}\
X_{\frac{\la}{cz+d}}\ q^{\frac{az+b}{cz+d}\,x^2/2} \for\\
&\ga= \left(
  \begin{array}{cc}
    a & b \\
    c & d \\
  \end{array}
\right)\in SL(2,\Z),\ \la \in \C 
\hbox{\ and for generic\, } z\in \C\,. \notag
\end{align}

\comment{
The mere fact of its existence makes it possible
to deduce  relation (\ref{xhtaux}) from
(\ref{xh0x}). Indeed, 
$\tau_+$ and $\tau_-$ become conjugations 
by $W$\~invariant series in terms of $\{X_b\}$ 
and $\{Y_b\}$ upon the
proper completion of $\HH$ . For instance,
we can set 
$$
Y_b^{(0)}=q^{y_b}=e^{-\partial_b} \for
y_b=-\frac{\partial_b}{\log q},\ 
\partial_b(x_a)=(b,a),
$$
where the notation $ x_a=(x,a)$ is used
($a\in \C^n$).
Then $\tau_-^{(0)}$ acts in $\tilde{\v}$ as
$\tau_-^\circ=q^{-y^2/2}=e^{-\frac{\partial^2}{2\log q}}$;
here $y^2$ and $\partial^2$
are defined for the standard inner product
$(\,,\,)$ as for $x$. 
}
%\smallskip

{\sf DAHA-Jones polynomials.} 
We will now switch to the $A_1$\~case. 
Let $\al=\al_1$, $s=s_1$ and $\om=\om_1$;
then $\alpha=\al_1=2\omega$ and
$\rho=\om$. The extended affine Weyl group is
$\widehat{W}=\lan s,\om\ran$.
The weights $b\om$ ($b\in \Z$) will be denoted 
simply by $b$.
\smallskip

The {\em double affine Hecke algebra
$\HH$\,} is generated by invertible elements
$Y=Y_{\om}, T=T_1, X=X_{\om}$
subject to the quadratic relation $(T-t^{1/2})(T+t^{-1/2})=0$
and the cross-relations:
\begin{align}\label{dahaone}
&TXT=X^{-1},\ T^{-1}YT^{-1}=Y^{-1},\ Y^{-1}X^{-1}YXT^2q^{1/2}=1.
\end{align}
%Setting $\pi\equal YT^{-1}$, the second relation 
%becomes $\pi^2=1$.
The field of definition will be $\Q(q^{1/4},t^{1/2})$.
% though the ring $\Z[q^{1/4},t^{1/2}]$ is sufficient 
%for DAHA-Jones polynomials.
Here $q^{\pm 1/4}$ is needed for the automorphisms $\tau_{\pm}$:
\begin{align}\label{taudef}
&\tau_+(X)=X,\ \ \tau_+(T)=T,\ \ \tau_+(Y)=q^{-1/4}XY,\ \\
%\tau_+(\pi)=q^{-1/4}X\pi,\\
%\label{tau-def}
&\tau_-(Y)=Y,\ \ \,\tau_-(T)=T,\ \ \,\tau_-(X)=q^{1/4}YX.\ 
%\tau_-(\pi)=\pi.
\notag
\end{align}
\smallskip

We will prove the Connection Conjecture 
for the {\em non-reduced\,} version of the
DAHA-Jones polynomials. Namely, we modify
(\ref{jones-dit}) as follows:

\begin{align}\label{nonrjone}
&J\!D^\#_{\,\vec\rr,\,\vec\ss\,}(b\,;\,q,t) \equal 
\Bigl\{\hat{\ga_{1}}\Bigr(
\cdots\Bigl(\hat{\ga}_{\ell-1}
\Bigl(\bigl(\hat{\ga}_\ell(P_b^{(k)})\bigr)\!\Downarrow
\Bigr)\!\Downarrow\Bigr) \cdots\Bigr)\Bigr\}^{(k)}, \\
&\for t \equal q^k\hbox{\,\, and\,\, } 
\ga=\ga_{r,s} \hbox{\,\, lifted to\,\, }
\hat{\ga}\in PSL_{\,2}^{\wedge}(\Z),
\notag
\end{align}
where we show explicitly the dependence of the $P$\~polynomial
on the parameter $k$ and also use ${}^{(k)}$ for the
corresponding coinvariant. 
Generally, the non-reduced DAHA-Jones polynomials are 
inconvenient because of their non-trivial $q,t$\~denominators. 
However they become $q$\~polynomials as $t=q$ and simplify 
the considerations. 

\smallskip

\subsection{\bf Jones polynomials}
We are going to use the
relation for the colored Jones polynomials of the cable 
$C\!ab(\aa,\rr)K$ with that 
of the initial knot $K$. See \cite{Mo,Tr} and references therein.
Note that $t$ in \cite{Tr} is our $q^{1/4}$,  
%we use $(\aa,\rr)$ parameters for the cabling to distinguish them 
%from the $(\rr,\ss)$ parameters used in the construction of the 
%of the DAHA-Jones polynomials.  
and $C\!ab(r,s)$ there is our $C\!ab(\aa,\rr)$  (in this order!). 

Assuming here and below that the $j$-summation has step $1$, one has:
\begin{align}\label{jones-formula}
&\j^{\#}_{C\!ab(\aa,\rr)K}(b\,;\,q) = q^{-\frac{1}{4}\aa\rr\, b(b+2)} 
\sum_{j=-b/2}^{b/2} q^{\aa\, j(j\,\rr +1)} \j^{\#}_{K}(2\rr j\,;\,q).
\end{align}

This formula gives a recursive definition of the Jones 
polynomials of iterated knots; the unknot is normalized as
follows:

\begin{align}
\j^{\#}_{\unknot}(b\,;\,q) = (q^{(b+1)/2}-q^{-(b+1)/2})/
(q^{1/2}-q^{-1/2}).
\end{align}
\smallskip

\begin{proposition}\label{RELATIONTOJONES}
For any (possibly negative)
$\{\rr,\ss\}$ from Theorem \ref{JONITER}, 
\begin{align}
&J\!D^\#_{\{\rr_1,\ldots \rr_\ell\},\,\{\ss_1,\ldots \ss_\ell\}\,}
(b\,;\,q,t\mapsto q)\\
= q^{-\frac{1}{4}\aa_\ell \rr_\ell\, b(b+2)}
\j&^{\#}_{T(\{\rr_1,\ldots 
\rr_\ell\},\,\{\ss_1,\ldots \ss_\ell\})\,}(b\,;\, q\mapsto q^{-1}),
\where 
\notag
%&= q^{-\frac{1}{4}\aa_\ell \rr_\ell b(b+2)} 
%\j^{\#}_{\bigl(C\!ab(\aa_\ell,\rr_\ell)
%\cdots C\!ab(\aa_1,\rr_1)\bigr)
%(\unknot )}(b\,;\, q\mapsto q^{-1}), \notag
\end{align}
$T(\{\rr_1,\ldots \rr_\ell\},\,\{\ss_1,\ldots \ss_\ell\}) = 
C\!ab(\vec\aa,\vec\rr)=\bigl(C\!ab(\aa_\ell,\rr_\ell)\cdots 
C\!ab(\aa_1,\rr_1)\bigr)(\unknot )$,
and $\aa_i$ are recursively defined by relations $\aa_{i+1} = 
\rr_i \rr_{i+1} \aa_i + \ss_{i+1}$ for $1 \leq i \leq \ell - 1$ 
and the initial condition $\aa_1 = \ss_1$. 
\end{proposition}
\smallskip

{\it Proof.}
As above, we formally set $X_\la=q^{\la x}$ for $\la\in \C$,
which extends $X_b=X_{b\om_1}=X^b$ for $b$ in $\Z$
(identified with the lattice $P=\Z\om_1$). The 
{\em free action\,} of $SL(2,\Z)$ 
on $X_\la\,q^{M x^2}$ from (\ref{sl2Mq})
will be used for   
$\ga=\begin{pmatrix} r & \star \\ s & \star \\ \end{pmatrix}$.
We need the following formula:

\begin{align}\label{garsmfunq}
&\ga^\circ\,\bigl(X_\la\,
(\ga^\circ)^{-1}(X_\mu)\,\bigr)
=\ q^{-rs\,\la^2/4-s\la\mu/2}\, X_{r\la+\mu}\,,
\end{align}
which results in

\begin{align}%\label{jonesf-1r}
&\ga^\circ \bigl(P_b^{(1)}\,(\ga^\circ)^{-1}(X\!-\!X^{-1})\bigr) 
= 
\ga^\circ \,\bigl(\sum_
{j=-b/2}^{b/2}X^{2j}\,(\ga^\circ)^{-1}
(X\!-\!X^{-1})
\,\bigr)\notag\\
&= \ssum_{j=-b/2}^{b/2}
\,q^{-(rj+1)sj}X^{2rj+1} - \ssum_{j=-b/2}^{b/2} 
\,q^{-(rj-1)sj}X^{2rj-1} \notag \\
&= \ssum_{j=-b/2}^{b/2} \,q^{-(rj+1)sj}X^{2rj+1} - 
\ssum_{j=-b/2}^{b/2} 
\,q^{-(rj+1)sj}X^{-2rj-1} \notag \\
&= \ssum_{j=-b/2}^{b/2} \,q^{-(rj+1)sj}\, (X - X^{-1})P_{2rj}^{(1)},
\notag
\end{align}
and then implies the relation

\begin{align} \label{DAHArecur}
&\Bigl(\hat{\ga}_{\ell}^{(1)}\,P_b^{(1)}\,
(\hat{\ga}_{\ell}^{(1)})^{-1}\Bigr)
\!\Downarrow\ =\, \bigl(
\x_1^{-1}\hat{\ga}^\circ_\ell \x_1 P_b^{(1)}\x_1^{-1} 
(\hat{\ga}^\circ_\ell)^{-1} \x_1 \bigr)\!\Downarrow\\
&= \ssum_{j=-b/2}^{b/2}\ q^{-(\rr_\ell j+1)\ss_\ell j} 
P_{2\rr_\ell j}^{(1)} \for \x_1 = X\!-\!X^{-1}. \notag
\end{align}

Thus we can switch from $k=1$ to the case of
$k=0$, where the free action of $PSL(2,\Z)$ can be used.
\smallskip

%\vskip 0.1cm
Using the latter relation we can proceed by induction
as follows:
\begin{align}
&\ \ \ \ \ 
J\!D^\#_{T(\{\rr_1,\ldots \rr_\ell\},\,\{\ss_1,\ldots \ss_\ell\})\,}
(b\,;\,q,t\mapsto q)\notag\\
 = \Bigl\{\hat{\ga}_{1}^{(1)}&\Bigr(
\cdots\Bigl(\hat{\ga}_{\ell-1}^{(1)}
\Bigl(\hat{\ga}_{\ell}^{(1)}P_b^{(1)}
(\hat{\ga}_{\ell}^{(1)})^{-1}\Bigr)
\!\Downarrow
(\hat{\ga}_{\ell-1}^{(1)})^{-1} \Bigr)\!\Downarrow\Bigr) \cdots\Bigr)
(\hat{\ga}_{1}^{(1)})^{-1}\Bigr\}^{(1)} \notag
\end{align}

\begin{align}
= \sum_{j=-b/2}^{b/2} q&^{-(\rr_\ell j+1)\,\ss_\ell j} 
\Bigl\{\hat{\ga}_{1}^{(1)}\Bigr(
\cdots\Bigl(\hat{\ga}_{\ell-1}^{(1)}
P_{2\rr_\ell j}^{(1)}
(\hat{\ga}_{\ell-1}^{(1)})^{-1} \Bigr)\!\Downarrow\Bigr) 
\cdots\Bigr)(\hat{\ga}_{1}^{(1)})^{-1}\Bigr\}^{(1)} \notag\\
 = \sum_{j=-b/2}^{b/2} q&^{-(\rr_\ell j+1)\,\ss_\ell j} 
J\!D^\#_{T(\{\rr_1,\ldots \rr_{\ell-1}\},\,\{\ss_1,\ldots 
\ss_{\ell-1}\})\,}(2\rr_\ell j\,;\,q,t\mapsto q) \notag\\
= \sum_{j=-b/2}^{b/2} q&^{-(\rr_\ell j+1)\,\ss_\ell\, j - 
\frac{1}{4}\aa_{\ell-1} \rr_{\ell-1}\, (2\rr_\ell j)\,
(2\rr_\ell j+2)}
\notag\\ 
&\ \ \ \ \ \ \ \ \ \ 
\times\j^{\#}_{T(\{\rr_1,\ldots \rr_{\ell-1}\},\,\{\ss_1,\ldots 
\ss_{\ell-1}\})\,}(2\rr_\ell j\,;\,q\mapsto q^{-1}) \notag
\end{align}
\begin{align}
 = \sum_{j=-b/2}^{b/2}\,\, q&^{-(\rr_\ell j+1)(\ss_\ell+ 
\aa_{\ell-1}\rr_\ell \rr_{\ell-1}) j } 
\j^{\#}_{T(\{\rr_1,\ldots \rr_{\ell-1}\},\,\{\ss_1,\ldots 
\ss_{\ell-1}\})\,}(2\rr_\ell j\,;\,q\mapsto q^{-1}) \notag
\end{align}
\vskip -0.3cm
\begin{align}
= \ssum_{j=-b/2}^{b/2}\,\, q&^{-(\rr_\ell j+1)
\aa_{\ell} j } \j^{\#}_{T(\{\rr_1,\ldots \rr_{\ell-1}\},
\,\{\ss_1,\ldots \ss_{\ell-1}\})\,}(2\rr_\ell j\,;
\,q\mapsto q^{-1})\notag
\end{align}
\vskip -0.3cm
\begin{align}
&= q^{-\frac{1}{4}\aa_\ell \rr_\ell\, b(b+2)} 
\j^{\#}_{T(\{\rr_1,\ldots \rr_\ell\},\,\{\ss_1,\ldots 
\ss_\ell\})\,}(b\,;\, q\mapsto q^{-1}). \notag
\end{align}

\comment{
{\it Remark.} As we see from the calculation parameters 
$\aa $ and $\ss $ appear in different types of recursive formulas. 
The parameter $\ss $ appears in relation (\ref{DAHArecur}), there 
is no common multiplier depending on the representation $b$. Due to 
topological reasons, it appears in (\ref{jones-formula}) and the 
correction of this dependence by multiplier 
$q^{-\frac{1}{4}\aa_\ell \rr_\ell b(b+2)}$ 
leads to shift from $\ss$ to $\aa$.\\
}

Finally,
\begin{align}
J\!D^\#_{\unknot }(b;\,q,t\mapsto q) = 
\frac{q^{(b+1)/2}-q^{-(b+1)/2}}{q^{1/2}-q^{-1/2}} = 
\j^{\#}_{\unknot}(b\,;\, q\mapsto q^{-1}),
\end{align}
which completes the proof.
\sq

\subsection{\bf Concluding remarks}
Concerning Part $(iii)$ of the Connection Conjecture, 
{\em adding colors\,} (arbitrary Young diagrams $\la$) 
to Conjecture 1.2 from \cite{ORS}
is a natural challenge.  This would provide 
a $t$\~extension of Theorem 1.2 from \cite{Ma}. 
Then our construction (which is for arbitrary colors), 
will be employed at full potential; we hope to extend
it from iterated torus knots (the setup of this paper)
to iterated {\em links}.

\smallskip
The switch to the curves $C_{\la}$ from \cite{Ma}, Section 1.3
supported on $C=\c_{\vec\rr,\vec\ss}$ and associated with Young 
diagrams $\la$ is natural here. This approach 
was initiated 
in \cite{DHS} (via the so-called {\em resolved conifold}).
Using (the germs of) such curves can be hopefully combined with
considering the {\em weight filtration\,} in 
the cohomology of the corresponding nested Hilbert schemes and 
Jacobian factors. 
\smallskip

On the other hand, DAHA uniformly manages arbitrary 
root systems and weights 
via the Macdonald polynomials; changing the curves  
$\,\c_{\vec\rr,\vec\ss}\,$ for incorporating the colors does 
not seem really necessary from this perspective.
Instead, the spaces Bun${}_G(\c_{\vec\rr,\vec\ss})$
of $G$\~bundles 
over $\,\c_{\vec\rr,\vec\ss}\,$ for the simple Lie group $G$
associated with the root system $R$ 
can be expected.  This is related to the interpretation of 
the Jacobian factors as Springer fibers due to \cite{La}; 
$\,\c_{\vec\rr,\vec\ss}\,$ would then become  
{\em spectral curves}. 

See also \cite{ChE}, Corollary 2.3 (and Theorem 4.10),
where this connection was established for 
torsion free sheaves over arbitrary (possibly singular)
coverings of rational or elliptic curves $E$ 
(assuming the stability of their direct images 
over $E$). The connection with {\em Baker functions\,} was 
established there (Proposition 2.4), which actually makes this 
construction a variant of the {\em Hitchin system\,} in the 
case of {\em factorizable Lie group schemes\,} over $E$ or $P^1$,
The latter are associated with (non-unitary) classical 
$r$\~matrices; see there and \cite{C101}, Section 1.7.
\smallskip
  
{\sf Khovanov-Rozansky homology.}
On the topological side, the {\em skein relations}, 
which play an important role in the proof 
of Theorem 1.2 from \cite{Ma}, are questionable 
upon the switch from the HOMFLYPT 
polynomials to their homology from \cite{KhR1,KhR2,Kh,Rou} 
(generalizing the Khovanov homology for the Jones polynomials). 
Adding colors (arbitrary Young diagrams, generally speaking dominant
weights) here 
is generally an unsettled problem, though the $QG$-based
{\em categorification theory\,} provides certain tools for this. 

Considering only {\em uncolored\,} Khovanov-Rozansky polynomials
is a significant simplification. In contrast to the Khovanov 
polynomials, there are only very few (uncolored) cases where 
the Khovanov-Rozansky polynomials and $K\!hR_{\hbox{\small stab}}$ 
are known, which are the Poincar\``e polynomials of the 
corresponding triply-graded homology. However
this theory works for {\em any links\,}; such a great universality 
is one of the reasons why this theory is so complicated. For us,
the Khovanov-Rozansky and Khovanov polynomials are the only way so 
far to interpret geometrically(-topologically) the
DAHA superpolynomials of {\em arbitrary\,} iterated torus knots, 
which are generally non-algebraic.
\smallskip

{\sf Betti numbers.}
The examples above confirm our conjecture that the DAHA 
superpolynomials
(among other expected applications) split the Euler numbers of 
the Jacobian factors of unibranch plane curve singularities into the 
corresponding 
Betti numbers. This is the simplest (and convincing) demonstration 
of the power of the refinement (adding "$t$" to the theory). This 
claim 
of our Connection Conjecture presumably follows from 
Conjecture 2 from \cite{ORS}; see the discussion of the 
{\em reduced\,} 
$\mathscr{P}_{\!\hbox{\tiny alg}}^{\hbox{\tiny min}}$ 
there. In the notation there (see the Overview and 
Section 4),
\begin{align}\label{overlinep}
&\mathscr{P}_{\!\hbox{\tiny alg}}=
\bigl(\frac{q_{st}}{a_{st}}\bigr)^{\mu}
\frac{1-q^2_{st}}{1+a_{st}^2 t_{st}}
\,\sum_{l,m\ge 0}q_{st}^{2l}\,a_{st}^{2m}\,
t_{st}^{m^2}\,\mathfrak{w}(\c_{\vec\rr,\vec\ss}^{[l\le l+m]}),
\end{align}
where $\mu$ is the Milnor number ($\mu=2\de$ in the unibranch
case) and $\mathfrak{w}$ is the {\em weight filtration\,} in the
compactly supported cohomology of the corresponding scheme 
(its existence was justified by Deligne).
 
The spectral sequence $E_1$ in terms of 
the weight filtration on the compactly supported cohomology
over $\Q$, which converges to 
the whole cohomology and degenerates at $E_2$ 
(for smooth not necessarily compact manifolds), can potentially
connect (\ref{overlinep}) with the Betti numbers.
Here one needs to know that the odd cohomology of the Jacobian 
factors vanish for unibranch plane curve singularities
(the van Straten -Warmt conjecture)
and that the corresponding {\em mixed Hodge filtration} 
is {\em pure\,} (a conjecture). We thank 
Vivek Shende for the explanations; see \cite{ORS}
(the Betti numbers are not explicitly discussed there). 

Not much is known concerning the Betti
numbers of the Jacobian factors, though there are 
quite a few known cases and natural conjectures for
the Puiseux exponents $(\rr,\ss)$ (i.e. for the 
torus knots) and for the series 
$(4,2u,v)$, where $4<2u<v$ for odd $u,v$. See \cite{Pi}.
\smallskip

{\sf Versal deformations.}
Last but not least, let us refer to \cite{GSh} for some 
relatively recent developments concerning the Hilbert
schemes of (locally) plane curves, including their role
in the theory of the Gopakumar-Vafa BPS invariants
(see \cite{PaT}) and various generalizations.

This paper is helpful to put our construction into perspective 
and to link it to the {\em versal deformations of plane curve 
singularities}, at least in the case of usual (non-nested) Hilbert 
schemes of curves. We do not touch this important direction
here. However let us mention that the theory
of {\em adjacent singularities\,} of the pairs $\{R,V\}$
for finite-dimensional irreducible representation of the Lie
group $G$ for $R$ is certainly related to the theory of 
the DAHA superpolynomials.

When $t_{st}=1$ and for the minimal possible degree of $a_{st}$
the sum in (\ref{overlinep}) essentially reduces to
the following summation
\begin{align}\label{PandT}
 \sum_{n\ge 0} q_{st}^{n+1-\de} e(\c^{[n]})=
\sum_{0\le i\le \de} n_{\c}{(i)}\bigl(\frac{q_{st}}
{(1-q_{st})^2}\bigr)
^{i+1-\de}, 
\ \c=\c_{\vec\rr,\vec\ss}\,, 
\end{align}
for the Euler numbers
of Hilbert schemes $\c^{[n]}$ of $n$ points on $\c$ of arithmetic
genus $\de$. 
See \cite{PaT,GSh}.
Importantly, $n_{\c}{(i)}\in \Z_+$ due to \cite{FGV} and \cite{Sh},
because these numbers
are the multiplicities of $\c$ in  
the closures of the strata in the space of its versal deformation, 
stratified with respect to the 
{\em geometric\,} genus of the deformation curves. 

In our notations, $t_{st}=1$ corresponds to $q=t$ and
therefore we can use our $\sqrt{q}$ instead of $q_{st}$ 
here upon $a=0$. 
Hopefully Conjecture 2 from \cite{ORS} and our construction
lead to a similar deformation interpretation of the right-hand side 
of (\ref{PandT}) upon its generalization to the nested 
Hilbert-schemes (which adds the parameter $a_{st}$) under the weight 
filtration (associated with $t_{st}$) and for arbitrary partitions
(generally, dominant weights). 
\medskip

{\bf Acknowledgements.}
Our special thanks go to Vivek Shende for clarifying discussions
on the topics related to \cite{ORS,GSh}. We thank Semen Artamonov 
for our using his unique software for calculating colored 
HOMFLYPT polynomials, Peter Samuelson for sending us his work 
before it was posted and him and Yuri Berest for discussions. 
The first author thanks Andras Szenes for the talks on the
algebraic-geometric aspects of our construction and the 
University of Geneva for the invitation and hospitality.

The paper was mainly written at RIMS. The
first author thanks Hiraku Nakajima and RIMS for
the invitation and hospitality and
the Simons Foundation (which made this visit possible);
the second author is grateful for the invitation to the 
school \& conference "Geometric Representation Theory" 
at RIMS (July 21- August 1, 2014) and generous RIMS' support.
I.D. acknowledges partial support from the RFBR grants 13-02-00478, 
14-02-31446-mol-a, NSh-1500.2014.2 and the common grant 
14-01-92691-Ind-a.

\vskip -2cm
%\medskip
\bibliographystyle{unsrt}

\end{document}